\documentclass[letterpaper,11pt]{article}
\usepackage[margin=1in]{geometry}  
\usepackage{setspace}  
\usepackage{natbib}
 \bibpunct[, ]{(}{)}{,}{a}{}{,}%

\usepackage{amsfonts}
\usepackage{amsthm}
\usepackage{mathrsfs}
\usepackage{amsmath}
\usepackage{amssymb}
\usepackage{latexsym}
\usepackage{indentfirst}
\usepackage{footmisc}
\usepackage{algorithm}
\usepackage{algorithmic}
\usepackage{float}
\usepackage{graphicx}
\usepackage{epstopdf}
\usepackage{longtable}
\usepackage{enumerate}
\usepackage{enumitem}
\usepackage{authblk}
\usepackage{bbm}
\usepackage{booktabs}
\usepackage{url}

\usepackage{multirow}
\usepackage{bm}
\usepackage{array}
\usepackage{subfigure}
\usepackage{caption}
\usepackage{color}

\def\x{\mathbf{x}}

\def\x{\mathbf{x}}

\def\c{\mathbf{c}}
\def\bkappa{\boldsymbol{\kappa}}

\def\Acal{\mathcal{A}}
\def\Bcal{\mathcal{B}}

\def\Ccal{\mathcal{C}}
\def\Ecal{\mathcal{E}}

\def\Ncal{\mathcal{N}}
\def\Xcal{\mathcal{X}}
\def\Scal{\mathcal{S}}
\def\Gcal{\mathcal{G}}

\def\Zcal{\mathcal{Z}}

\def\Ucal{\mathcal{U}}
\def\Vcal{\mathcal{V}}

\def\Xcal{\mathcal{X}}
\def\Kcal{\mathcal{K}}
\def\Pcal{\mathcal{P}}
\def\Hcal{\mathcal{H}}
\def\CSO{\mathcal{CSO}}
\def\CS{\mathcal{S}}
\def\CP{\mathcal{P}}
\def\Tcal{\mathcal{T}}

\def\Pb{ \mathbb{P} }
\def\Eb{ \mathbb{E} }
\def\Rb{ \mathbb{R} }
 
\def\1{\mathbb{I}}

\newtheorem{assumption}{ASSUMPTION}
\newtheorem{lemma}{LEMMA}

\newtheorem{theorem}{THEOREM}
\newtheorem{proposition}{PROPOSITION}
\newtheorem{corollary}{COROLLARY}

\newtheorem{definition}{DEFINITION}

\allowdisplaybreaks

\title{Understanding the Variance Dichotomy in Continuous Simulation Optimization: A Minimax Lower Bound Perspective}

\author[a]{Jianzhong Du}
\author[b]{L. Jeff Hong}
\affil[a]{Faculty of Business for Science and Technology, School of Management, University of Science and Technology of China, dujianzhong@ustc.edu.cn}
\affil[b]{Department of Industrial and Systems Engineering,
University of Minnesota, lhong@umn.edu}

\date{\vspace{-8ex}}

\begin{document}

\maketitle

\begin{abstract}
This paper studies the variance dichotomy in continuous simulation optimization (CSO). Existing literature shows a sharp contrast between deterministic CSO and stochastic CSO, with convergence rates in stochastic settings appearing insensitive to the magnitude of the noise variance. However, this asymptotic view does not fully explain the behavior of CSO under finite simulation budgets, especially in low-noise settings.
To address this gap, this work develops a minimax lower-bound analysis and shows that the complexity is decided by the maximum of a variance-dependent term and a variance-independent term.  When the simulation budget is not very large and the noise variance is low, the variance-independent term dominates, implying that low-noise stochastic CSO has essentially the same complexity as deterministic CSO. As the budget increases, the variance-dependent term becomes dominant, and the convergence behavior of stochastic CSO transitions to a slower regime determined jointly by the noise variance and the simulation budget.
\end{abstract}




\section{Introduction}\label{sec:intro}

Continuous simulation optimization (CSO) concerns a class of continuous stochastic optimization problems in which the objective function can only be evaluated through simulation experiments. Such problems arise naturally in many complex systems where analytical expressions of the objective are unavailable and performance must be assessed via simulation. CSO has been extensively studied in the simulation literature, see the review articles of \cite{fu2002optimization,hong2009brief}, and \cite{fan2025review}, and it plays an important role in a wide range of applications, including supply chain management \citep{wang2023large}, transportation \citep{jia2020simulation}, and healthcare \citep{gong2022managing}. CSO is also closely connected to black-box optimization or zeroth-order optimization \citep{liu2020primer,larson2019derivative,grill2015black}, which has recently attracted significant attention in machine learning, particularly in hyperparameter tuning \citep{Jiang2024OpenBox}. In addition, it is related to continuum-armed bandit models \citep{agrawal1995continuum}, although the underlying objectives in these settings differ.

Over the past decades, a variety of algorithms have been proposed to solve CSO problems. These include gradient-based methods \citep{kiefer1952stochastic,spall2000adaptive}, random search methods \citep{hu2007model,andradottir2014review,kiatsupaibul2018single}, and model-based approaches \citep{barton2006metamodel,hong2021surrogate,chen2023pseudo}. To evaluate and compare these methods, a standard criterion is the convergence rate, which characterizes how fast the optimality gap, also known as the simple regret \citep{bubeck2018}, decreases to zero as the number of simulation observations, or budget, increases. Beyond its asymptotic interpretation, the convergence rate also provides a principled framework for assessing algorithmic performance under finite simulation budgets, which is often the regime of practical interest. 

However, standard convergence rate analyses often overlook the role of the variance of simulation observations in determining the performance of CSO algorithms, as it is typically absorbed into constant terms. In contrast, variance is well known to have a substantial impact on the finite-time performance of these algorithms, yet there remains limited understanding of how it influences algorithmic behavior in a systematic way. Early work in the CSO literature has attempted to mitigate the impact of variance by incorporating ranking-and-selection (R\&S) procedures into the optimization process. For example, \cite{pichitlamken2006sequential} incorporated R\&S into neighborhood comparisons, and \cite{hong2007selecting} applied R\&S throughout the entire optimization procedure. However, these approaches have generally not led to significant performance improvements. While it is often conjectured that allocating computational effort to explore more solutions is more effective than reducing estimation errors caused by simulation variance, the precise reasons behind this phenomenon remain not well understood.

In this paper, we examine an apparent dichotomy in standard convergence rate results that is driven by the variance of simulation noise, and we argue that this phenomenon provides key insight into the role of variance in CSO performance. We assume that the simulation noise follows subgaussian distributions with a common upper bound on the variance proxy, denoted by $\sigma^2$. A striking feature of existing results is that, as long as $\sigma^2 > 0$, the convergence rates of CSO algorithms are typically unaffected by the magnitude of $\sigma^2$, with its influence appearing only in constant terms. However, this behavior changes fundamentally when $\sigma^2 = 0$, where markedly different convergence rates can arise. For example, when the objective function is strongly concave, the optimal convergence rate for stochastic CSO is polynomial, of order $O(n^{-1/2})$ \citep{hu2024convergence,shamir2013complexity}. In contrast, in the absence of simulation noise, the convergence rate can become exponential, of order $O(e^{-cn})$ for some $c>0$ \citep{munos2014}. When the objective function is less regular and only assumed to be Lipschitz continuous, the best achievable rate for stochastic CSO further deteriorates to $\tilde{O}(n^{-1/(d+2)})$ \citep{yakowitz2000global}, up to logarithmic factors. In the corresponding deterministic setting, the convergence rate improves to $O(n^{-1/d})$, although it still deteriorates rapidly with the dimension $d$.

These observations point to what we term a {\it variance dichotomy} in CSO. Specifically, under a given smoothness assumption on the objective function, there are two distinct regimes of convergence: one corresponding to deterministic CSO with $\sigma^2 = 0$, and the other to stochastic CSO with $\sigma^2 > 0$. Moreover, within the stochastic regime, the convergence rate is largely insensitive to the magnitude of the noise variance, in the sense that problems with low noise levels and those with high noise levels share the same rate up to constant factors.

\begin{figure}[tbh]
	\centering
	\includegraphics[width=0.5\textwidth]{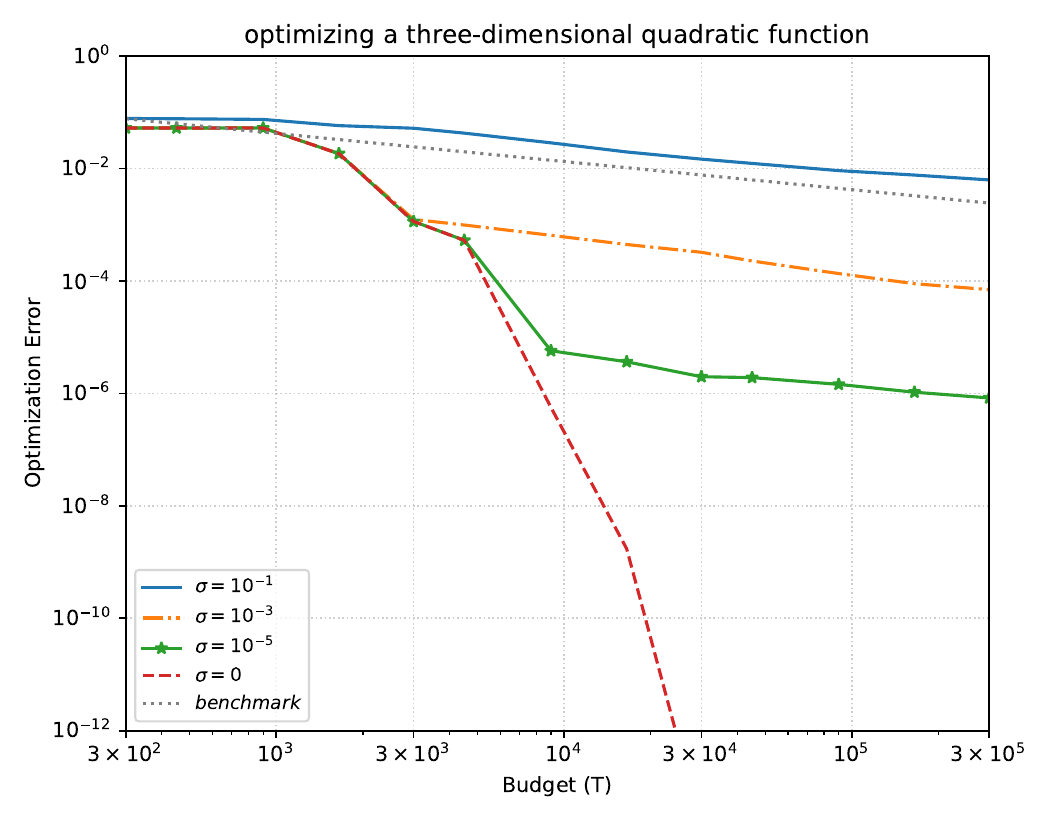}
	\caption{The decrease of optimization error when the observation's variance is different.}\label{fig:intro}
\end{figure}

To better understand this dichotomy, consider the problem of maximizing the three-dimensional strongly concave function $y(\x) = -\| \x - \c\|_2^2$ over $\x \in [0,1]^3$ where $\c$ is the vector with all coordinates equal to $e^{-1}$. We apply the StroquOOL algorithm \citep{bartlett2019simple} to this problem under normally distributed observations with variance $\sigma^2$. Figure~\ref{fig:intro} plots the resulting optimization errors for different values of $\sigma^2$, along with a benchmark line proportional to $n^{-1/2}$. Two key observations emerge. First, when the budget is sufficiently large, for example exceeding $10^4$, the variance dichotomy becomes evident: the optimization errors in the stochastic cases ($\sigma^2 > 0$) decrease at similar rates across different noise levels, yet these rates are significantly slower than that of the deterministic case ($\sigma^2 = 0$). Second, in the small-budget regime, the performance of the stochastic cases closely tracks that of the deterministic case up to a variance-dependent switching point. This switching occurs earlier when the variance is high and later when the variance is low.

The switching behavior is critical for understanding the variance dichotomy as well as the role of variance in CSO algorithms. First, although asymptotic convergence rates exhibit a dichotomy as the budget goes to infinity, the switching point may vary continuously with $\sigma$, indicating that variance affects the performance of CSO algorithms in a continuous manner. Second, it suggests that when the simulation budget is limited and the variance is not too large, variance may have only a minor impact on algorithm performance. This observation helps explain why incorporating R\&S procedures to mitigate variance may not lead to performance gains, as observed in \cite{pichitlamken2006sequential} and \cite{hong2007selecting}.

In this paper, we adopt a minimax lower bound approach to understand the switching behavior. This perspective allows us to derive results that are independent of specific CSO algorithms and instead characterize the best achievable convergence rates. We develop two types of lower bounds, one variance-dependent and the other variance-independent, corresponding to the two poles of the variance dichotomy under different smoothness conditions. The variance-dependent lower bound arises from stochastic noise and therefore vanishes in the deterministic case, while the variance-independent lower bound does not depend on $\sigma^2$ and remains the same for both stochastic and deterministic settings. Consequently, the effective lower bound on the optimization error is given by the maximum of these two bounds, which naturally captures the switching between stochastic and deterministic regimes.

When the budget is small, the variance-independent lower bound dominates, and the behavior of the lower bound coincides with that of the deterministic case where $\sigma^2 = 0$. As the budget increases, the variance-dependent lower bound eventually dominates, and the convergence rate aligns with those established for stochastic CSO in the asymptotic regime. In contrast, for deterministic CSO, the effective lower bound always coincides with the variance-independent one, implying that deterministic CSO remains intrinsically easier across all budget regimes. In this sense, the deterministic case can also be understood within our framework as a limiting scenario in which no switching occurs for any finite budget.

\subsection{Literature Review}

There is a substantial literature on lower bounds for CSO under various structural assumptions. For stochastic CSO with strongly concave objectives, \cite{shamir2013complexity} and \cite{akhavan2020exploiting} established a lower bound of order $\Omega(n^{-1/2})$, and showed that faster rates may be achievable under higher-order smoothness assumptions. \cite{wang2018optimization} studied stochastic CSO problems whose objective functions possess geometric structures similar to a reference function with known smoothness properties, and derived lower bounds in the asymptotic regime where the simulation budget $n$ is sufficiently large. For stochastic CSO with Lipschitz continuous and potentially highly irregular objective functions, \cite{singh2021continuum} obtained a lower bound of order $\Omega(n^{-1/(d+2)})$.

The near-optimality dimension \citep{munos2014,grill2015black,bubeck2018,kleinberg2008multi} provides a more refined way to characterize problem difficulty in stochastic CSO. In particular, \cite{locatelli2018adaptivity} established a lower bound of order $\Omega(n^{-1/(D+2)})$ for stochastic CSO problems whose near-optimality dimension is less than $D$, where $D$ typically ranges between $0$ and $d$. In the extreme case of $D=0$, this bound coincides with that for strongly concave stochastic CSO; where when $D=d$, it approaches the lower bound for Lipschitz stochastic CSO. We will introduce the notion of near-optimality dimension in more detail in Section~\ref{sec:form}.

Lower bounds have also been derived for deterministic CSO. For example, \cite{singh2021continuum} showed that for deterministic CSO with Lipschitz continuous objectives, the lower bound is $\Omega(n^{-1/d})$. Using a reproducing kernel Hilbert space characterization of smoothness, \cite{bull2011convergence} derived a lower bound of order $\Omega(n^{-\nu/d})$, where $\nu$ is a smoothness parameter.

These lower bounds are commonly used to establish the rate optimality of CSO algorithms. For instance, the bound $\Omega(n^{-1/2})$ implies that the convergence rate $O(n^{-1/2})$ achieved by stochastic approximation methods \citep{kiefer1952stochastic,hu2024convergence} is optimal for strongly concave stochastic CSO. Similarly, the bound $\Omega(n^{-1/(d+2)})$ indicates that grid search methods \citep{yakowitz2000global} are nearly rate-optimal for Lipschitz stochastic CSO. The lower bound $\Omega(n^{-1/(D+2)})$ further supports the near-optimality of several tree-based algorithms \citep{munos2014,grill2015black,shang2019general,bartlett2019simple} for stochastic CSO problems with near-optimality dimension less than $D$.

Most existing lower-bound analyses and algorithmic results treat stochastic and deterministic CSO as separate problems. To our knowledge, lower-bound analyses that simultaneously capture both poles of the variance dichotomy within a unified framework are limited. An exception is \cite{singh2021continuum}, which derived lower bounds for both stochastic and deterministic CSO by characterizing function smoothness using Besov spaces generalized from the class of Lipschitz functions. Their results imply that pure exploration methods, such as grid search, are optimal in this setting, but at the cost of convergence rates that suffer severely from the curse of dimensionality. While the class of Lipschitz functions contains challenging optimization instances, it also includes relatively easy ones. In contrast, the smoothness framework adopted in this work is capable of distinguishing easy-to-optimize functions from difficult ones, leading to lower bounds that adapt to the underlying smoothness and yield faster rates for simpler problem classes.

From an algorithmic perspective, \cite{bartlett2019simple} proposed an algorithm that achieves fast convergence rates for both deterministic and stochastic CSO. They observed that when the simulation budget is moderate, the convergence rate for low-noise stochastic CSO matches that of deterministic CSO, a phenomenon illustrated in Figure~\ref{fig:intro}. However, it remains unclear whether these rates are optimal. In this work, we derive lower bounds that characterize the variance dichotomy and nearly match the convergence rates obtained in \cite{bartlett2019simple}. This agreement establishes the tightness of our lower bounds and confirms the near-optimality of the their algorithmic rates.

The rest of this paper is organized as follows. Section \ref{sec:form} formulates the CSO problem and introduces the general form of the minimax lower bounds. Section \ref{sec:lem} provides three preliminary lemmas for the minimax analysis. Section \ref{sec:var_dep} establishes the variance-dependent lower bound, and Sections \ref{sec:var_ind_uneq} and \ref{sec:var_ind_eq} establish the variance-independent lower bounds. Section \ref{sec:lb_overall} presents the overall lower bounds for the variance dichotomy. Section \ref{sec:num} presents illustrative numerical experiments. Section \ref{sec:con} concludes the paper, and the Electronic Companion contains all technical proofs.

\section{Problem Formulation}\label{sec:form}

We consider the following CSO problem:
\begin{align}\label{prob:cso}
	\max_{\x \in \Xcal} y(\x) := \Eb [ Y(\x,\omega) ],
\end{align}
where $\x=(x_1,\dots,x_d)$ denotes the vector of decision variables, $\Xcal$ is a $d$-dimensional compact set representing the feasible solution space, and $y(\x)$ is the (unknown) objective function with $y(\x) = \Eb [ Y(\x,\omega) ]$. Here, $Y(\x,\omega)$ is a simulation observation of the objective function at $\x$, and $\omega$ is a random vector representing the stochastic uncertainty of the underlying system during simulation. Without loss of generality, we take $\Xcal = [0,1]^d$ for the theoretical analysis. Throughout, we assume that the objective function can be evaluated only through simulation.

We make the following assumption on Problem (\ref{prob:cso}). The first is on the uniqueness of the optimal solution, while the second is on the noise term $\varepsilon(\x)=Y(\x)-y(\x)$.

\begin{assumption}\label{ass:unique}
	The optimal solution to $\max_{\x \in \Xcal} y(\x)$ is unique, which is denoted by $\x^*$.
\end{assumption}

\begin{assumption} \label{ass:subgauss}
	Given $\x$, the noise term $\varepsilon(\x)$ is independent of everything else and $\sigma$-subgaussian with $\Eb \left[\exp (\lambda\varepsilon(\x))\right] \le \exp ( \lambda^2 \sigma^2 /2 )$ for $\lambda \in \Rb$ where $\sigma^2 \ge 0$ is the variance constant.
\end{assumption}
Both assumptions are quite standard in the SO literature. It is worthwhile pointing out that the $\sigma$-subgaussian assumption of Assumption \ref{ass:subgauss} permits the distribution to be normal, which is also common in the literature, or any bounded distribution. In the special case where $\sigma^2 = 0$, the simulation observation is deterministic.

In rest of this sections, we provide the shape constraint on the objective function, discuss the general setting of CSO algorithms, and introduce the form of lower bounds studied in this work.

\subsection{Shape Constraint}

Let $\| \x \|_{\infty} =  \max_{l=1,2,\dots,d} | x_l | $ denote the maximum norm of vector $\x$. In convergence analysis, this norm is equivalent to the Euclidean norm $\| \x \|_{2} = \left( \sum_{l=1}^d  x_l ^2 \right)^{1/2}$ in finite-dimensional spaces, in the sense that there exist positive constants $c$ and $C$ such that $c \| \x \|_{\infty} \le \| \x \|_{2} \le C \| \x \|_{\infty}$ for all $\x \in \Xcal$ \citep{buhler2018functional}. The use of the maximum norm is primarily for analytical convenience and simplifies the presentation in this paper. We use the following assumption based on the maximum norm to characterize the shape of objective function.

\begin{assumption}\label{ass:smo}
	There exist smoothness parameters $\alpha$ and $\beta$ with $0 < \alpha \le \beta < \infty$ such that 
	\begin{align}\label{ineq:localsmo_ass2}
		\tilde{M} \| \x^* - \x \|_{\infty}^\beta \le \left| y(\x^*) - y(\x) \right| \le M \| \x^* - \x \|_{\infty}^\alpha, \ \forall \x \in \Xcal
	\end{align}
	where $M$ and $\tilde{M}$ are strictly positive constants.
\end{assumption}

Assumption \ref{ass:smo} uses two polynomial envelopes of orders $\alpha$ and $\beta$ to upper and lower bound the shape of the objective function. The upper bound in \eqref{ineq:localsmo_ass2} ensures that $y(\x)$ approaches $y(\x^*)$ as $\x$ approaches $\x^*$, while the lower bound prevents the function from being excessively flat near $\x^*$. Throughout the analysis, the primary quantities of interest are the parameters $\alpha$ and $\beta$, which capture the geometry of the objective function around $\x^*$. The constants $M$ and $\tilde{M}$ can be chosen flexibly, provided that they do not depend on the distance $\|\x^*-\x\|_{\infty}$ or the simulation budget $n$.

\begin{figure}[t]
	\centering
	\includegraphics[width=0.99\textwidth]{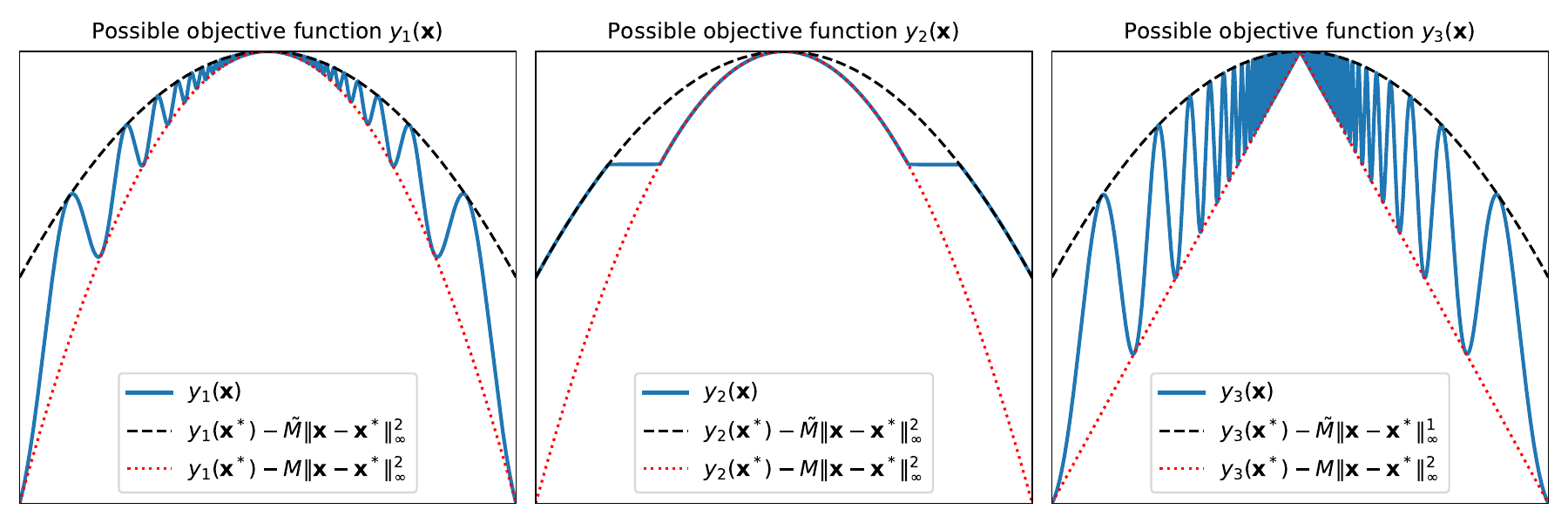}
	\caption{Plot of several possible functions under Assumption \ref{ass:smo}.}\label{fig:fun_ex}
\end{figure}

Assumption \ref{ass:smo} is general enough to incorporate many objective functions. For example, if $y(\x)$ is strongly concave \citep{boyd2004convex} such that its Hessian matrix satisfies $\tilde{M} I \preceq \nabla^2 y(\x) / 2 \preceq M I$ and the optimal solution $\x^*$ is stationary, then Assumption \ref{ass:smo} is satisfied naturally with $\alpha = \beta = 2$. If $y(\x)$ is strongly concave but $\x^*$ is non-stationary, then Assumption \ref{ass:smo} may be satisfied with $\alpha = \beta = 1$. Many non-concave functions also fit Assumption \ref{ass:smo}. Figure \ref{fig:fun_ex} provides three one-dimensional examples. The first panel shows a function whose value changes between the upper and lower polynomial envelopes of order two. Although it is non-concave, it satisfies Assumption \ref{ass:smo} with $\alpha=\beta=2$. The second panel plots a function that remains constant in some sub-regions. Still, it satisfies Assumption \ref{ass:smo} with $\alpha=\beta=2$. The value of $\alpha$ and $\beta$ can be different to characterize bumpier functions. The third panel provides a function whose value jumps between the two polynomial envelopes of order one and two, respectively.  This function satisfies Assumption \ref{ass:smo} with $\alpha=1$ and $\beta=2$.

Parameter $\alpha$ is related to the local smoothness at the optimal solution because the upper bound in \eqref{ineq:localsmo_ass2} 
basically says that the objective function is smooth locally at $\x^*$. This requirement is milder than the global smoothness assumption such as the Lipschitz condition $\left| y(\x) - y(\x') \right| \le M \| \x - \x' \|_{\infty}^1$, $\forall \x,\x' \in \Xcal$. Meanwhile, parameter $\beta$ is related to how fast the volume of the $\epsilon$-optimal region $\{ \x \in \Xcal: y(\x^*) - y(\x) \le \epsilon \}$ decreases as $\epsilon \to 0$. Due to the lower bound in \eqref{ineq:localsmo_ass2}, the $\epsilon$-optimal region must be a subset of the ball $\{ \x \in \Xcal: \| \x - \x^* \|_{\infty} \le (\epsilon / \tilde{M})^{1/\beta} \}$, which implies as $\epsilon \to 0$, the region's volume decreases to zero with order $d/\beta$.

Similar assumptions to Assumption \ref{ass:smo} have been adopted in the literature. \cite{locatelli2018adaptivity} also used two parameters $\alpha_{LC}$ and $\beta_{LC}$ to characterize the shape of the objective function. Under Assumption \ref{ass:smo}, their $\alpha_{LC}$ and $\beta_{LC}$ are equivalent to $\alpha$ and $d/\beta$ of this work. \cite{grill2015black} used the near-optimality dimension to characterize the function shape and smoothness. In Section \ref{sec:near_opt}, we show that the near-optimality dimension is less than $d(1/\alpha-1/\beta)$ under Assumption \ref{ass:smo}.

\subsection{CSO Algorithms}

CSO problems are typically solved using algorithms that adaptively determine which solutions to simulate and, after exhausting the simulation budget, output an estimated optimal solution $\hat{\x}_n^*$. Throughout the paper, the budget refers to the total number of simulation observations and is denoted by $n$. The minimax analysis in this paper is conducted over a broad class of CSO algorithms and establishes lower bounds that hold uniformly for all algorithms in this class. These algorithms are assumed to satisfy the following general assumption. 

\begin{assumption}\label{ass:adaptive_algo}
	The CSO algorithm with total budget $n$ satisfies the following properties.
	\begin{enumerate}
		\item For $t=0$, the first simulated solution $\x_{1}$ is either deterministic or a measurable function of an independent random vector $U_{1}$.  
		
		\item For $t=1,2,\dots,n-1$, the $(t+1)$-th simulated solution $\x_{t+1}$ is a measurable function of the historical simulated solutions and their simulation observations $\x_1,Y(\x_1),\dots,\x_t,Y(\x_t)$, together with an independent random vector $U_{t+1}$.
		
		\item The estimated optimal solution $\hat{\x}_n^*$ is a measurable function of the historical simulated solutions and their simulation observations $\x_1,Y(\x_1),\dots,\x_n,Y(\x_n)$ and an independent random vector $U_{n}^*$.
		
		\item $\hat{\x}_n^* \in \{\x_1,\dots,\x_n\}$.
	\end{enumerate}
\end{assumption}

Assumption~\ref{ass:adaptive_algo} characterizes how a CSO algorithm selects simulated solutions in the first two items and produces the estimated optimal solution in the third item. Similar assumptions have been used explicitly by \citet{wang2018optimization} and implicitly by \citet{locatelli2018adaptivity}. This framework accommodates a broad class of adaptive algorithms, including also static algorithms in which the simulated solutions are pre-specified.

The independent random vectors $\{U_{1},\dots,U_{n},U_{n}^*\}$ are introduced to capture the randomness employed by certain algorithms, such as random search methods for selecting simulated solutions \citep{wang2023gaussian,kiatsupaibul2018single} and randomized recommendation methods for choosing the estimated optimal solution \citep{bubeck2018}. In contrast, other algorithms--such as Bayesian optimization \citep{bull2011convergence} and tree search methods \citep{munos2014}--employ deterministic sampling and recommendation strategies. In these cases, the simulated solutions and $\hat{\x}_n^*$ are invariant with respect to $\{U_{1},\dots,U_{n},U_{n}^*\}$, yet still satisfy the first three items of Assumption~\ref{ass:adaptive_algo}. Hence, the assumption is sufficiently general to encompass many popular CSO algorithms.

Finally, the fourth item of Assumption~\ref{ass:adaptive_algo} requires that the estimated optimal solution $\hat{\x}_n^*$ be one of the simulated solutions. This condition is imposed purely for technical convenience. For algorithms that may recommend a solution not previously simulated, one may construct an equivalent algorithm with budget $n+1$ whose first $n$ simulated solutions coincide with those of the original algorithm and whose $(n+1)$-th simulated solution equals the original recommendation. The performance of the original algorithm can then be analyzed through this modified procedure.

\subsection{Minimax Lower Bound}\label{subsec:minimax}

The estimated optimal solution $\hat{\x}_n^*$ is generally random and rarely coincides exactly with the true optimal solution $\x^*$. Its distribution depends on both the algorithm employed and the underlying CSO problem. Let $\CSO$ denote a CSO problem satisfying Assumptions~\ref{ass:unique}--\ref{ass:smo}, and let \texttt{alg} denote a CSO algorithm satisfying Assumption~\ref{ass:adaptive_algo}. The optimization error of Algorithm \texttt{alg} on Problem $\CSO$ is defined as
\begin{align*}
\Eb_{\CSO}^{\texttt{alg}} \bigl[ y(\x^*) - y(\hat{\x}_n^*) \bigr],
\end{align*}
where the expectation is taken with respect to the distribution of $\hat{\x}_n^*$. An efficient CSO algorithm should yield a small optimization error. For notational convenience, we occasionally suppress the superscript $\texttt{alg}$ or the subscript $\CSO$ when the context is clear.

The difficulty of CSO problems varies substantially, and it is therefore standard to focus on a class of problems satisfying certain regularity conditions. Let $\Scal$ denote the set of CSO problems of interest, and let $\Acal$ denote the set of admissible CSO algorithms. A function $B(n,\sigma^2)$ is called a minimax lower bound on the optimization error if
\begin{align*}
\inf_{\texttt{alg} \in \Acal} \sup_{\CSO \in \Scal} \Eb_{\CSO}^{\texttt{alg}} \bigl[ y(\x^*) - y(\hat{\x}_n^*) \bigr] \ge B(n,\sigma^2).
\end{align*}
That is, for any algorithm in $\Acal$, there exists at least one problem in $\Scal$ for which the optimization error is no smaller than $B(n,\sigma^2)$. Typically, $B(n,\sigma^2)$ decreases to zero as $n \to \infty$. In this work, we explicitly emphasize the dependence on $\sigma^2$ to highlight the variance dichotomy.

Minimax lower bounds are commonly established by constructing a subset of CSO problems, $\Tcal \subset \Scal$, and showing that, for any algorithm \texttt{alg}$ \in \Acal$,
\begin{align*}
\sup_{\CSO \in \Tcal} \Eb_{\CSO}^{\texttt{alg}} \bigl[ y(\x^*) - y(\hat{\x}_n^*) \bigr] \ge B(n,\sigma^2).
\end{align*}
In other words, the optimization error must exceed the lower bound for at least one problem in $\Tcal$. Although $\Scal$ may contain infinitely many problems, the subset $\Tcal$ used in the analysis is often finite and selected for analytical tractability. For instance, under Assumptions~\ref{ass:unique}--\ref{ass:smo}, one may restrict attention to CSO problems with normally distributed simulation noise rather than merely subgaussian noise, or fix the constants $M$ and $\tilde{M}$ in Assumption~\ref{ass:smo} to specific values, provided they remain independent of the budget $n$ and the decision variable $\x$. This flexibility naturally raises the question of how tight the resulting lower bounds are.

A minimax lower bound $B(n,\sigma^2)$ is said to be tight if there exists an algorithm whose optimization error matches the bound up to constant or logarithmic factors. Because the exact optimization error of a given algorithm is typically intractable to derive, asymptotic analysis is commonly used to characterize its order in terms of the budget $n$ (and $\sigma^2$ in this work). We therefore adopt the standard Big-$O$, Big-$\Omega$, and $\tilde{O}$ notations.
\begin{definition}
	Let $a(n,\sigma^2)$ and $b(n,\sigma^2)$ be positive-valued functions. Then:
	\begin{itemize}
		\item $a(n,\sigma^2) = O(b(n,\sigma^2))$ if there exists a constant $c_{o,1} > 0$ such that 
		$a(n,\sigma^2) \le c_{o,1} b(n,\sigma^2)$.
		
		\item $a(n,\sigma^2) = \Omega(b(n,\sigma^2))$ if there exists a constant $c_{o,2} > 0$ such that 
		$a(n,\sigma^2) \ge c_{o,2} b(n,\sigma^2)$.
		
		\item $a(n,\sigma^2) = \tilde{O}(b(n,\sigma^2))$ if there exist constants $c_{o,3}, c_{o,4} > 0$ such that 
		$a(n,\sigma^2) \le c_{o,3} (\log n)^{c_{o,4}} b(n,\sigma^2)$.
	\end{itemize}
\end{definition}

Suppose the minimax lower bound $B(n,\sigma^2)$ satisfies $B(n,\sigma^2) = \Omega(b(n,\sigma^2))$, up to multiplicative constants. If there exists an algorithm $\texttt{alg}^*$ such that
\begin{align*}
\sup_{\CSO \in \Scal} \Eb_{\CSO}^{\texttt{alg}^*} \bigl[ y(\x^*) - y(\hat{\x}_n^*) \bigr]
= O(b(n,\sigma^2)) \quad \text{or} \quad \tilde{O}(b(n,\sigma^2)),
\end{align*}
then the lower bound is said to be tight or near-tight, and the algorithm $\texttt{alg}^*$ is optimal or near-optimal, respectively.

\section{Three Preliminary Lemmas}\label{sec:lem}

To elucidate the variance dichotomy through a lower-bound analysis, we establish two types of minimax lower bounds: one that explicitly depends on the noise variance and another that is variance-independent. The overarching strategy for deriving these bounds is similar in both cases. First, we construct a benchmark CSO problem that may or may not satisfy Assumptions~\ref{ass:unique} and~\ref{ass:smo}, but for which certain properties of the estimated optimal solution $\hat{\x}_n^*$ can be quantified regardless of the algorithm employed (provided the algorithm satisfies Assumption~\ref{ass:adaptive_algo}). We then construct a sequence of target CSO problems that satisfy Assumptions~\ref{ass:unique}--\ref{ass:smo}, and show that, given the behavior of the benchmark problem, at least one of the target problems must exhibit an unfavorable optimality gap that defines the minimax lower bound. A key step in this analysis is to transport probabilistic behaviors from the benchmark problem to the target problems. The three lemmas presented in this section facilitate this transportation.

We define two CSO problems as follows:
\begin{itemize}
	\item {Problem 1:} the objective function is $y_{1}(\x)$, and $\varepsilon(\x) \sim N(0,\sigma^2)$;
	
	\item {Problem 2:} the objective function is $y_{2}(\x)$, and $\varepsilon(\x) \sim N(0,\sigma^2)$.
\end{itemize}
No assumptions are imposed on $y_{1}(\x)$ and $y_{2}(\x)$ with respect to Assumptions~\ref{ass:unique} and~\ref{ass:smo}, ensuring that the lemmas apply to a broad class of problems. In the subsequent analysis, Problem~1 will serve as the benchmark problem, which may violate these assumptions, whereas Problem~2 will represent the target problem that satisfies them. Throughout, we fix the simulation noise $\varepsilon(\x)=Y(\x)-y(\x)$ to follow the distribution $N(0,\sigma^2)$ for every solution $\x$ and every CSO problem considered. This Gaussian setting satisfies Assumption~\ref{ass:subgauss}, simplifies the exposition of the variance dichotomy, and admits closed-form expressions for likelihoods and divergences.

Since all lemmas hold for any algorithm satisfying Assumption~\ref{ass:adaptive_algo}, we omit the superscript $\texttt{alg}$ in expectation and probability operators for notational simplicity.

\subsection{A Lemma on Probability Transportation}

This subsection presents a lemma that is primarily applicable to stochastic CSO problems with $\sigma^2 > 0$. Note that the distribution of the estimated optimal solution $\hat{\x}_n^*$ generally differs across CSO problem instances. Suppose that $\hat{\x}_n^*$ is the estimated optimal solution obtained by applying a given algorithm to Problem~1. We seek to characterize the likelihood that the same algorithm, when applied to  Problem~2, yields the same solution. The following lemma addresses this question. Its detailed proof is provided in Section \ref{subsec:proof_lemlb}.

\begin{lemma}\label{lem:lb}
	Suppose that $\Zcal$ is a subregion of $\Xcal$ and that a CSO algorithm satisfying Assumption~\ref{ass:adaptive_algo} is used to solve the CSO problems. Let $\x_1,Y(\x_1),\dots,\x_n,Y(\x_n)$ denote the trajectory of the algorithm when solving CSO Problem~1, and let $\hat{\x}_n^*$ be the resulting estimated optimal solution. Then, for any such $\Zcal$ and any such algorithm,
	\begin{align}\label{ineq:probtrans}
		\Pb_{2}\bigl(\hat{\x}_n^* \in \Zcal\bigr)
		\le \sqrt{\frac{\Eb_{1}(L_n)}{2}}
		+ \Pb_{1}\bigl(\hat{\x}_n^* \in \Zcal\bigr),
	\end{align}
    where $L_{n}  = \sum_{t=1}^n \left[ y_{1}(\x_t) - y_{2}(\x_t)\right]^2  /(2 \sigma^2)$, and $\Pb_i(\cdot)$ and $\Eb_i(\cdot)$ denote probability and expectation, respectively, under CSO Problem~$i$, for $i=1,2$, when using the algorithm.
\end{lemma}

Note that $\hat{\x}_n^*$ can be viewed as a measurable function of the algorithmic trajectory $\x_1,Y(\x_1),\dots,$ $\x_n,Y(\x_n)$. The proof of Lemma \ref{lem:lb} is based on a change-of-measure argument that relates the probability of this trajectory under CSO Problem 1 to that under CSO Problem 2, with $L_n$ serving as the log-likelihood ratio between the two induced distributions. It is also interesting to note that $\Eb_{1}(L_n)$ is also the Kullback-Leibler (KL) divergence between the two distributions. Then, if two problems are close to each other, $\Eb_{1}(L_n)$ is small.

This proof is partially inspired by the transportation lemma of \citet{kaufmann2016complexity}, originally developed for the best arm identification problem with a finite set of alternatives. In contrast, CSO features a continuous solution space with infinitely many feasible solutions, which requires additional technical arguments beyond those used in the finite-armed setting. Related tools have also been proposed for stochastic CSO; for example, combining the transportation lemma of \citet{kaufmann2016complexity} with a minimax argument in \citet{locatelli2018adaptivity} yields a lower bound for the special case $\sigma^2=1$, but the explicit dependence on $\sigma^2$ for general noise levels remains unclear. Alternative approaches based on Fano's inequality \citep{singh2021continuum} or Pinsker's inequality \citep{wang2018optimization} may also be applicable, although they typically require introducing suitable probability metrics (see Chapter~2 of \citet{tsybakov2010}), leading to more technically involved analyses.

\subsection{Two Lemmas on First Hitting Time}

Lemma~\ref{lem:lb} can be used to establish the variance-dependent lower bound associated with stochastic CSO in the variance dichotomy. However, as $\sigma^2$ decreases to zero, this bound also vanishes and becomes loose, particularly in the deterministic CSO regime. This motivates the need for a variance-independent lower bound. The final lower bound in the variance dichotomy will be given by the maximum of these two bounds. The derivation of the variance-independent lower bound relies on the following definition and lemmas, which concern the first hitting time a simulated solution enters a given subregion.

\begin{definition}\label{def:stopping_rule}
	For any subset $\Zcal \subseteq \Xcal$ and any algorithm satisfying Assumption~\ref{ass:adaptive_algo}, define
	\begin{align*}
		\tau(\Zcal) = 
		\begin{cases}
		\min\{t=1,\dots,n: \x_t \in \Zcal\}, & \text{if } \{\x_1,\dots,\x_n\} \cap \Zcal \neq \emptyset, \\
		n+1, & \text{otherwise}.
		\end{cases}
	\end{align*}
\end{definition}

The random variable $\tau(\Zcal)$ represents the first time at which the algorithm simulates a solution in the subregion $\Zcal$. For completeness, we set $\tau(\Zcal)=n+1$ if none of the simulated solutions $\x_1,\dots,\x_n$ lies in $\Zcal$. Clearly, the distribution of $\tau(\Zcal)$ depends on the underlying CSO problem.

The following lemma shows that, under suitable conditions, the distribution of $\tau(\Zcal)$ remains unchanged when the same algorithm is applied to Problems~1 and~2. The detailed proof of the lemma is provided in Section \ref{subsec:proof_lemlb2}.

\begin{lemma}\label{lem:lb2}
	Suppose that $\Zcal$ is a subregion of $\Xcal$ and that $y_{1}(\x)=y_{2}(\x)$ for all $\x \notin \Zcal$ in Problems~1 and~2. Then, for any algorithm satisfying Assumption~\ref{ass:adaptive_algo}, $\Pb_{2}\bigl(\tau(\Zcal)>t\bigr)=\Pb_{1}\bigl(\tau(\Zcal)>t\bigr)$ for all $t$.
\end{lemma}

Lemma~\ref{lem:lb2} assumes that the two objective functions coincide outside the subregion $\Zcal$. As a result, the simulation observations at any solution outside $\Zcal$ have the same distribution under both problems. If none of the first $t$ simulated solutions lies in $\Zcal$, then the joint distribution of $(Y(\x_1),\dots,Y(\x_t))$ is identical under the two problems, which in turn implies that the distribution of $\tau(\Zcal)$ is the same by Assumption~\ref{ass:adaptive_algo}.

In addition to Lemma~\ref{lem:lb2}, which compares the distribution of the first hitting time of a fixed subregion across two different CSO problems, we also require the following lemma, which analyzes first hitting times across different subregions within a single CSO problem. The detailed proof of the lemma is provided in Section \ref{subsec:proof_lemlb3}.

\begin{lemma}\label{lem:lb3}
	Consider $n_1$ mutually disjoint sub-regions $\{ \Zcal_{\kappa}, \kappa =1,\dots,n_1\}$. When an algorithm satisfying Assumption \ref{ass:adaptive_algo} is applied to any fixed CSO problem, for any integer $n_2 \le \min\{n_1/2, n \}$, there must exist a sub-region $\Zcal_{\kappa_1} \in \{ \Zcal_{\kappa}, \kappa =1,\dots,n_1\}$ such that $\Pb( \tau(\Zcal_{\kappa_1}) > n_2 ) \ge 1/2$. 
\end{lemma}

Lemma~\ref{lem:lb3} fixes a CSO problem and considers $n_1$ mutually disjoint subregions of the solution space. By definition, the event $\{\tau(\Zcal_{\kappa}) > n_2\}$ means that none of the first $n_2$ simulated solutions lies in the subregion $\Zcal_{\kappa}$. The lemma therefore examines whether each subregion can receive at least one simulated solution within the first $n_2$ iterations of the algorithm. It shows that there must exist a subregion $\Zcal_{\kappa_1}$ such that, with probability at least $1/2$, none of the first $n_2$ simulated solutions falls in $\Zcal_{\kappa_1}$. 

Since the estimated optimal solution must be one of the simulated solutions by Assumption~\ref{ass:adaptive_algo}, if $n_1 \ge 2n$ and $n_2 = n$, Lemma~\ref{lem:lb3} implies that there exists a subregion that does not contain the estimated optimal solution with probability at least $1/2$. When $n_2 < n$, the lemma instead yields a lower bound on the expectation of the hitting time, and if this bound is sufficiently large, it again implies a non-negligible probability that $\hat{\x}_n^*$ lies outside the subregion.

\section{Variance-Dependent Lower Bounds}\label{sec:var_dep}

In this section, we derive a variance-dependent bound that is closely related to those in the existing literature, except that it explicitly accounts for the dependence on the noise variance.

\subsection{Prototype of Objective Functions}\label{sec:lb_glob_cons}

Selecting an appropriate family of objective functions is central to constructing CSO problems for the lower-bound analysis. In this subsection we modify the family used by \cite{locatelli2018adaptivity} to fit into our context. Let $\rho = \left[M/(2^{2\beta+1} \tilde{M})\right]^{\frac{1}{\beta}}$ and let $\zeta \le 1$ be a strictly positive constant whose value will be specified later. Let $\c_0$ denote a solution in the interior of domain $\Xcal=[0,1]^d$ such that $\{ \x: \| \x - \c_0 \|_{\infty} \le \rho \zeta^{\alpha / \beta} \} \subset \Xcal$. Let $\c$ denote a solution satisfying $\| \c - \c_0 \|_{\infty} \le \rho \zeta^{\alpha / \beta} - \zeta $. Given $\c,\c_0,\zeta$, we define the objective function
\begin{equation}\label{eq:glob_lb_fund}
	\begin{aligned}
		y_{\Gcal} (\x | \c,\c_0,\zeta) = \left\{ \begin{array}{ll}
			\frac{M}{2} \zeta^{\alpha}-\frac{M}{2} \| \x - \c \|_{\infty}^\alpha, &  \| \x - \c \|_{\infty}   \le \zeta \\
			0,     &  \| \x - \c \|_{\infty}   > \zeta  \text{ and }  \| \x - \c_{0} \|_{\infty}   \le \rho \zeta^{\alpha/\beta} \\
			2^{\beta}\tilde{M} \rho^{\beta} \zeta^{\alpha} - 2^{\beta}\tilde{M} \| \x - \c_{0} \|_{\infty}^\beta,      &  \| \x - \c \|_{\infty}   > \zeta  \text{ and }  \| \x - \c_{0} \|_{\infty}   > \rho \zeta^{\alpha/\beta} \\
		\end{array} \right.  
	\end{aligned}  
\end{equation}

\begin{figure}[thbp]
	\centering
	\includegraphics[width=0.99\textwidth]{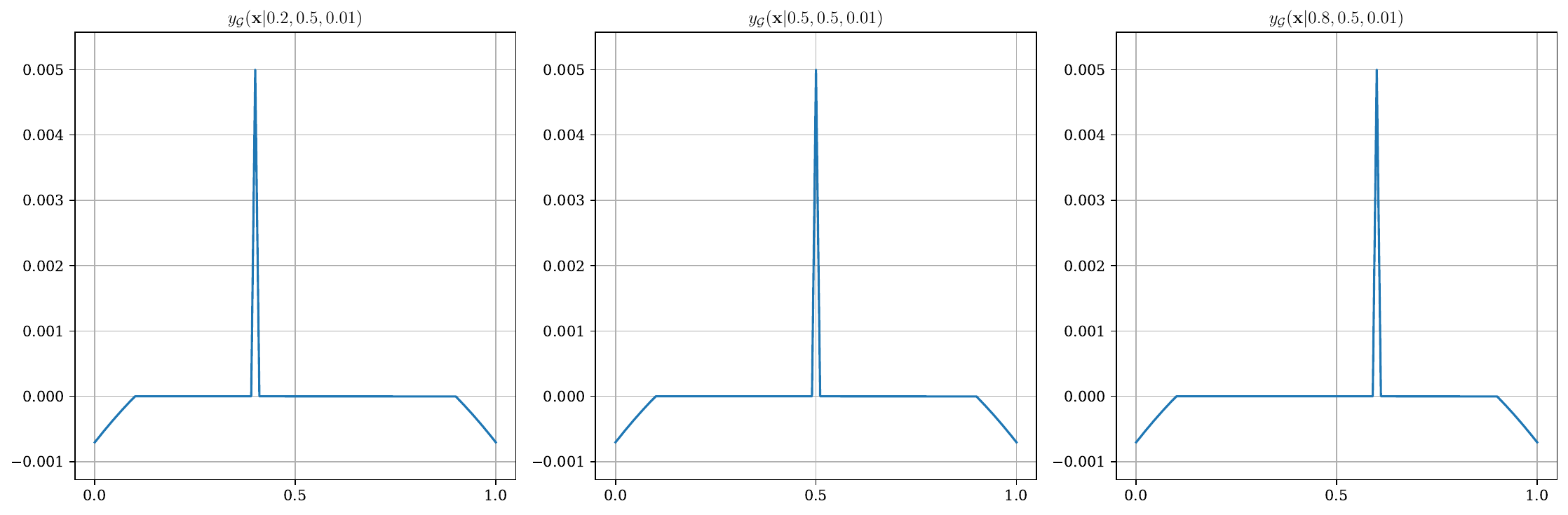}
	\caption{Plot of $y_{\Gcal} (\x | \c,\c_0,\zeta)$ when $\c=0.2,0.5,0.8$, $\c_0=0.5$,  $\zeta=0.01$, $\alpha = 1$, $\beta=2$, $M=1$ and $\tilde{M} = \frac{1}{512}$.}\label{fig:fun_g}
\end{figure}

Figure~\ref{fig:fun_g} illustrates three one-dimensional examples of $y_{\Gcal}(\x \mid \c,\c_0,\zeta)$. The function is continuous and uniquely maximized at $\c$. Its structure involves two local regions: a small neighborhood of radius $\zeta$ centered at $\c$, and a larger neighborhood of radius $\rho \zeta^{\alpha/\beta}$ centered at $\c_0$. The function exhibits three regimes. Within the small neighborhood of $\c$, it decreases at rate $\|\x-\c\|_\infty^{\alpha}$. In the annular region outside the neighborhood of $\c$ but inside the larger neighborhood of $\c_0$, it remains flat. Outside the larger neighborhood of $\c_0$, it decreases at rate $\|\x-\c_0\|_\infty^{\beta}$.

As noted in Section~\ref{sec:form}, the constants $M$ and $\tilde{M}$ may be chosen flexibly, provided they are independent of both the total budget and the decision variable $\x$. Here, we require them to satisfy
\begin{align}\label{ineq:lb_M}
	\frac{M}{2^{4\beta+1} \tilde{M} } \ge 1.
\end{align}

We now summarize a key property of this construction in the following lemma, which shows that the constructed functions satisfy Assumption \ref{ass:smo}. We want to emphasize that, even though the construction of the functions are closely related to that of \cite{locatelli2018adaptivity}, this key property is not explicitly established in their work. The detailed proof of the lemma is provided in Section \ref{subsec:glob_fun1_prop}.
\begin{lemma}\label{lem:glob_fun1_prop}
	The function $y_{\Gcal}(\x \mid \c,\c_0,\zeta)$ satisfies Assumption~\ref{ass:smo} for $\alpha \le \beta$.
\end{lemma}

Throughout the rest of this section, the objective function of each constructed CSO problem will take the form $y_{\Gcal}(\x \mid \c,\c_0,\zeta)$, with $\c$, $\c_0$, and $\zeta$ chosen based on the types of lower bounds.

\subsection{Variance-Dependent Lower Bound}\label{sec:lb_glob}

This subsection establishes a variance-dependent minimax lower bound on the optimization error. The resulting bound matches the convergence rates with respect to the budget $n$ obtained in the existing literature, while explicitly characterizing the dependence on the noise variance $\sigma^2$, which has not previously been made explicit; for example, the lower bound derived by \citet{locatelli2018adaptivity} is obtained under the normalization $\sigma^2=1$. Our analysis follows a different technical route, relying on Lemma~\ref{lem:lb} to accommodate general noise levels in a unified and rigorous manner.

We proceed as follows to establish the variance-dependent lower bound. First, we construct a collection of CSO problems, denoted by $\Gcal_{\alpha,\beta,1}$. Let $\psi_1$ be a positive integer (defined precisely in the proof) that is of order $(\sigma^2/n)^{-1/[\alpha(D+2)]}$, where $D = d(1/\alpha-1/\beta)$, and define $\zeta_1 = (2\psi_1)^{-1}$. Let $\bkappa=(\kappa_1,\dots,\kappa_d)$ denote a multi-index. We partition the domain $[0,1]^d$ into $\psi_1^d$ equal-sized and mutually disjoint subregions indexed by $\bkappa$:
\[
\{\Zcal_{\bkappa}:\kappa_l=1,2,\dots,\psi_1,\; l=1,2,\dots,d\}.
\]
Let $\c_{\bkappa}$ denote the center of subregion $\Zcal_{\bkappa}$, and let $\Zcal_{\bkappa_0}$ be a subregion located at the center of $[0,1]^d$. We focus on those subregions satisfying the distance condition
\begin{align}\label{ineq:lb_glob_dist1}
	\|\c_{\bkappa}-\c_{\bkappa_0}\|_{\infty} \le \rho \zeta_1^{\alpha/\beta}-\zeta_1.
\end{align}
For each such subregion $\Zcal_{\bkappa}$, we construct a CSO problem $\CSO_{\bkappa}$ whose objective function is given by $y_{\Gcal}(\x \mid \c_{\bkappa},\c_{\bkappa_0},\zeta_1)$ as defined in \eqref{eq:glob_lb_fund}. The resulting collection of problems forms the set $\Gcal_{\alpha,\beta,1}$. A key feature of these problems is that the optimal solution of $\CSO_{\bkappa}$ coincides with the center $\c_{\bkappa}$ of the corresponding subregion $\Zcal_{\bkappa}$. Accordingly, the estimated optimal solution $\hat{\x}_n^*$ for $\CSO_{\bkappa}$ achieves satisfactory performance if $\hat{\x}_n^*\in\Zcal_{\bkappa}$.

Next, we construct a benchmark CSO problem. Define the function
\begin{align}\label{eq:glob_bench}
	y_{\Bcal}(\x \mid \c_0,\zeta)=
	\begin{cases}
	0, & \|\x-\c_0\|_\infty \le \rho \zeta^{\alpha/\beta},\\[4pt]
	2^{\beta}\tilde{M}\rho^{\beta}\zeta^{\alpha}
	-2^{\beta}\tilde{M}\|\x-\c_0\|_\infty^{\beta}, 
	& \|\x-\c_0\|_\infty > \rho \zeta^{\alpha/\beta}.
	\end{cases}
\end{align}
The objective function in \eqref{eq:glob_bench} differs from that in \eqref{eq:glob_lb_fund} only within a localized neighborhood of the maximizer $\c$. We define the benchmark CSO problem, denoted by $\CSO_{\mathbf{0}}$, by setting $\c_0=\c_{\bkappa_0}$ and $\zeta=\zeta_1$ in \eqref{eq:glob_bench}. The set of optimal solutions for $\CSO_{\mathbf{0}}$ consists of all points lying within distance $\rho \zeta_1^{\alpha/\beta}$ of $\c_{\bkappa_0}$. Since any solution in a subregion satisfying the distance condition \eqref{ineq:lb_glob_dist1} lies within this radius, all points in these subregions are optimal under the benchmark objective.

Finally, we analyze the event $\{\hat{\x}_n^* \in \Zcal_{\bkappa}\}$ for each problem $\CSO_{\bkappa}$. Among the subregions satisfying \eqref{ineq:lb_glob_dist1}, we can identify a subregion $\Zcal_{\bar{\bkappa}}$ such that the probability of $\{\hat{\x}_n^* \in \Zcal_{\bar{\bkappa}}\}$ is at most $1/2$ and the expected number of observations in this subregion is minimal under the benchmark CSO. To see this, consider the two subregions with the smallest expected numbers of observations. If both were selected with probability strictly greater than $1/2$, the sum of their probabilities would exceed one, which is impossible. Hence, at least one of these minimally sampled subregions must satisfy $\Pb\bigl(\hat{\x}_n^* \in \Zcal_{\bar{\bkappa}}\bigr) \le 1/2$.

We now transition from the benchmark problem $\CSO_{\mathbf{0}}$ to the perturbed problem $\CSO_{\bar{\bkappa}}$, in which the true optimizer is relocated to $\Zcal_{\bar{\bkappa}}$. Applying Lemma~\ref{lem:lb}, we show that under $\CSO_{\bar{\bkappa}}$, $\Pb\bigl(\hat{\x}_n^* \in \Zcal_{\bar{\bkappa}}\bigr) \le 2/3$. Thus, with probability at least $1/3$, the estimated optimal solution lies outside the true optimal subregion, yielding an optimality gap of order at least $(\sigma^2/n)^{1/(D+2)}$. The following proposition formalizes this variance-dependent lower bound. The detailed proof is provided in Section \ref{subsec:proof_prop:lb_local_large}.

\begin{proposition}\label{prop:lb_local_large}
	For any adaptive algorithm satisfying Assumption~\ref{ass:adaptive_algo}, there exists a constant $b_{\Gcal,1}$, independent of $n$ and $\sigma^2$, such that for all $n \ge b_{\Gcal,1}\sigma^2$,
	\begin{align*}
		\sup_{\CSO \in \Gcal_{\alpha,\beta,1}} 
		\Eb_{\CSO}\bigl[y(\x^*)-y(\hat{\x}_n^*)\bigr]
		\ge b_{\Gcal,2}\left({\sigma^2}/{n}\right)^{\frac{1}{D+2}},
	\end{align*}
	where $D = d(1/\alpha-1/\beta)$ is the effective dimension and $b_{\Gcal,2}$ is a positive constant independent of $n$ and $\sigma^2$.
\end{proposition}

A related result was previously established by \citet{locatelli2018adaptivity} for the special case $\sigma^2=1$. In contrast, the lower bound $\Omega\bigl((\sigma^2/n)^{1/(D+2)}\bigr)$ in Proposition~\ref{prop:lb_local_large} explicitly characterizes the dependence on the noise variance. The simulation budget $n$ and the variance $\sigma^2$ jointly determine the magnitude of the bound; in particular, the bound is invariant under simultaneous scaling of $\sigma^2$ and $n$ by the same factor. Consequently, CSO problems with smaller noise variance are effectively equivalent, in terms of optimization difficulty, to problems with larger variance but proportionally larger simulation budgets.

Furthermore, we see that when $\alpha=\beta$, the minimax lower bound becomes $\Omega\bigl((\sigma^2/n)^{1/2}\bigr)$ and the convergence rate becomes independent of the problem dimension. However, this dimension-independent rate does not imply that dimensionality is irrelevant. As indicated by \eqref{ineq:glob_bg2_lb} in the proof of Proposition~\ref{prop:lb_local_large}, the constant $b_{\Ccal,2}$ scales proportionally to $2^{d/2}$, exhibiting an exponential dependence on the dimension $d$. Consequently, although the asymptotic rate remains $(\sigma^2/n)^{1/2}$, the magnitude of the lower bound becomes substantially larger in high-dimensional stochastic CSO problems.

In the low-noise regime $\sigma^2 \to 0$, however, the bound in Proposition~\ref{prop:lb_local_large} vanishes and thus becomes loose, providing no meaningful performance guarantee. Although variance-dependent lower bounds of this form are prevalent in the CSO literature, the role of the noise level itself has typically not been examined carefully. Existing analyses often implicitly assume a fixed variance or focus on asymptotic regimes where the budget $n$ is sufficiently large, whereas settings in which the variance is small or the simulation budget is limited have received far less attention. This practically relevant regime motivates the variance-independent analysis developed in the next subsection.

\section{Variance-Independent Lower Bound when $\alpha<\beta$}\label{sec:var_ind_uneq}

To derive the variance-independent lower bounds on the optimization error, it turns out that the distinction between $\alpha<\beta$ and $\alpha=\beta$ plays a critical role. In particular, we have to construct different types of objective functions to derive the lower bounds. In this section, we focus on the case of $\alpha<\beta$ and leave the case of $\alpha=\beta$ to the next section.

To derive the variance-independent minimax lower bound for the case of $\alpha<\beta$, we use the same family of objective functions developed in Section \ref{sec:lb_glob_cons}. The construction follows a similar three-step procedure to that used for the variance-dependent bound, with key modifications in the parameter choices. In particular, the domain $[0,1]^d$ is partitioned into $\psi_2^d$ subregions, where $\psi_2$ is a positive integer of order $n^{1/(\alpha D)}$. In contrast to $\psi_1$, which depends on both the simulation budget and the noise variance, $\psi_2$ depends solely on the budget $n$. This change in scaling is the essential ingredient that removes the dependence on $\sigma^2$ in the resulting lower bound.

Let $\zeta_2 = (2\psi_2)^{-1}$. Following the same construction as in the variance-dependent case, we select a central subregion $\Zcal_{\bkappa_0}$ among the $\psi_2^d$ subregions. For each subregion $\Zcal_{\bkappa}$ satisfying
\[
\|\c_{\bkappa}-\c_{\bkappa_0}\|_\infty \le \rho \zeta_2^{\alpha/\beta}-\zeta_2,
\]
we define a CSO problem $\CSO_{\bkappa}$ with objective function $y_{\Gcal}(\x \mid \c_{\bkappa},\c_{\bkappa_0},\zeta_2)$ as in \eqref{eq:glob_lb_fund}. The resulting family of problems constitutes the set $\Gcal_{\alpha,\beta,2}$. We then construct the corresponding benchmark CSO problem by setting $\c_0=\c_{\bkappa_0}$ and $\zeta=\zeta_2$ in \eqref{eq:glob_bench}.

The third step in establishing the variance-independent lower bound again analyzes the event $\{\hat{\x}_n^* \in \Zcal_{\bkappa}\}$ for each $\CSO_{\bkappa}$, but now relies on Lemmas~\ref{lem:lb2} and~\ref{lem:lb3} rather than Lemma~\ref{lem:lb}. For sufficiently large $n$, the set $\Gcal_{\alpha,\beta,2}$ contains at least $2n$ CSO problems whose corresponding subregions are mutually disjoint. By Lemma~\ref{lem:lb3}, there exists a subregion $\Zcal_{\tilde{\bkappa}}$ such that, under the benchmark CSO problem, its first hitting time exceeds $n$ with probability greater than $1/2$. 

Since the objective function of $\CSO_{\tilde{\bkappa}}$ differs from that of the benchmark problem only within the subregion $\Zcal_{\tilde{\bkappa}}$, Lemma~\ref{lem:lb2} implies that the probability of the first hitting time exceeding $n$ remains greater than $1/2$ under $\CSO_{\tilde{\bkappa}}$. Consequently, with probability exceeding $1/2$, none of the $n$ simulated solutions lies in $\Zcal_{\tilde{\bkappa}}$. Under Assumption~\ref{ass:adaptive_algo}, the estimated optimal solution $\hat{\x}_n^*$ must be one of the simulated solutions and therefore cannot lie in $\Zcal_{\tilde{\bkappa}}$ on this event. This yields an optimality gap bounded below by a constant multiple of $n^{-1/D}$, leading to a variance-independent minimax lower bound, which is formalized in the following proposition. The detailed proof is provided in Section \ref{subsec:lb_glob_fun2}.

\begin{proposition}\label{prop:lb_local_large2}
	For any adaptive algorithm satisfying Assumption~\ref{ass:adaptive_algo}, there exists a constant $b_{\Gcal,3}$, independent of $n$ and $\sigma^2$, such that for all $n \ge b_{\Gcal,3}$,
	\begin{align*}
		\sup_{\CSO \in \Gcal_{\alpha,\beta,2}} 
		\Eb_{\CSO}\bigl[y(\x^*)-y(\hat{\x}_n^*)\bigr] 
		\ge b_{\Gcal,4}\cdot n^{-{\frac{1}{D}}},
	\end{align*}
	where $b_{\Gcal,4}$ is a constant independent of $n$ and $\sigma^2$.
\end{proposition}

The lower bound $\Omega(n^{-1/D})$ in Proposition~\ref{prop:lb_local_large2} is independent of the noise variance $\sigma^2$ and therefore remains robust in the low-noise regime. For a fixed budget $n$, the bound is governed by the effective dimensionality parameter $D$, which incorporates both the problem dimension $d$ and the shape parameters $\alpha$ and $\beta$ of the objective function. Notably, this rate scales as $n^{-1/D}$ rather than the $n^{-1/(D+2)}$ rate that typically arises in variance-dependent settings. To the best of our knowledge, this variance-independent lower bound is new in the CSO literature.

\section{Variance-Independent Lower Bounds When $\alpha=\beta$}\label{sec:var_ind_eq}

The variance-independent lower bound established in Proposition~\ref{prop:lb_local_large2} for the case $\alpha<\beta$ cannot be directly extended to $\alpha=\beta$. When $\alpha=\beta$, the effective dimension parameter $D$ equals zero, rendering the expression $\Omega(n^{-1/D})$ ill-defined. As $D \to 0$, the rate $n^{-1/D}$ decays faster than any polynomial rate $n^{-p}$ for fixed $p>0$, suggesting that the variance-independent lower bound in this regime follows a super-polynomial, and in fact exponential, convergence behavior. 

We therefore develop a separate construction to establish an exponential variance-independent lower bound. The proof differs from the previous case in two key aspects: (i) it employs a hierarchical partition of the search space into nested subregions across multiple resolution levels; and (ii) instead of relying on a single benchmark CSO problem, it constructs a sequence of benchmark problems corresponding to successive refinement levels.

\subsection{Objective Functions}

The objective functions for the CSO problems are built upon this hierarchical decomposition. Let $\bar{a}$ be a positive integer that scales with the budget $n$ (specified precisely in Section \ref{subsec:detail_partition}). For each level $a \in \{0,1,\dots,\bar{a}\}$, define $n_a = 5^a$. The domain $[0,1]^d$ is partitioned into $n_a^d$ mutually disjoint and congruent subregions at level $a$. Let $\Zcal_{\bkappa_a,a}$ denote a subregion at level $a$, where the multi-index $\bkappa_a=(\kappa_{1,a},\dots,\kappa_{d,a})$ satisfies $\kappa_{l,a}\in\{1,\dots,n_a\}$. These subregions form a nested hierarchy: for each subregion at level $a$, there exists a unique parent subregion at level $a-1$, indexed by $\Pcal(\bkappa_a)$, such that $\Zcal_{\bkappa_a,a} \subset \Zcal_{\Pcal(\bkappa_a),a-1}$. More generally, we denote by $\Pcal^{\ell}(\bkappa_a)$ the ancestor index at level $a-\ell$, which tracks the lineage of each subregion across successive refinement levels. 

The construction of objective functions is based on a hierarchical selection of subregions across levels $a=0,1,2,\dots,\bar{a}$. For each selected subregion, we define a local component function, and the overall objective function is obtained by summing these components. The selection procedure follows a recursive nesting structure.

Let $\Ucal_a$ denote the collection of selected subregions at level $a$. At the initial level $a=0$, we set $\Ucal_0=\{\Zcal_{\mathbf{1},0}\}$. For each subsequent level $a=1,2,\dots,\bar{a}$, a subregion $\Zcal_{\bkappa_a,a}$ is included in $\Ucal_a$ if and only if its parent subregion $\Zcal_{\Pcal(\bkappa_a),a-1}$ belongs to $\Ucal_{a-1}$ and the distance between their centers does not exceed $2\gamma_a$, where $\gamma_a = (2n_a)^{-1}$.

To illustrate the construction, consider the one-dimensional case in which $\Zcal_{1,0}=(0,1)$. At level $a=1$, the domain is partitioned into five subregions. The distances between the centers of the subregions $(1/5,2/5)$, $(2/5,3/5)$, and $(3/5,4/5)$ and the center of $\Zcal_{1,0}$ are all within the threshold $2\gamma_1$. In contrast, the subregions $(0,1/5)$ and $(4/5,1)$ lie beyond this threshold. Consequently, only the three central subregions are included in $\Ucal_1$.

\subsubsection{Component Functions.}

We now describe the construction of the individual component functions associated with each selected subregion. Analogous to the condition in \eqref{ineq:lb_M} for the case $\alpha<\beta$, we impose the following requirement in this section:
\begin{align}\label{ineq:lb_M3}
	\frac{M}{\tilde{M}} \ge (1-5^{-\alpha})^{-1}\left(1+\frac{1}{5^{\alpha}-3^{\alpha}}\right)80^{\alpha}.
\end{align}
We further define the scaling constant $\hat{M} = \frac{40^{\alpha}}{5^{\alpha}-3^{\alpha}}\,\tilde{M}$, and let $\c_{\bkappa_a,a}$ denote the center of the subregion $\Zcal_{\bkappa_a,a}$. 

For the initial subregion $\Zcal_{\mathbf{1},0}\in\Ucal_0$, we define the corresponding component function as
\begin{equation}\label{eq:scv_fun20}
	y_{\CS}(\x \mid \mathbf{1},0)=
	\begin{cases}
	\hat{M}\bigl(5^{\alpha}-3^{\alpha}\bigr)\gamma_1^{\alpha}, 
	& \quad \|\x-\c_{\mathbf{1},0}\|_\infty \le 3\gamma_1,\\
	\hat{M}5^{\alpha}\gamma_1^{\alpha}-\hat{M}\|\x-\c_{\mathbf{1},0}\|_\infty^{\alpha}, 
	& \quad 3\gamma_1 \le \|\x-\c_{\mathbf{1},0}\|_\infty \le 5\gamma_1 .
	\end{cases}
\end{equation}
Since $\c_{\mathbf{1},0}$ is the center of the unit hypercube $[0,1]^d$ and $\gamma_1=1/10$, we have $[0,1]^d=\{\x:\|\x-\c_{\mathbf{1},0}\|_\infty \le 5\gamma_1\}$. The function defined in \eqref{eq:scv_fun20} exhibits two distinct spatial regimes. Within the neighborhood of radius $3\gamma_1$ centered at $\c_{\mathbf{1},0}$, the function remains constant. Outside this region, it decreases proportionally to $\|\x-\c_{\mathbf{1},0}\|_\infty^{\alpha}$. This structure is directly analogous to the benchmark objective function in \eqref{eq:glob_bench}, which underpins the constructions used in Propositions~\ref{prop:lb_local_large} and~\ref{prop:lb_local_large2}.

For subregions $\Zcal_{\bkappa_a,a}\in\Ucal_a$ at intermediate levels, the corresponding component functions are defined as
\begin{equation}\label{eq:scv_fun2}
	y_{\CS}(\x \mid \bkappa_a,a)=
	\begin{cases}
	\hat{M}\bigl(5^{\alpha}-3^{\alpha}\bigr)\gamma_{a+1}^{\alpha}, 
	& \|\x-\c_{\bkappa_a,a}\|_\infty \le 3\gamma_{a+1},\\
	\hat{M}5^{\alpha}\gamma_{a+1}^{\alpha}
	-\hat{M}\|\x-\c_{\bkappa_a,a}\|_\infty^{\alpha}, 
	& 3\gamma_{a+1} < \|\x-\c_{\bkappa_a,a}\|_\infty \le 5\gamma_{a+1},\\
	0, 
	& \|\x-\c_{\bkappa_a,a}\|_\infty > 5\gamma_{a+1}.
	\end{cases}
\end{equation}
The function in \eqref{eq:scv_fun2} is a rescaled and localized version of the base function in \eqref{eq:scv_fun20}. It preserves the same structural decay pattern while being supported on progressively smaller regions corresponding to finer levels of the hierarchical partition.

At the finest refinement level, each subregion $\Zcal_{\bkappa_{\bar{a}},\bar{a}}\in\Ucal_{\bar{a}}$ is assigned a terminal component function defined by
\begin{equation}\label{eq:scv_fun3}
	y_{\mathcal{D}}(\x \mid \bkappa_{\bar{a}},\bar{a})=
	\begin{cases}
	\hat{M}\gamma_{\bar{a}}^{\alpha}
	-\hat{M}\|\x-\c_{\bkappa_{\bar{a}},\bar{a}}\|_\infty^{\alpha}, 
	& \|\x-\c_{\bkappa_{\bar{a}},\bar{a}}\|_\infty \le \gamma_{\bar{a}},\\
	0, 
	& \|\x-\c_{\bkappa_{\bar{a}},\bar{a}}\|_\infty > \gamma_{\bar{a}}.
	\end{cases}
\end{equation}
Within the subregion $\Zcal_{\bkappa_{\bar{a}},\bar{a}}$, this function decreases proportionally to $\|\x-\c_{\bkappa_{\bar{a}},\bar{a}}\|_\infty^{\alpha}$ and vanishes outside the subregion.

\subsubsection{Objective Functions.}

For each selected subregion $\Zcal_{\bkappa_a,a}$ at levels $a\le \bar{a}-1$, we define the objective function
\begin{align}\label{eq:lb_scv_benchfun}
	y_{\Ccal}(\x|\bkappa_{a},a) =  y_{\CS} (\x | \bkappa_{a},a) + y_{\CS} (\x | \Pcal(\bkappa_{a}),a-1) + \dots + y_{\CS} (\x | \Pcal^{a-1}(\bkappa_{a}),1) + y_{\CS} (\x | \mathbf{1},0) ,
\end{align}
which is the sum of $a+1$ component functions defined in \eqref{eq:scv_fun20} and \eqref{eq:scv_fun2}. By construction, the subregions associated with these component functions form a nested hierarchy, with each parent subregion containing the target subregion $\Zcal_{\bkappa_a,a}$.

For each selected subregion $\Zcal_{\bkappa_{\bar{a}},\bar{a}}$ at the finest level $\bar{a}$, we define
\begin{align}\label{eq:lb_scv_fun3}
	y_{\Ccal}(\x|\bkappa_{\bar{a}},\bar{a}) = y_{\mathcal{D}} (\x | \bkappa_{\bar{a}},\bar{a}) + y_{\CS} (\x | \Pcal(\bkappa_{\bar{a}}),\bar{a}-1) + \dots + y_{\CS} (\x | \Pcal^{\bar{a}-1}(\bkappa_{\bar{a}}),1) + y_{\CS} (\x | \mathbf{1},0),
\end{align}
which consists of $\bar{a}+1$ component functions defined in \eqref{eq:scv_fun20} to \eqref{eq:scv_fun3}. The component subregions again follow a nested structure, and the function $y_{\Ccal}(\x \mid \bkappa_{\bar{a}},\bar{a})$ attains its unique maximum at the center $\c_{\bkappa_{\bar{a}},\bar{a}}$ of the target subregion.

\begin{figure}[htbp]
	\centering
	\includegraphics[width=0.82\textwidth]{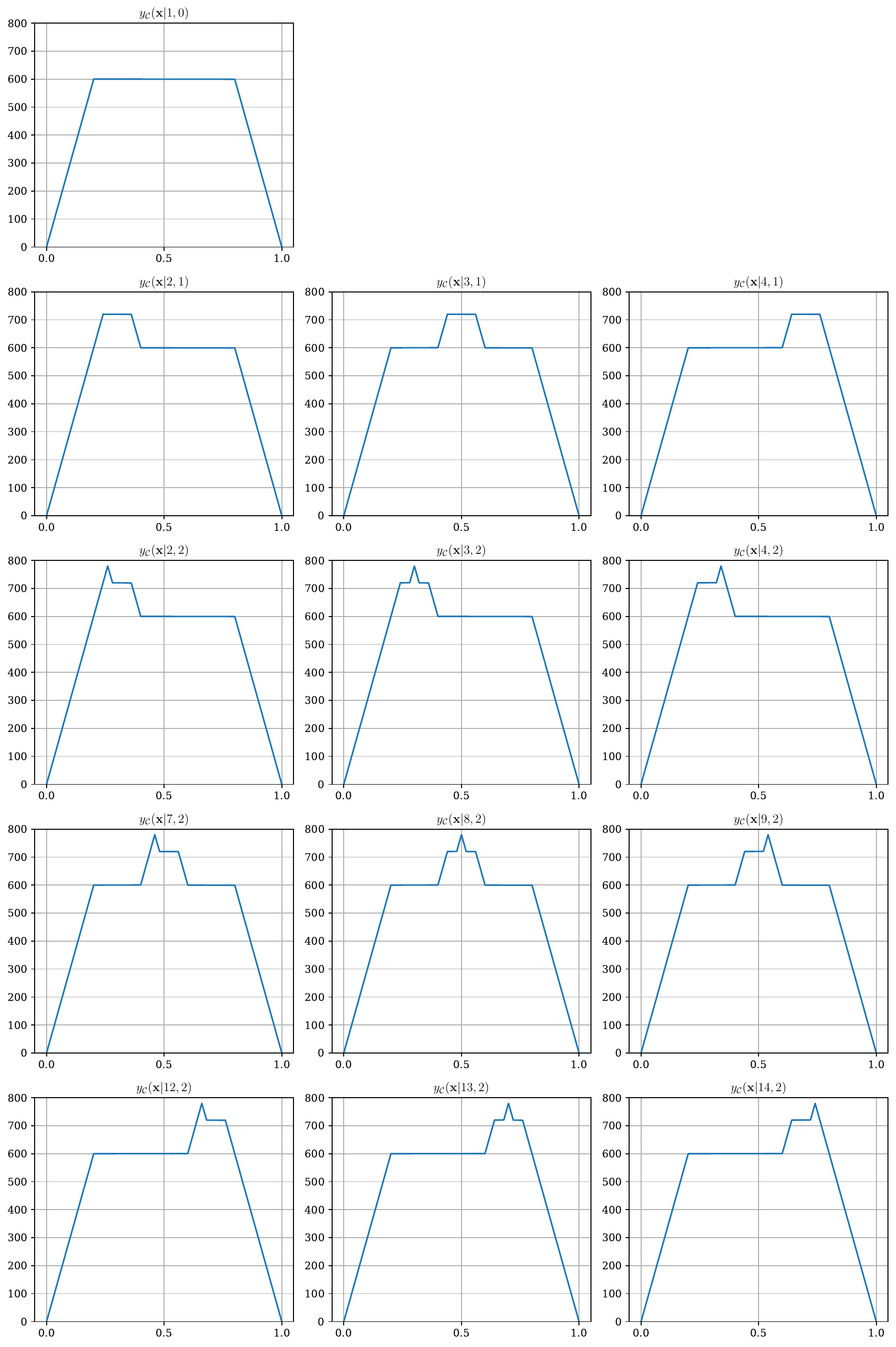}
	\caption{Plot of $y_{\Ccal} (\x | \bkappa_a,a)$ when the highest partition level is $\bar{a}=2$.}\label{fig:fun_c}
\end{figure}

Figure~\ref{fig:fun_c} provides one-dimensional illustrations of the constructed objective functions for the case $\bar{a}=2$. The top two rows correspond to functions of the form \eqref{eq:lb_scv_benchfun} at levels $a=0$ and $a=1$, while the bottom three rows display examples of \eqref{eq:lb_scv_fun3} at the finest level. These plots visually demonstrate the hierarchical nesting of component regions and how successive refinements localize the maximizer while preserving the overall smooth structure.

\subsubsection{Key Property.}

The following lemma establishes the essential regularity property of the constructed objective functions. A detailed proof is provided in Section \ref{subsec:lem:scv_fun1_prop}.

\begin{lemma}\label{lem:scv_fun1_prop}
	For any subregion $\Zcal_{\bkappa_{\bar{a}},\bar{a}}\in\Ucal_{\bar{a}}$, the function $y_{\Ccal}(\x \mid \bkappa_{\bar{a}},\bar{a})$ satisfies Assumption~\ref{ass:smo} with $\alpha=\beta$.
\end{lemma}

For each selected subregion, let $\CSO_{\bkappa_a,a}$ denote the CSO problem with objective function $y_{\Ccal}(\x \mid \bkappa_a,a)$. We define $\Ccal_{\alpha,2}$ as the class of CSO problems corresponding to subregions in $\Ucal_{\bar{a}}$ at the finest level $\bar{a}$. The CSO problems associated with subregions in $\Ucal_a$ for $a<\bar{a}$ serve as benchmark problems in the subsequent analysis.

\subsection{Variance-Independent Lower Bound}\label{sec:lb_scv2}

The fundamental strategy for deriving the variance-independent lower bound is to identify, at each refinement level of the hierarchical partition, a subregion whose expected first hitting time is sufficiently large, and to propagate this property across levels using the nested structure. To this end, for each $\Zcal_{\bkappa_a,a}\in\Ucal_a$, we define the set of child subregions $\Hcal_{\bkappa_a,a}$ as the collection of subregions in $\Ucal_{a+1}$ whose parent is $\Zcal_{\bkappa_a,a}$. By construction, $|\Hcal_{\bkappa_a,a}|=3^d$. The analysis begins with the benchmark problem $\CSO_{\mathbf{1},0}$ at the coarsest level and proceeds recursively through the hierarchy.

Let $\Delta=\lfloor 3^d/2 \rfloor$. When solving the benchmark problem $\CSO_{\mathbf{1},0}$, Lemma~\ref{lem:lb3} guarantees the existence of a child subregion $\Zcal_{\bar{\bkappa}_1,1}\in\Hcal_{\mathbf{1},0}$ such that its first hitting time exceeds $\Delta$ with probability greater than $1/2$. Consequently, the expected hitting time satisfies $\Eb[\tau(\Zcal_{\bar{\bkappa}_1,1})]>\Delta/2$ under the probability measure induced by $\CSO_{\mathbf{1},0}$. Moreover, since the objective functions $y_{\Ccal}(\x\mid \bar{\bkappa}_1,1)$ and $y_{\Ccal}(\x\mid \mathbf{1},0)$ coincide outside the subregion $\Zcal_{\bar{\bkappa}_1,1}$, Lemma~\ref{lem:lb2} implies that the same lower bound on the expected hitting time holds under the problem $\CSO_{\bar{\bkappa}_1,1}$, namely $\Eb[\tau(\Zcal_{\bar{\bkappa}_1,1})]>\Delta/2$.

We now extend this argument to the next refinement level. When solving the problem $\CSO_{\bar{\bkappa}_1,1}$, Lemma~\ref{lem:lb3} again guarantees the existence of a child subregion $\Zcal_{\bar{\bkappa}_2,2}\in\Hcal_{\bar{\bkappa}_1,1}$ such that its first hitting time exceeds $\tau(\Zcal_{\bar{\bkappa}_1,1})+\Delta$ with probability greater than $1/2$. Consequently,
\[
\Eb\bigl[\tau(\Zcal_{\bar{\bkappa}_2,2})\bigr]
\ge \Eb\bigl[\tau(\Zcal_{\bar{\bkappa}_1,1})\bigr] + \Delta/2
\ge \Delta,
\]
and this inequality continues to hold under the problem $\CSO_{\bar{\bkappa}_2,2}$. Repeating this argument recursively across refinement levels yields a sequence of subregions such that, at level $a$,
\[
\Eb\bigl[\tau(\Zcal_{\bar{\bkappa}_a,a})\bigr] \ge a\Delta/2
\]
under the corresponding problem $\CSO_{\bar{\bkappa}_a,a}$. For sufficiently large $a$, the expected hitting time exceeds the total budget $n$, implying that, with high probability, none of the $n$ simulated solutions falls within the associated subregion.

The rigorous proof is necessarily more delicate, as it must ensure that all statements involving hitting times remain well defined under the finite simulation budget. To this end, we restrict attention to a sequence of suitably defined ``regular'' events on which hitting times are properly characterized, and carry out the analysis conditional on these events. Within this framework, the lower bounds on the expected hitting times are expressed as sums of the probabilities of these regular events across successive refinement levels. 

We show that the probability of the regular event at the finest level must be small; otherwise, the resulting lower bound on the expected hitting time would exceed $n+1$, contradicting the definition of the stopping time. Consequently, there must exist a subregion $\Zcal_{\bar{\bkappa}_{\bar{a}},\bar{a}}$ at level $\bar{a}$ such that, when solving the problem $\CSO_{\bar{\bkappa}_{\bar{a}},\bar{a}}$, none of the $n$ simulated solutions falls within $\Zcal_{\bar{\bkappa}_{\bar{a}},\bar{a}}$ with high probability. This absence of sampling within the subregion leads to a substantial optimization error. The following proposition formalizes the resulting variance-independent lower bound; its proof is provided in Section \ref{subsec:lb_scv_fun2}.

\begin{proposition}\label{prop:lb_scv2}
	For any adaptive algorithm satisfying Assumption~\ref{ass:adaptive_algo}, there exists a constant $b_{\Ccal,3}$ such that for all sufficiently large $n \ge b_{\Ccal,3}$,
	\begin{align*}
		\sup_{\CSO \in \Ccal_{\alpha,2}} 
		\Eb_{\CSO}\bigl[y(\x^*)-y(\hat{\x}_n^*)\bigr] 
		\ge b_{\Ccal,4}\exp(-b_{\Ccal,5}n),
	\end{align*}
	where $b_{\Ccal,4}$ and $b_{\Ccal,5}$ are positive constants independent of $n$ and $\sigma^2$.
\end{proposition}

Proposition~\ref{prop:lb_scv2} establishes a variance-independent minimax lower bound that decays exponentially with the simulation budget $n$. Such exponential lower bounds are sparse in the stochastic CSO literature (where $\sigma^2>0$), but they are known to arise under substantially stronger conditions in deterministic optimization (where $\sigma^2=0$). In particular, for strongly concave deterministic optimization problems with access to accurate gradient information, exponential convergence and corresponding lower bounds are well established \citep{nemirovskij1983problem,arjevani2016lower}. These results, however, rely critically on strong concavity and gradient availability, and therefore do not apply to nonconcave objectives or to CSO settings where gradients are not observed. In deterministic CSO with only Lipschitz continuity, the classical minimax lower bound is $\Omega(n^{-1/d})$, reflecting the intrinsic difficulty of exploration-based optimization \citep{bull2011convergence,malherbe2017global,singh2021continuum,fourati2025every}. Compared with this baseline, the exponential lower bound in Proposition~\ref{prop:lb_scv2} represents a dramatic improvement when $\alpha=\beta$.

Beyond the lower bound literature, it has been shown that under additional shape conditions, specifically $\alpha=1$ and $\beta\ge1$, certain algorithms can achieve faster polynomial rates than $n^{-1/d}$ and, in some favorable cases, even exponential decay of the optimization error in deterministic settings \citep{munos2014,malherbe2017global}. The variance-independent lower bounds established in Proposition~\ref{prop:lb_local_large2} (for $\beta>\alpha=1$) and Proposition~\ref{prop:lb_scv2} (for $\beta=\alpha=1$) match these algorithmic rates, thereby providing the first minimax justification for their optimality under such shape constraints.

Finally, we emphasize that the exponential lower bound in Proposition~\ref{prop:lb_scv2} does not eliminate the curse of dimensionality. The constant $b_{\Ccal,5}$ governing the exponential rate scales as $3^{-d}$ and thus decays exponentially with the dimension $d$. Consequently, although the asymptotic convergence rate is exponential in $n$, high dimensionality substantially weakens the effective bound, embedding the curse of dimensionality within the constant factor.

\section{Overall Lower Bound and Variance Dichotomy}\label{sec:lb_overall}

In this section, we establish minimax lower bounds that characterize the variance dichotomy in the regimes $\alpha < \beta$ and $\alpha=\beta$. 
In each regime, we combine the variance-dependent and variance-independent bounds derived in the preceding sections to obtain a complete characterization of the variance dichotomy.

\subsection{Overall Lower Bound When $\alpha<\beta$}\label{subsec:over_lb_uneq}

Let $\Gcal_{\alpha,\beta,0}$ denote the class of CSO problems satisfying Assumptions~\ref{ass:unique}--\ref{ass:smo} for fixed parameters $\alpha$ and $\beta$ with $\alpha < \beta$. By construction, both $\Gcal_{\alpha,\beta,1}$ with $\alpha < \beta$ and $\Gcal_{\alpha,\beta,2}$ are subsets of $\Gcal_{\alpha,\beta,0}$, which represents the full class of CSO problems characterized by the shape parameters $\alpha$ and $\beta$. Since the minimax lower bound over $\Gcal_{\alpha,\beta,0}$ is defined as the supremum of the optimization error across all problems in this class, it must be no smaller than the corresponding bounds over any of its subsets. Consequently, Propositions~\ref{prop:lb_local_large} and~\ref{prop:lb_local_large2} together yield the following theorem, which establishes the overall minimax lower bound governing the variance dichotomy.

\begin{theorem}\label{th:lb_local_large_s}
	For any adaptive algorithm satisfying Assumption~\ref{ass:adaptive_algo}, there exist constants $b_{\Gcal,1}$ and $b_{\Gcal,3}$, independent of $n$ and $\sigma^2$, such that for all $n \ge \max\{b_{\Gcal,1}\sigma^2,\, b_{\Gcal,3}\}$,
	\begin{align*}
		\sup_{\CSO \in \Gcal_{\alpha,\beta,0}} 
		\Eb_{\CSO}\bigl[y(\x^*)-y(\hat{\x}_n^*)\bigr]
		\ge b_{\Gcal,5}\cdot
		\max\left\{
		\left({\sigma^2}/{n}\right)^{\frac{1}{D+2}},\,
		n^{-\frac{1}{D}}
		\right\},
	\end{align*}
	where $b_{\Gcal,5}=\max\{b_{\Gcal,2},\, b_{\Gcal,4}\}$.
\end{theorem}

The minimax lower bound in Theorem~\ref{th:lb_local_large_s} is given by the maximum of two terms, reflecting the two regimes of the variance dichotomy. As a special case, for deterministic CSO with $\sigma^2=0$, the bound reduces to $\Omega(n^{-1/D})$. In the stochastic setting, when the budget $n$ is sufficiently large (i.e., $n \ge \max\{b_{\Gcal,1}\sigma^2, b_{\Gcal,3}\}$) but remains in a moderate range, the variance-dependent term $\Omega\bigl((\sigma^2/n)^{1/(D+2)}\bigr)$ may be dominated by the variance-independent term $\Omega(n^{-1/D})$. In this regime, the effective lower bound coincides with that of deterministic CSO. As $n$ increases further, the variance-dependent term, which decays at rate $n^{-1/(D+2)}$, eventually dominates the variance-independent term, causing the effective lower bound to transition to $\Omega\bigl((\sigma^2/n)^{1/(D+2)}\bigr)$. Thus, Theorem~\ref{th:lb_local_large_s} fully characterizes the variance dichotomy, capturing the initial deterministic-like behavior, the regime transition, and the asymptotic stochastic rate.

Moreover, the lower bound in Theorem~\ref{th:lb_local_large_s} is near-tight and therefore minimax optimal. In particular, \citet{bartlett2019simple} showed in Corollaries~9 and~11 that for CSO problems with $\alpha<\beta$ (so that $D>0$), their algorithm achieves convergence rates of $\tilde{O}\bigl((\sigma^2/n)^{1/(D+2)}\bigr)$ in the high-variance regime and $\tilde{O}(n^{-1/D})$ when the noise variance is negligible or zero. This yields an overall rate of
\[
\max\left\{		\tilde{O}\left({\sigma^2}/{n}\right)^{\frac{1}{D+2}},\; \tilde{O} (n^{-\frac{1}{D})} \right\}
,
\]
which matches our minimax lower bound up to logarithmic factors. Other existing methods attain partial matches to this result; for example, several tree search algorithms \citep{shang2019general,grill2015black,munos2014} achieve $\tilde{O}(n^{-1/(D+2)})$ in stochastic CSO. However, these methods typically require prior knowledge of an upper bound on $\sigma^2$, and their explicit dependence on the noise variance is not characterized.

Finally, since $D=d(1/\alpha-1/\beta)$, Theorem~\ref{th:lb_local_large_s} also reveals the pronounced impact of dimensionality on optimization error when $\alpha<\beta$, indicating that high-dimensional CSO is inherently more challenging. For example, when $d=1$, $\alpha=1$, and $\beta=2$, the lower bound becomes $\max\left\{\Omega\bigl((\sigma^2/n)^{2/5}\bigr),\; \Omega(n^{-2})\right\}$. In contrast, when the dimension increases to $d=10$, the rate deteriorates to $\max\left\{\Omega\bigl((\sigma^2/n)^{1/7}\bigr),\; \Omega(n^{-1/5})\right\}$, illustrating the severe curse of dimensionality inherent in this regime.

\subsection{Overall Lower Bound When $\alpha=\beta$}

Let $\Ccal_{\alpha,0}$ denote the class of CSO problems satisfying Assumptions~\ref{ass:unique}--\ref{ass:smo} for the given parameters with $\alpha=\beta$. By construction, both $\Gcal_{\alpha,\beta,1}$ with $\alpha = \beta$ and $\Ccal_{\alpha,2}$ are subsets of $\Ccal_{\alpha,0}$, which represents the full class of CSO problems with shape parameters satisfying $\alpha=\beta$. Analogous to Section \ref{subsec:over_lb_uneq}, the minimax lower bound over $\Ccal_{\alpha,0}$ must be no smaller than the corresponding bounds over its subsets. Consequently, Propositions~\ref{prop:lb_local_large} and~\ref{prop:lb_scv2} together yield the following overall minimax lower bound.

\begin{theorem}\label{th:lb_scv_large_s}
	For any adaptive algorithm satisfying Assumption~\ref{ass:adaptive_algo}, there exist constants $b_{\Ccal,1}$ and $b_{\Ccal,3}$ such that for all $n \ge \max\left\{b_{\Ccal,1}\sigma^2,\, b_{\Ccal,3}\right\}$,
	\begin{align*}
		\sup_{\CSO \in \Ccal_{\alpha,0}} 
		\Eb_{\CSO}\bigl[y(\x^*)-y(\hat{\x}_n^*)\bigr]
		\ge b_{\Ccal,6}\cdot
		\max\left\{
		\Bigl({\sigma^2}/{n}\Bigr)^{1/2},\,
		\exp(-b_{\Ccal,5}n)
		\right\},
	\end{align*}
	where $b_{\Ccal,6}$ is a positive constant independent of $n$ and $\sigma^2$.
\end{theorem}

The lower bound in Theorem~\ref{th:lb_scv_large_s} is again expressed as the maximum of two terms, reflecting the two regimes of the variance dichotomy when $\alpha=\beta$. In the deterministic setting ($\sigma^2=0$), the bound reduces to a purely exponential decay. In the stochastic setting, however, the dominant term depends on the simulation budget. When $n$ is sufficiently large (i.e., $n\ge \max\{b_{\Ccal,1}\sigma^2,\,b_{\Ccal,3}\}$) but remains in a moderate range, the variance-dependent term $\Omega\bigl((\sigma^2/n)^{1/2}\bigr)$ may be dominated by the exponential term, which therefore governs the effective lower bound. As $n$ increases further, the polynomial term $n^{-1/2}$ eventually dominates the exponential decay, leading to a transition in the effective lower bound to $\Omega\bigl((\sigma^2/n)^{1/2}\bigr)$. Thus, Theorem~\ref{th:lb_scv_large_s} provides a complete characterization of the variance dichotomy in the regime $\alpha=\beta$.

Moreover, the lower bound in Theorem~\ref{th:lb_scv_large_s} is near-tight in its order of magnitude. \citet{bartlett2019simple} showed that their algorithm achieves a convergence rate of $\tilde{O}\bigl((\sigma^2/n)^{1/2}\bigr)$ in the high-variance regime and a near-exponential rate of $O\left(\exp\bigl(-b_{\Ccal,7}n/(\log n)^2\bigr)\right)$ in the low-variance regime, where $b_{\Ccal,7}>0$ is a constant. Consequently, their method attains an overall rate $\max\big\{\tilde{O}\bigl((\sigma^2/n)^{1/2}\bigr),$ $ O\left(\exp\bigl(-b_{\Ccal,7}n/(\log n)^2\bigr)\right)\big\}$,  which nearly matches our minimax lower bound $\max\big\{\Omega\bigl((\sigma^2/n)^{1/2}\bigr),$ $ \Omega\bigl(\exp(-b_{\Ccal,5}n)\bigr)\big\}$ up to logarithmic factors. The constants governing the exponential rates, namely $b_{\Ccal,5}$ in our lower bound and $b_{\Ccal,7}$ in the algorithmic upper bound, do not coincide exactly; reconciling this gap remains an interesting direction for future research.

Other existing methods further corroborate this characterization by providing partial matches to the bound. In stochastic CSO, the effective lower bound inevitably transitions to the slower polynomial rate $\Omega(n^{-1/2})$ as the simulation budget grows. A variety of algorithms \citep{shang2019general,bubeck2018,wang2018optimization,grill2015black} asymptotically achieve the corresponding upper rate $\tilde{O}(n^{-1/2})$, confirming that this polynomial rate constitutes the fundamental precision limit in the presence of stochastic noise.

While exponential convergence is typically associated with deterministic optimization, our framework demonstrates that such rates are also attainable for CSO in low-variance regimes under suitable shape conditions. In particular, strongly concave objectives, which are widely studied in the literature \citep{hong2020finite,nemirovski2009robust}, often satisfy our smoothness assumptions corresponding to $\alpha=\beta$ in Theorem~\ref{th:lb_scv_large_s}. The theorem indicates that, although CSO algorithms tailored to stochastic strongly concave CSO can achieve $O(n^{-1/2})$ convergence \citep{shamir2013complexity,akhavan2024gradient,yu2024stochastic}, substantially faster exponential rates are theoretically possible when the noise level is sufficiently small.

This improvement, however, comes with an important trade-off. The general variance-independent lower bounds established here exhibit an exponential dependence on the dimension, reflecting a pronounced curse of dimensionality. In contrast, specialized results for strongly concave optimization that exploit gradient structure, such as the $O(d\,n^{-1/2})$ rates in \citet{shamir2013complexity,akhavan2024gradient,yu2024stochastic}, scale only linearly with $d$. Consequently, while our results reveal the potential for dramatically faster convergence in low-dimensional, low-noise CSO, algorithms tailored specifically to strong concavity remain more favorable in high-dimensional settings.

\section{Numerical Evidence and Additional Insights}\label{sec:num}

In this section, we use numerical experiments on synthetic CSO problems to illustrate the variance dichotomy and to obtain additional insights on the effect of variance. Let $\c$ denote a fixed solution. We consider two objective functions:
\begin{enumerate}
	\item Test function 1: $y_1(\x) = 1 - \| \x - \c \|_{2}^2 $.
	
	\item Test function 2: 
	\begin{align*}
		y_2(\x)= 1 - \left[\| \x - \c \|_{2} - \frac{ 1-\cos(4\pi/\| \x - \c \|_{2})}{2} \left(0.5 \| \x - \c \|_{2}^2 - \| \x - \c \|_{2}\right) \right].
	\end{align*}
\end{enumerate}
For both objective functions, the unique optimal solution is $\x^* = \c$. The function $y_1(\x)$ satisfies Assumption \ref{ass:smo} with $\alpha = \beta = 2$. In contrast, $y_2(\x)$ oscillates between two polynomial envelopes of orders one and two, and therefore satisfies Assumption \ref{ass:smo} with $\alpha = 1$ and $\beta = 2$. The third panel of Figure \ref{fig:intro} plots $y_2(\x)$ for the case $d = 1$.
In this experiment, the problem dimension and the optimal solution $\c$ are set to $d \in \{1,3,5\}$ and $(e^{-1},\dots,e^{-1})$, respectively. For each combination of problem dimension and objective function, we construct four synthetic CSO problems with sampling variances $\sigma^2 \in \{ (10^{-1})^2, (10^{-3})^2, (10^{-5})^2,0 \}$.

\begin{figure}[htbp]
	\centering    \includegraphics[width=1.0\textwidth]{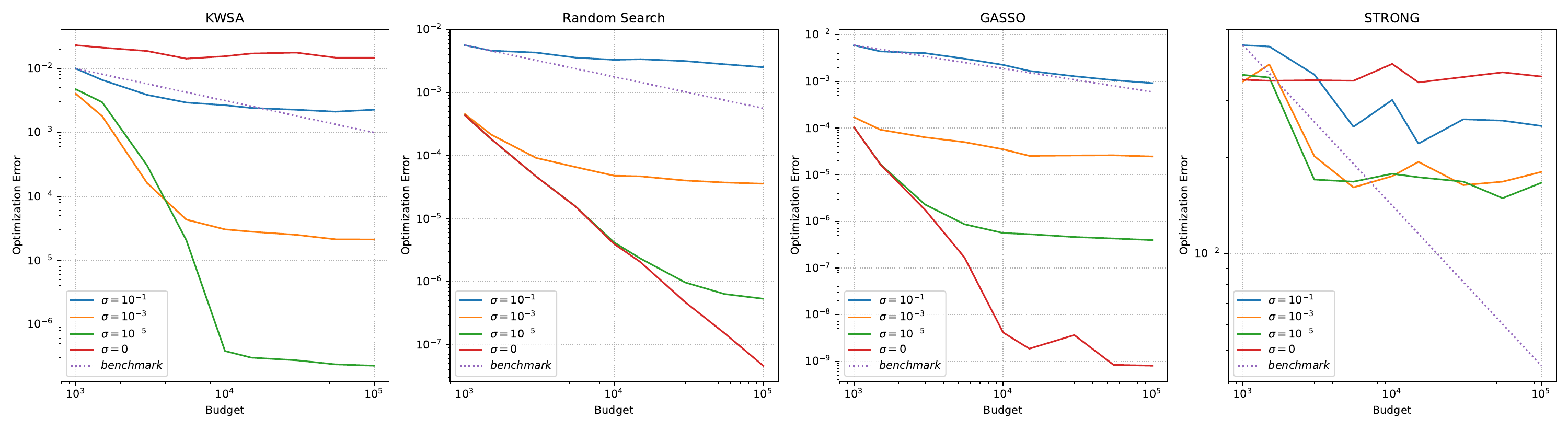}
	\caption{When the objective function is $y_1(\x)$ and problem dimension is one.}\label{fig:AlgTest11}
\end{figure}

\begin{figure}[htbp]
	\centering	\includegraphics[width=1.0\textwidth]{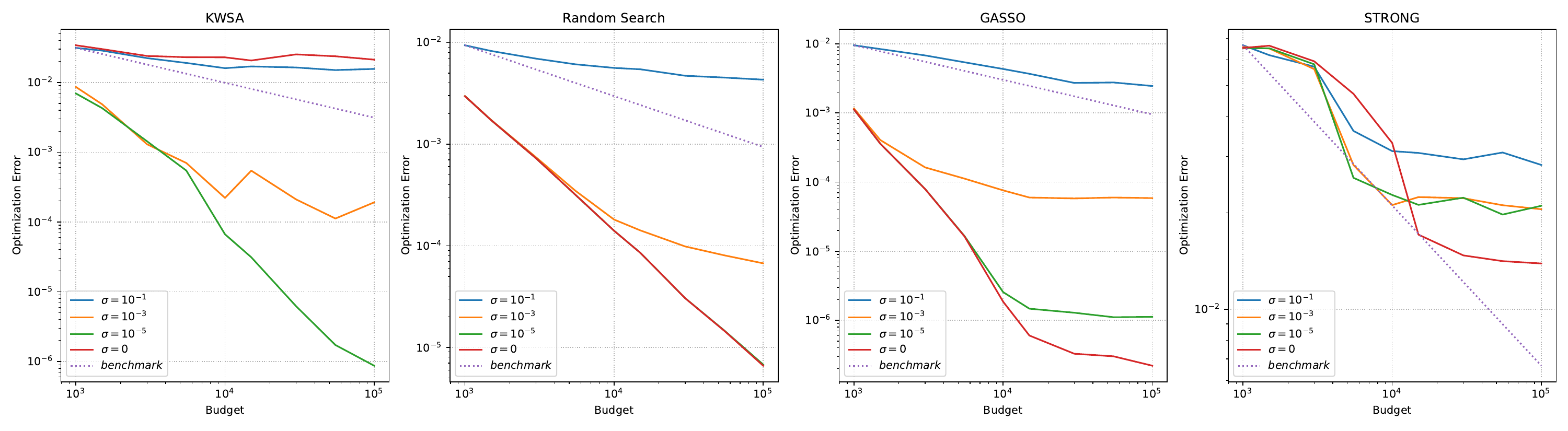}
	\caption{When the objective function is $y_2(\x)$ and problem dimension is one.}\label{fig:AlgTest21}
\end{figure}

We first examine the behavior of existing CSO algorithms under the variance dichotomy. The benchmark methods include the stochastic approximation algorithm KWSA \citep{kiefer1952stochastic}, two random search methods, Uniform Search \citep{chia2013limit} and GASSO \citep{zhou2017gradient}, and a model-based method, STRONG \citep{chang2013stochastic}. The implementations of these four methods are from \cite{Eckman_SimOpt}. For the one-dimensional problem with objective $y_1(\x)$, Figure \ref{fig:AlgTest11} reports the optimization errors estimated via 2500 macro-replications for the four synthetic CSO problems sharing the same objective but differing in variance. Each panel corresponds to one algorithm. For reference, we include the rate $r(n) = \hat{c} n^{-1/2}$ with appropriately chosen $\hat{c}$. All axes are on logarithmic scales. Figure \ref{fig:AlgTest21} presents the analogous results for $y_2(\x)$.

Two observations emerge from Figures \ref{fig:AlgTest11}-\ref{fig:AlgTest21}. First, the random search methods benefit from variance reduction. Although GASSO consistently outperforms Uniform Search in terms of error magnitude, both exhibit similar convergence patterns across variance levels. In particular, performance improves as $\sigma$ decreases, with the deterministic case ($\sigma = 0$) achieving the smallest errors. For small but positive variance ($\sigma = 10^{-3},10^{-5}$), both methods initially exhibit rapid error decay, followed by a transition to slower rates as the budget increases. Second, KWSA and STRONG do not exhibit improved performance in the deterministic setting. KWSA performs best when $\sigma$ is small but strictly positive, which can be attributed to its step-size mechanism: the step size depends on variance estimates, and when $\sigma = 0$, the resulting updates become overly conservative. In contrast, STRONG constructs surrogate models to guide sampling and is less sensitive to the noise level.

These experiments indicate that random-search methods that adapt to the noise level can perform superior empirically. However, CSO methods that explicitly account for the noise level remain relatively limited in the literature. The lower-bound analysis in this paper suggests that StroquOOL \citep{bartlett2019simple} is near-optimal under the variance dichotomy. StroquOOL employs a hierarchical partition of the search space, similar to the decomposition introduced in Section \ref{sec:var_ind_eq}, and adaptively allocates samples across refinement levels based on simple sample means. In the following experiments, we use StroquOOL as the representative method to illustrate the main findings of this paper.

\begin{figure}[htbp]
	\centering
	\caption{The optimization errors for synthetic CSO problems with objective function $y_1(\x)$.} \label{fig:AlgTest}
	\includegraphics[width=0.99\textwidth]{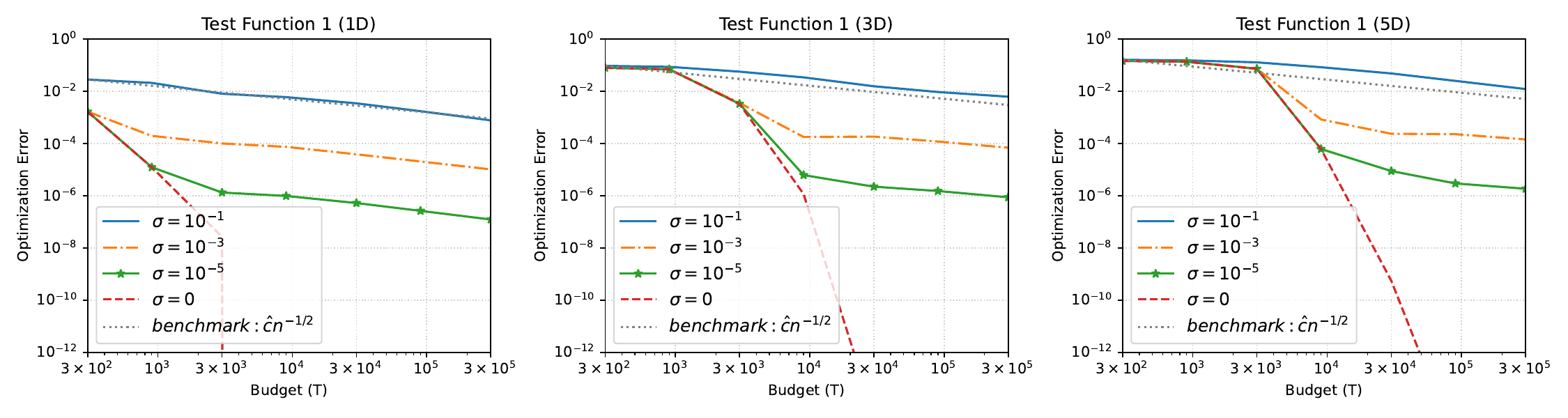}	
\end{figure}

We apply StroquOOL to CSO problems with budgets $n \in \{ 3\times10^2,10^3,3\times10^3,10^4,3\times10^4,10^5,3\times10^5\}$. Figure \ref{fig:AlgTest} contains three panels, each corresponding to a different problem dimension. Each panel reports the optimization errors estimated using 2500 macro-replications for the four CSO problems sharing the same objective function $y_1(\x)$. For comparison, we also plot the reference rate $r(n) = \hat{c} n^{-1/2}$, where $\hat{c}$ is chosen properly.

The variance dichotomy is clearly observed in each panel of Figure \ref{fig:AlgTest}. When the simulation budget $n$ is sufficiently large, the optimization error in the deterministic case decreases at a substantially faster rate than in the stochastic cases. Although smaller noise variances lead to uniformly smaller errors, the long-run convergence rates of the stochastic CSO problems are nearly identical. In particular, their error curves are approximately parallel to the benchmark line $\hat{c} n^{-1/2}$, indicating that variance primarily affects the constant factor rather than the asymptotic decay rate. A distinctive two-phase behavior emerges for low-noise stochastic problems. In the initial stage, the optimization errors decrease at a rate comparable to that of the deterministic problem. As the budget increases, however, the decay switches to the slower stochastic regime. The switching time depends on the noise level: smaller variances delay the onset of the stochastic regime. For example, in the three-dimensional setting (second panel), the switching budget is approximately $10^4$ when $\sigma = 10^{-5}$, compared with roughly $3 \times 10^3$ when $\sigma = 10^{-3}$.

\begin{figure}[htbp]
	\centering
	\caption{The normalized optimization errors for synthetic CSO problems with objective function $y_1(\x)$.} \label{fig:AlgTest1n}
	\includegraphics[width=0.99\textwidth]{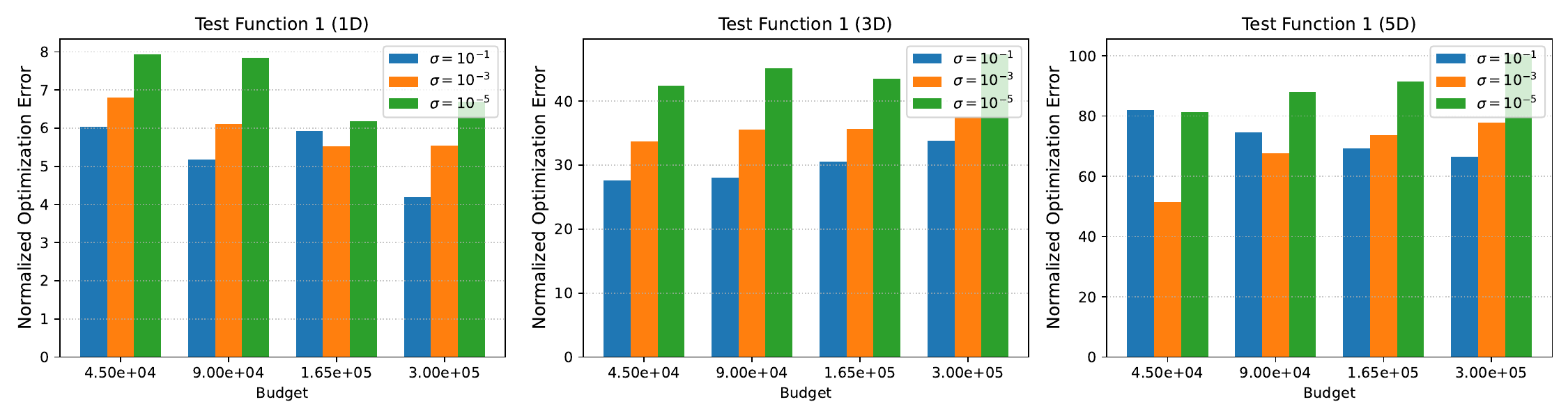}	
\end{figure}

The variance-dependent lower bounds indicate that the optimal convergence rates for stochastic CSO in the long run scale with $\sigma^2/n$, rather than $n$ alone. To further illustrate the long-run impact of $\sigma^2$ in stochastic CSO problems, Figure \ref{fig:AlgTest1n} reports the \emph{normalized} optimization error
\begin{align*}
\Eb[ y(\x^*) - y(\hat\x_n^*) ] / (\sigma^2/n)^{1/(D+2)}
\end{align*}
for budgets $n \in \{ 4.50\times10^4, 9.00\times10^4,1.65\times10^5,3.00\times10^5\}$, which are sufficiently large such that the convergence rates have switched to the variance-dependent regime. The normalized errors are of comparable magnitude across different values of $\sigma^2$ and $n$. This confirms that the convergence rate is jointly governed by the variance $\sigma^2$ and the budget $n$.

\begin{figure}[htbp]
	\centering
	\caption{The optimization errors for synthetic CSO problems with objective function $y_2(\x)$.} \label{fig:AlgTest2}
	\includegraphics[width=0.99\textwidth]{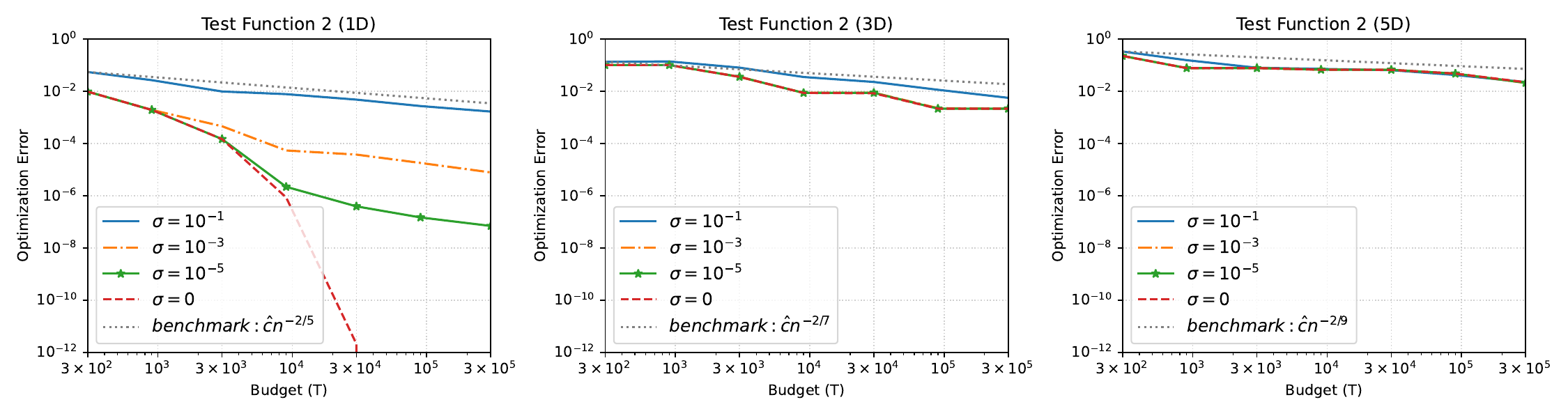}	
\end{figure}

We further apply StroquOOL to CSO problems with objective $y_2(\x)$ where $\alpha < \beta$. Figure \ref{fig:AlgTest2} presents results analogous to Figure \ref{fig:AlgTest}. The variance dichotomy remains evident in the first panel (one-dimensional setting), where the minimax lower bound takes the form $\max\{\Omega((\sigma^2/n)^{ 2/5 }), \Omega(n^{- 2 }) \}$. The convergence rates of the stochastic CSO problems closely follow the benchmark rate $n^{-2/5}$, consistent with the variance-dependent lower bound. In contrast, the deterministic CSO exhibits a rate faster than $n^{-2}$ (not plotted for clarity). This does not contradict the theory: minimax lower bounds characterize worst-case performance over a class and do not preclude faster rates for specific instances.

The dependence on dimension is also consistent with the theory. In Figure \ref{fig:AlgTest} with $\alpha = \beta$, multi-dimensional problems exhibit larger optimization errors despite eventually achieving the same asymptotic rate $n^{-1/2}$. 
The influence of dimension is more pronounced in Figure \ref{fig:AlgTest2}, where $\alpha < \beta$ and the minimax rates deteriorate with dimension. Specifically, the bounds become $\max\{\Omega((\sigma^2/n)^{ 2/7 }), \Omega(n^{- 2/3 }) \}$ when $d=3$ and $\max\{\Omega((\sigma^2/n)^{ 2/9 }),$ $\Omega(n^{- 2/5 }) \}$ when $d=5$. The exponents of the variance-dependent and variance-independent terms are close, which leads to similar convergence behavior for deterministic and stochastic problems in the second and third panels. In addition, in such cases, the switching time required for the variance-dependent term to dominate becomes large. Thus, when $\alpha < \beta$ and the dimension is moderately large, a substantial simulation budget is required to achieve variance-adaptive performance and outperform methods that do not adapt to the variance.

\section{Conclusion}\label{sec:con}

This paper studies the variance dichotomy in CSO through a minimax lower-bound analysis. Under a general smoothness framework, we show that the effective lower bound is given by the maximum of a variance-dependent term and a variance-independent term. This type of lower bound captures the switching behavior between deterministic and stochastic regimes. When the simulation budget is limited and the variance is not very large, the variance-independent bound dominates, so low-noise stochastic CSO has essentially the same problem complexity as deterministic CSO. As the budget increases, the variance-dependent bound of stochastic CSO eventually becomes dominant, and the convergence behavior switches to the slower stochastic regime. 

These findings refine the asymptotic view in the literature by showing that the effect of noise depends not only on whether the variance is positive, but also on the available simulation budget. In particular, when the budget is limited and the variance-independent term dominates, the budget is the primary driver of complexity. When the budget is sufficiently large and the variance-dependent term dominates, the complexity of stochastic CSO is governed jointly by the noise level and the budget, rather than by the budget alone. These observations provide a theoretical foundation for the development of more efficient CSO algorithms.

\bibliographystyle{poms}
\bibliography{reference}

\newpage
\appendix
{\noindent \LARGE \textbf{E-Companion to ``Understanding the Variance Dichotomy in Continuous Simulation Optimization: A Minimax Lower Bound Perspective"}}

\section{Near-Optimality Dimension Under Assumption \ref{ass:smo}}\label{sec:near_opt}

The definition of near-optimality dimension relies on a hierarchical partition of $\Xcal$. Let $h$ denote the level of a subregion in the hierarchical structure, $h=0,1,2,\dots$. At each level, the solution space is divided into mutually disjoint subregions $\{\Xcal_{h,i}\}$, $1 \le i \le I_h$, where $I_h$ is the total number of subregions at level $h$. The subregion at a higher level is a subset of one subregion at a lower level. As the level increases, $\Xcal$ is divided into smaller subregions. 

For each partition level, there exists one unique subregion, denoted by $\Xcal_{h}^*$, that contains the optimal solution. Suppose for any $\x \in \Xcal_h^*$ and any $h=0,1,2,\dots$, there exist constants $\nu > 0$ and $\rho \in (0,1)$ so that 
\begin{align}\label{ineq:lit_loc_smo}
	y(\x) \ge y(\x^*) - \nu \rho^{ h}.
\end{align}
Let $\Ncal_h(\epsilon)$ denote the number of subregions $\Xcal_{h,i}$ at level $h$ satisfying $\sup_{\x \in \Xcal_{h,i}} y(\x) \ge y(\x^*) - \epsilon$. The near-optimality dimension \citep{grill2015black} is defined as
\begin{align}\label{ineq:lit_near_dim}
d_{\text{opt}} = \inf \Bigl\{ d' \in \mathbb{R}^+: \exists C > 0, \forall h \ge 0, \Ncal_h(2\nu \rho^{ h}) \le C \rho^{-d' h} \Bigr\}.
\end{align}

Consider a regular hierarchical structure where the space $\Xcal = [0,1]^d$ is divided into $L^{dh}$ equal-size and mutually disjoint subregions at level $h$ where $L \ge 2$ is an integer. The side length of each subregion at level $h$ is $L^{-h}$. Then for each level $h=0,1,2,\dots$, we have 
\begin{align*}
\| \x - \x^* \|_{\infty} \le L^{-h}, \ \forall \x \in \Xcal_h^*. 
\end{align*}
By the upper bound $|y(\x^*) - y(\x)| \le  M \| \x^*-\x \|_{\infty}^{\alpha}$ of Assumption \ref{ass:smo}, \eqref{ineq:lit_loc_smo} holds with $\nu = M$ and $\rho = L^{-\alpha}$. By the lower bound $|y(\x^*) - y(\x)| \ge  \tilde{M} \| \x^*-\x \|_{\infty}^{\beta}$ of Assumption \ref{ass:smo}, if $y(\x) \ge y(\x^*) - 2\nu \rho^{ h}$, then
\begin{align}\label{ineq:lit_eff_reg}
\| \x^*-\x \|_{\infty} \le (2 \nu / \tilde{M})^{1/\beta} \rho^{h/\beta}
\end{align}
must hold. Since the side length of subregions at level $h$ is $L^{-h} = \rho^{h/\alpha}$, there are at most
\begin{align*}
\left( \left\lceil \frac{(2 \nu / \tilde{M})^{1/\beta} \rho^{h/\beta}}{\rho^{h/\alpha}} \right\rceil + 1 \right)^d \le \left(  (2 \nu / \tilde{M})^{1/\beta} \rho^{h/\beta-h/\alpha}  + 2 \right)^d
\end{align*}
subregions at level $h$ having a solution satisfying \eqref{ineq:lit_eff_reg}. Let $C' = \max\bigl\{ (2 \nu / \tilde{M})^{1/\beta},2 \bigr\}  $. We have $\Bigl(  (2 \nu / \tilde{M})^{1/\beta} \rho^{h/\beta-h/\alpha}  + 2 \Bigr)^d \le \bigl( 2C'\rho^{h/\beta-h/\alpha} \bigr)^d = ( 2C' )^d \rho^{-hd(1/\alpha-1/\beta)}$. Thus
\begin{align*}
\Ncal_h(2\nu \rho^{ h}) \le ( 2C' )^d \rho^{-hd(1/\alpha-1/\beta)}.
\end{align*}
By \eqref{ineq:lit_near_dim}, the near-optimality dimension $d_{\text{opt}}$ is less than $d(1/\alpha-1/\beta)$.

\section{Proofs in Section \ref{sec:lem}}\label{sec:proof_lem}

\subsection{Proof of Lemma \ref{lem:lb}}\label{subsec:proof_lemlb}

\textbf{Part I: deterministic sampling and recommendation methods.} First consider the algorithm under Assumption \ref{ass:adaptive_algo} whose sampling and recommendation method are deterministic. For such algorithm, the $(t+1)$-th simulated solution  $\x_{t+1}$ is a measurable function of $\x_1,Y(\x_1),\dots,\x_t,Y(\x_t)$, $t = 1,2,\dots, n-1$. We show $\x_{t+1}$ can also be expressed as a measurable function of $Y(\x_1),Y(\x_2),\dots,$ $Y(\x_t)$ (that is, the dependency on $\x_1,\dots,\x_t$ can be omitted) by induction:
\begin{itemize}
	\item For $t=1$, item 1 of Assumption \ref{ass:adaptive_algo} shows that $\x_{2}$ is a measurable function of $\x_1,Y(\x_1)$. Let $m\big(\x_1,Y(\x_1)\big)$ denote the function. Since $\x_1$ is deterministic, we can define $m'\big(Y(\x_1)\big) = m\big(\x_1,Y(\x_1)\big)$ such that $\x_2 = m'\big(Y(\x_1)\big)$ is a measurable function of $Y(\x_1)$.
	
	\item Suppose $\x_{s+1}$ is a measurable function of $Y(\x_1),Y(\x_2),\dots,Y(\x_{s})$ for $s=1,2,\dots,t-1$, which means that all of $\x_2$, $\x_3,\dots,\x_t$ are measurable functions of $Y(\x_1),Y(\x_2),\dots,Y(\x_{t-1})$. Combining this result with the fact that $\x_{t+1}$ is a measurable function of $\x_1,Y(\x_1),\dots,$ $\x_t,Y(\x_t)$ by Assumption \ref{ass:adaptive_algo}, we have that solution $\x_{t+1}$ is a measurable function of $\x_1,Y(\x_1),Y(\x_2),\dots,$ $Y(\x_t)$. Since $\x_1$ is deterministic, the dependency on $\x_1$ can be neglected.
\end{itemize}
Similarly, the estimated optimal solution $\hat{\x}_n^*$ is a measurable function of  $Y(\x_1),Y(\x_2),\dots,Y(\x_{n})$. Let $\Xi_n\big( Y(\x_1),Y(\x_2),\dots,Y(\x_{n}) \big)$ denote the measurable function such that $\hat{\x}_n^* = \Xi_n\big( Y(\x_1),Y(\x_2),$ $\dots,Y(\x_{n}) \big)$. 
The event $ \big\{\hat{\x}^*_n \in \Zcal \big\}$ equals to $\left\{ \big( Y(\x_1),Y(\x_2),\dots,Y(\x_n) \big) \in \Vcal \right\}$ where $\Vcal = \Xi_n^{-1}( \Zcal )$. Then
\begin{align*}
	\Pb_{2} \left( \hat{\x}^*_n \in \Zcal  \right) 
	=& \Eb_{2} \left[ \1 \Big( \big( Y(\x_1),Y(\x_2),\dots,Y(\x_n) \big) \in \Vcal \Big) \right]  \\
	=& \Eb_{1} \left[ \1 \Big( \big( Y(\x_1),Y(\x_2),\dots,Y(\x_n) \big) \in \Vcal \Big) \frac{f_{2} \big( Y(\x_1),Y(\x_2),\dots,Y(\x_n)\big) }{ f_{1}\big(  Y(\x_1),Y(\x_2),\dots,Y(\x_n) \big) } \right],
\end{align*}
where we let $f_{i}$ denote the density under the CSO Problem $i$, $i=1,2$.

For notation simplicity, let $L_{n} = \log f_{1}\big( Y(\x_1),\dots,Y(\x_n) \big) - \log f_{2}\big( Y(\x_1),\dots,Y(\x_n) \big) $ and $\mathcal{E} = \big\{ \big( Y(\x_1),Y(\x_2),\dots,Y(\x_n) \big) \in \Vcal \big\}$. Then,
\begin{align*}
	\Pb_{2} \left( \mathcal{E}  \right) =& \Eb_{1} \bigl[ \1 \left( \mathcal{E} \right) \exp (-L_{n}) \bigr]  \\
	\ge& \Eb_{1} \Bigl[ \1 \left( \mathcal{E} \right)   \exp \Big(-\Eb_{1} \big[L_{n} \big| \1 \left( \mathcal{E} \right)\big] \Big)   \Bigr] \\
    =& \Pb_{1} (\mathcal{E}) \exp \Big(-\Eb_{1} \big[L_{n} \big| \1 \left( \mathcal{E} \right) = 1\big]  \Big),
\end{align*}
which yields $\Eb_{1} \big[L_{n} \big| \1 \left( \mathcal{E} \right) = 1\big] \ge \log \frac{ \Pb_{1} (\mathcal{E}) }{ \Pb_{2} \left( \mathcal{E}  \right) }$. Similarly, $\Eb_{1} \big[L_{n} \big| \1 \left( \mathcal{E} \right) = 0\big] \ge \log \frac{ 1-\Pb_{1} (\mathcal{E}) }{ 1- \Pb_{2} \left( \mathcal{E}  \right) }$.
We have that
\begin{align*}
	\Eb_{1} [ L_{n} ] =& \Eb_{1} \big[L_{n} \big| \1 \left( \mathcal{E} \right) = 1\big] \Pb_{1} (\mathcal{E}) +  \Eb_{1} \big[L_{n} \big| \1 \left( \mathcal{E} \right) = 0\big] (1-\Pb_{1} (\mathcal{E}))  \\
	\ge& \Pb_{1} (\mathcal{E}) \log \frac{ \Pb_{1} (\mathcal{E}) }{ \Pb_{2} \left( \mathcal{E}  \right) }  + \bigl(1-\Pb_{1} (\mathcal{E})\bigr) \log \frac{ 1-\Pb_{1} (\mathcal{E}) }{ 1- \Pb_{2} \left( \mathcal{E}  \right) }   \\
	\ge& 2 \bigl( \Pb_{1} (\mathcal{E}) - \Pb_{2} (\mathcal{E}) \bigr)^2,
\end{align*}
where the last inequality holds by Pinsker's inequality. Then
\begin{align}\label{ineq:deter_lem_lb}
	\Pb_{2} (\mathcal{E}) \le \sqrt{\Eb_{1} [ L_{n} ]/2} + \Pb_{1} (\mathcal{E}).
\end{align}
By the chain rule for densities, we have 
\begin{align*}
	f_{1} \Big( Y(\x_1),Y(\x_2),\dots,Y(\x_n)\Big)  
	= f_1\big(Y(\x_1)\big) \cdot f_1\Big(Y(\x_2)\big|Y(\x_1)\Big) \cdot \dots \cdot f_1\Big(Y(\x_n)\big|Y(\x_1),\dots,Y(\x_{n-1})\Big) ,
\end{align*}
We have shown that $\x_{t}$ is a function of $Y(\x_1),Y(\x_2),\dots,Y(\x_{t-1})$. Since $Y(\x_{t}) = y(\x_t) + \varepsilon(\x_t)$ and by Assumption \ref{ass:subgauss}, conditional on $\x_t$, the noise term $\varepsilon(\x_t)$ is independent of everything else. It follows that $Y(\x_{t})$ is independent of $Y(\x_1),Y(\x_2),\dots,Y(\x_{t-1})$ conditional on $\x_t$.
Then 
\begin{align*}
f_{1} \big( Y(\x_{t}) \big| Y(\x_1),Y(\x_2),\dots,Y(\x_{t-1}) \big) = f_{1} \big( Y(\x_{t}) \big| \x_{t} \big)
\end{align*}
can be simplified as  
\begin{align*}
	\Eb_{1} [ L_{n} ]  
	=& \Eb_{1} \left[ \log \frac{f_{1} \big( Y(\x_1),Y(\x_2),\dots,Y(\x_n)\big)}{f_{2} \big( Y(\x_1),Y(\x_2),\dots,Y(\x_n)\big)}  \right]  \nonumber  \\
	=& \Eb_{1} \left[ \log \frac{ f_1\big(Y(\x_1)\big) \cdot f_1\big(Y(\x_2)|Y(\x_1)\big) \cdot \dots \cdot f_1\big(Y(\x_n)|Y(\x_1),\dots,Y(\x_{n-1})\big)  }{ f_2\big(Y(\x_1)\big) \cdot f_2\big(Y(\x_2)|Y(\x_1)\big) \cdot \dots \cdot f_2\big(Y(\x_n)|Y(\x_1),\dots,Y(\x_{n-1})\big) }  \right]  \nonumber  \\
	=& \Eb_{1} \left[ \log \frac{ f_1\big(Y(\x_1)\big) \cdot f_1\big(Y(\x_2)|\x_2\big) \cdot \dots \cdot f_1\big(Y(\x_n)|\x_n\big)  }{ f_2\big(Y(\x_1)\big) \cdot f_2\big(Y(\x_2)|\x_2\big) \cdot \dots \cdot f_2\big(Y(\x_n)|\x_n\big) }  \right]  
	= \Eb_{1} \left[ \sum_{t=1}^n \log \frac{  f_{1}\big(Y(\x_t)|\x_t\big) }{ f_{2}\big(Y(\x_t)|\x_t\big)}  \right],
\end{align*}
where the last equality holds because $\x_1$ is fixed.

Since $Y(\x_t)$ given $\x_t$ is normally distributed with mean $y(\x_t)$ and variance $\sigma^2$, the conditional expectation $\Eb_{1} \left[ \log   \Big( f_{1}\big(Y(\x_t)|\x_t\big) /  f_{2}\big(Y(\x_t)|\x_t\big) \Big) \Big| \x_t \right]$ is the KL divergence of $\Ncal\big( y_{1}(\x_t),\sigma^2 \big)$ from $\Ncal\big( y_{2}(\x_t),\sigma^2 \big)$ such that
\begin{align*}
	\Eb_{1} \left[ \log \frac{  f_{1}\big(Y(\x_t)|\x_t\big) }{  f_{2}\big(Y(\x_t)|\x_t\big)}  \right] = \Eb_{1} \left[ \Eb_{1} \left[ \log \frac{  f_{1}\big(Y(\x_t)|\x_t\big) }{  f_{2}\big(Y(\x_t)|\x_t\big)} \Big| \x_t\right] \right] = \Eb_{1} \left[ \frac{ \big(y_{1}(\x_t) - y_{2}(\x_t)\big)^2 }{2 \sigma^2  } \right].
\end{align*}
Thus, $\Eb_{1} [ L_{n} ] = \Eb_1 \left[ \sum_{t=1}^n  \big(y_{1}(\x_t) - y_{2}(\x_t)\big)^2 / (2 \sigma^2)   \right]$.

\textbf{Part II: stochastic sampling and recommendation methods.} Next, we consider the algorithm whose sampling and recommendation method are stochastic. For the algorithm satisfying this requirement, $\x_{t+1}$ is a deterministic function of $\x_1,Y(\x_1),\dots,\x_t,Y(\x_t),U_{t+1}$, $t=1,2,\dots, n-1$. We can show by induction that $\x_{t+1}$ can be expressed as a function of $U_1,Y(\x_1),\dots,U_{t},Y(\x_t),U_{t+1}$. That is, the role of $\x_1,\x_2,\dots,\x_t$ can be replaced by $U_1,U_2,\dots,U_t$. 
\begin{itemize}
	\item For $t=0$, $\x_1$ is a measurable function of $U_1$ by Assumption \ref{ass:adaptive_algo}.
	
	\item For $t=1$, by Assumption \ref{ass:adaptive_algo}, $\x_{2}$ is a measurable function of $\x_1,Y(\x_1),U_2$ and thus a function of $U_1,Y(\x_1),U_2$. 
	
	\item Suppose $\x_{s+1}$ is a function of $U_1,Y(\x_1),\dots,U_{s},Y(\x_{s}),U_{s+1}$ for $s=0,1,\dots,t-1$. Then all of $\x_1,\x_2,\dots,\x_t$ are functions of $U_1,Y(\x_1),\dots,U_{t-1},Y(\x_{t-1}),U_{t}$. Combining this result with the fact that $\x_{t+1}$ is a measurable function of $\x_1,Y(\x_1),\dots,\x_t,Y(\x_t),U_{t+1}$, solution $\x_{t+1}$ is a measurable function of $U_1,Y(\x_1),\dots,U_t,Y(\x_t),U_{t+1}$.
\end{itemize}
Similarly, since $\hat{\x}_{n}^*$ is a measurable function of $\x_1,Y(\x_1),\dots,\x_n,$ $Y(\x_n),U_{n}^*$, the estimated optimal solution  $\hat{\x}_{n}^*$ can also be expressed as a function of $U_1,Y(\x_1),\dots,U_{n},Y(\x_n),U_{n}^*$.

Let $f_{i}( u_1, y_1, \dots,u_n,y_n,u_n^* )$ denote the density of $\big( U_1,Y(\x_1),\dots,U_n,Y(\x_n),U_n^* \big) $ under the CSO Problem $i$, $i=1,2$. Then $\hat{\x}_n^* \in \Zcal$ if and only if 
\begin{align*}
\Big( U_1,Y(\x_1),\dots,U_{n},Y(\x_n),U_{n}^*\Big) \in \tilde{\Vcal}
\end{align*}
where $\tilde{\Vcal}$ is a subset in the value space. Notice that 1) $U_1,\dots,U_{n},U_n^*$ are independent random vectors and 2) conditional on the value of $U_1,\dots,U_{n},U_n^*$, the stochastic sampling and recommendation methods reduces to the deterministic method again. Then we can obtain by \eqref{ineq:deter_lem_lb} that
\begin{align*}
	\Pb_{2} \big(\tilde{\mathcal{E}} \big| U_1,\dots, U_n, U_n^*\big) \le \sqrt{\Eb_{1} \big[ \tilde{L}_{n} \big| U_1, \dots, U_n, U_n^* \big]/2} + \Pb_{1} \big(\tilde{\mathcal{E}} \big| U_1, \dots, U_n, U_n^*\big),
\end{align*}
where 
\begin{align*}
&\tilde{L}_{n} = \log \Big( f_{1}\Big( U_1, Y(\x_1), \dots,U_n,Y(\x_n),U_n^* \Big) / f_{2}\Big( U_1, Y(\x_1), \dots,U_n,Y(\x_n),U_n^* \Big) \Big), \\
&\tilde{\mathcal{E}} = \Big\{ \big( U_1,Y(\x_1),\dots,U_n,Y(\x_n),U_n^* \big) \in \tilde{\Vcal} \Big\}.
\end{align*}
Taking the expectation with respect to $U_1,\dots, U_n, U_n^*$, we obtain
\begin{align*}
	\Pb_{2} (\tilde{\mathcal{E}} ) \le \Eb \left[ \sqrt{\Eb_{1} \big[ \tilde{L}_{n} \big| U_1, \dots, U_n, U_n^* \big]/2} \right] + \Pb_{1} (\tilde{\mathcal{E}}),
\end{align*} 
which can be simplified by Jensen's inequality as $\Pb_{2} (\tilde{\mathcal{E}} ) \le  \sqrt{\Eb_{1} [ \tilde{L}_{n} ]/2}  + \Pb_{1} (\tilde{\mathcal{E}})$.
By the property of conditional probability, we have 
\begin{align}
	&f_{1} \Big(U_1, Y(\x_1),\dots,U_n,Y(\x_n),U_n^*\Big) \nonumber \\  
	=& f (U_1) \cdot f_1\Big(Y(\x_1)|U_1\Big) \cdot \dots \cdot f(U_n) \cdot f_1\Big(Y(\x_n)\big|U_1,Y(\x_1),\dots,Y(\x_{n-1}),U_n\Big) \cdot f(U_n^*) , \label{ineq:prob_decomp2}
\end{align}
where $f(U_{t})$ and $f(U_n^*)$ are the density of $U_{t}$, $t=1,2,\dots,n$, and $U_n^*$. Then $\Eb_{1} [ \tilde{L}_{n} ]$ can be similarly simplified as 
\begin{align*}
	\Eb_{1} [ \tilde{L}_{n} ] =& \Eb_{1} \left[ \log \frac{f_{1} \Big(U_1, Y(\x_1),\dots,U_n,Y(\x_n),U_n^*\Big)}{f_{2} \Big(U_1, Y(\x_1),\dots,U_n,Y(\x_n),U_n^*\Big)}  \right]  = \Eb_{1} \left[ \sum_{t=1}^n \frac{ \big(y_{1}(\x_t) - y_{2}(\x_t)\big)^2 }{2 \sigma^2 } \right].
\end{align*}
This lemma is proved.

\subsection{Proof of Lemma \ref{lem:lb2}}\label{subsec:proof_lemlb2}

\textbf{Part I: deterministic sampling and recommendation methods.} Consider the algorithm whose sampling and recommendation method are deterministic. 

For $t=1$, $\x_1$ is deterministic, hence independent of whether Problem 1 or Problem 2 is the true problem. Therefore, $\Pb_{1}\big(\tau(\Zcal) > 1\big) = \Pb_{2}\big(\tau(\Zcal) > 1\big)$. 

For $t = 2,3,\dots,n-1$, as shown in the proof of Lemma \ref{lem:lb}, solution $\x_{t+1}$ can be expressed as a measurable function of $Y(\x_1),Y(\x_2),\dots,Y(\x_t)$, $t = 1,2,\dots, n-1$. Then we have 
\begin{align*}
	\Pb_{2}  \big(\tau(\Zcal ) > t\big) = \Pb_{2}  \Big(\bigcap\nolimits_{s=1}^{t} \{ \x_{s} \notin \tilde\Xcal  \}\Big)  = \1 \Big( \big( Y(\x_1),Y(\x_2),\dots,Y(\x_{t-1}) \big) \in \bar{\Vcal} \Big) \cdot \1 (\x_1 \notin \tilde\Xcal ) 
\end{align*}
where $\bar{\Vcal}$ is the subset such that $\bigcap_{s=2}^{t} \{ \x_{s} \notin \tilde\Xcal  \}$ if and only if $\big( Y(\x_1),Y(\x_2),\dots,Y(\x_{t-1})\big) \in \bar{\Vcal}$. When $\1 \left( \big( Y(\x_1),Y(\x_2),\dots,Y(\x_{t-1}) \big) \in \bar{\Vcal} \right) \cdot \1 (\x_1 \notin \tilde\Xcal ) =1$, since none of $\x_{1},\dots,\x_{t-1}$ is in $\tilde\Xcal $, we have $y_{1}(\x_{s}) = y_{2}(\x_{s})$ and the distribution of $Y(\x_{s})$ conditional on $\x_{s}$ is same in the two CSO problems, $s=1,\dots,t-1$. Let $f_{i}$ denote the density in CSO Problem $i$, $i=1,2$. We have
\begin{align*}
	\Pb_{2}  \big(\tau(\Zcal ) > t \big)=&\Eb_1 \left[ \1 \left( \tau(\Zcal ) > t \right)  \cdot \frac{f_{2} \big( Y(\x_1),Y(\x_2),\dots,Y(\x_{t-1})\big) }{ f_{1}\big(  Y(\x_1),Y(\x_2),\dots,Y(\x_{t-1}) \big) } \right]\\
	=& \Eb_1 \left[ \1 \left( \tau(\Zcal ) > t \right)  \cdot \frac{f_2\big(Y(\x_1)\big) \cdot f_2\big(Y(\x_2)|\x_2\big) \cdot \dots \cdot f_2\big(Y(\x_{t-1})|\x_{t-1}\big) }{ f_1\big(Y(\x_1)\big) \cdot f_1\big(Y(\x_2)|\x_2\big) \cdot \dots \cdot f_1\big(Y(\x_{t-1})|\x_{t-1}\big)  }  \right] \nonumber  \\
	=& \Eb_1 \left[ \1 \left( \tau(\Zcal ) > t \right)  \right].
\end{align*}

\textbf{Part II: stochastic sampling and recommendation methods.} Next, we consider the algorithm whose sampling and recommendation method are stochastic. 

For $t=1$, the probability of $\x_1 \in \tilde\Xcal $ depends on an independent random vector $U_1$ by Assumption \ref{ass:adaptive_algo}, which means $\Pb_{1}\big(\tau(\Zcal ) > 1\big) = \Pb_{2}\big(\tau(\Zcal ) > 1\big)$. 

For $t = 2,3,\dots,n-1$, as shown in the proof of Lemma \ref{lem:lb}, $\x_1,\dots,\x_t$ are functions of $\big( U_1,Y(\x_1),\dots,$ $U_{t-1},Y(\x_{t-1}),U_{t}\big)$. We have
\begin{align*}
	\Pb_{2}  \big(\tau(\Zcal ) > t\big) = \Pb_{2}  \Big(\bigcap\nolimits_{s=1}^{t} \{ \x_{s} \notin \tilde\Xcal  \}\Big)  = \Pb_{2} \Big[ \big( U_1,Y(\x_1),\dots,U_{t-1},Y(\x_{t-1}),U_{t}\big) \in \check{\Vcal} \Big]
\end{align*} 
where $\check{\Vcal}$ is a subset such that $\bigcap_{s=1}^{t} \{ \x_{s} \notin \tilde\Xcal  \}$ if and only if $\big(  U_1,Y(\x_1),\dots,U_{t-1},Y(\x_{t-1}),U_{t}\big) \in \check{\Vcal}$. When $\1 \big( (  U_1,Y(\x_1),\dots,U_{t-1},Y(\x_{t-1}),U_{t}) \in \check{\Vcal} \big)=1$, since none of $\x_{s}$, $s=1,\dots,t-1$, is in $\tilde\Xcal $, we have $y_{1}(\x_{s}) = y_{2}(\x_{s})$ and the distribution of $Y(\x_{s})$ conditional on $\x_{s}$ is same in the two CSO problems, $s=1,\dots,t-1$. Then 
\begin{align*}
	\Pb_{2}  \big(\tau(\Zcal ) > t\big) =& \Eb_1 \left[ \1 \left( \tau(\Zcal ) > t \right) \frac{f_{2} \big( U_1,Y(\x_1),\dots,U_{t-1},Y(\x_{t-1}),U_{t}\big) }{ f_{1}\big(  U_1,Y(\x_1),\dots,U_{t-1},Y(\x_{t-1}),U_{t} \big) }  \right] \\
	=& \Eb_1 \left[ \1 \left( \tau(\Zcal ) > t \right) \frac{ f_1\big(Y(\x_1)|\x_1\big)  \cdot \dots \cdot f_1\big(Y(\x_{t-1})|\x_{t-1}\big)  }{ f_2\big(Y(\x_1)|\x_1\big)  \cdot \dots \cdot f_2\big(Y(\x_{t-1})|\x_{t-1}\big) }  \right] \nonumber  \\
	=& \Eb_1 \left[ \1 \left( \tau(\Zcal ) > t \right) \right]
\end{align*}
where the second equality holds by Equation \eqref{ineq:prob_decomp2}. 
This lemma is proved.

\subsection{Proof of Lemma \ref{lem:lb3}}\label{subsec:proof_lemlb3}
Let $\kappa_1 \in \mathop{\arg\min}_{\kappa = 1,\dots,n_1}  \Eb \left[ \sum_{t=1}^{n_2} \1 ( \x_t \in \Zcal_{\kappa} ) \right]$. We must have $ \Eb \left[ \sum_{t=1}^{n_2} \1 (\x_t \in \Zcal_{\kappa_1})  \right] \le \frac{1}{2} $. Otherwise, if $\Eb \left[ \sum_{t=1}^{n_2} \1 ( \x_t \in \Zcal_{\kappa} ) \right] > \frac{1}{2}$ for all $\kappa =1,\dots,n_1$, then 
\begin{align}\label{ineq:lb3_ass}
	\sum\nolimits_{\kappa = 1}^{n_1} \Eb \left[ \sum\nolimits_{t=1}^{n_2} \1 ( \x_t \in \Zcal_{\kappa} ) \right] >  n_1/2.
\end{align}
Meanwhile, since $\Zcal_{\kappa}$, $\kappa = 1,\dots,n_1$, are mutually disjoint, we have $\sum_{\kappa = 1}^{n_1} \1 ( \x_t \in \Zcal_{\kappa} ) \le 1$ for any $t=1,\dots,n_2$, which leads to $\sum_{\kappa = 1}^{n_1} \sum_{t=1}^{n_2} \1 ( \x_t \in \Zcal_{\kappa} ) \le n_2 \le n_1/2$. This is contradictory to \eqref{ineq:lb3_ass}.

We can further show $\Pb \left( \bigcup_{t=1}^{n_2} \{ \x_t \in \Zcal_{\kappa_1} \} \right) \le \frac{1}{2}$ because if $\Pb \left( \bigcup_{t=1}^{n_2} \{ \x_t \in \Zcal_{\kappa_1} \} \right) > \frac{1}{2}$, then
\begin{align*}
	\Eb \Big[ \sum\nolimits_{t=1}^{n_2} \1 ( \x_t \in \Zcal_{\kappa_1} ) \Big] 
	\ge& \Eb \Big[ \sum\nolimits_{t=1}^{n_2} \1 ( \x_t \in \Zcal_{\kappa_1} ) \Big| \1 \Big(\bigcup\nolimits_{t=1}^{n_2} \{ \x_t \in \Zcal_{\kappa_1} \} \Big) = 1  \Big] \cdot \Pb \Big( \bigcup\nolimits_{t=1}^{n_2} \{ \x_t \in \Zcal_{\kappa_1} \} \Big) \\
	\ge& \Pb \Big( \bigcup\nolimits_{t=1}^{n_2} \{ \x_t \in \Zcal_{\kappa_1} \} \Big)  
	> \frac{1}{2}.
\end{align*}
Notice that $\bigcup_{t=1}^{n_2} \{ \x_t \in \Zcal_{\kappa_1} \}$ is equivalent to $\big\{ \tau (\Zcal_{\kappa_1}) \le n_2 \big\}$. Thus
$\Pb \big( \tau (\Zcal_{\kappa_1}) \le n_2 \big) \le \frac{1}{2}$. This lemma is proved.

\section{Proofs in Section \ref{sec:var_dep}}

\subsection{Proof of Lemma \ref{lem:glob_fun1_prop}}\label{subsec:glob_fun1_prop}

Notice that by \eqref{ineq:lb_M}, since $
\| \x - \c \|_{\infty}^{\beta - \alpha} \le 1 \le M/(2^{4\beta+1} \tilde{M})  \le M/(2^{2\beta+1} \tilde{M}) 
$ for all $\x \in [0,1]^d$, we have 
\begin{align}\label{ineq:lb_Mprop}
	2^{2\beta} \tilde{M} \| \x - \c \|_{\infty}^{\beta} \le \frac{M}{2} \| \x - \c \|_{\infty}^{\alpha}, \ \forall \x \in [0,1]^d.
\end{align}
Also by \eqref{ineq:lb_M}, we have $\rho =   \big(M/(2^{2\beta+1} \tilde{M})\big)^{\frac{1}{\beta}} \ge (2^{2 \beta})^{\frac{1}{\beta}} = 4$. 

To reduce clutter, let $y(\x) = y_{\Gcal} (\x| \c,\c_0,\zeta)$ in this proof.
First, when $\| \x - \c \|_{\infty} \le \zeta$, it is obvious that $|y (\c) - y (\x)| = \frac{M}{2} \| \x - \c \|_{\infty}^\alpha \le M \| \x - \c \|_{\infty}^{\alpha}$ is satisfied. Moreover, since $\zeta \le 1$, we have $\| \x - \c \|_{\infty}^{\beta-\alpha} \le \zeta^{\beta-\alpha} \le 1  \le M/(2\tilde{M})$ and thus, $|y (\c) - y (\x)| = \frac{M}{2} \| \x - \c \|_{\infty}^\alpha \ge \tilde{M} \| \x - \c \|_{\infty}^{\beta}$.

Second, when $\| \x - \c \|_{\infty}   > \zeta $ and $\| \x - \c_{0} \|_{\infty}   \le \rho \zeta^{\alpha/\beta}$, we have $y (\c) - y (\x) = \frac{M}{2} \zeta^{\alpha} \le \frac{M}{2} \| \x - \c \|_{\infty}^{\alpha}$. Notice that 
\begin{align*}
\| \x - \c \|_{\infty}^{\beta} \le \Big(\| \x - \c_{0} \|_{\infty} + \| \c_{0} - \c \|_{\infty}\Big)^{\beta} \le \big(2\rho \zeta^{\alpha/\beta}\big)^{\beta} = (2\rho )^{\beta} \zeta^{\alpha}.
\end{align*}
Meanwhile, by definition, $\rho =   \big(M/(2^{2\beta+1} \tilde{M})\big)^{\frac{1}{\beta}} \le \frac{1}{2} \big(M/(2\tilde{M})\big)^{\frac{1}{\beta}}$. The above inequality yields
\begin{align*}
\zeta^{\alpha} \ge (2\rho )^{-\beta} \| \x - \c \|_{\infty}^{\beta} \ge (2 \tilde{M}/M) \| \x - \c \|_{\infty}^{\beta}
\end{align*}
such that $y (\c) - y (\x) = \frac{M}{2} \zeta^{\alpha} \ge \tilde{M} \| \x - \c \|_{\infty}^{\beta}$.

Third, when $\| \x - \c \|_{\infty}   > \zeta$ and $\| \x - \c_{0} \|_{\infty}   > \rho \zeta^{\alpha/\beta}$, we have $y (\c) - y (\x) = \frac{M}{2} \zeta^{\alpha} -  2^{\beta}\tilde{M} \rho^{\beta} \zeta^{\alpha} + 2^{\beta}\tilde{M} \| \x - \c_{0} \|_{\infty}^\beta$. We prove the upper bound by further considering three cases.
\begin{enumerate}
	\item[(1)] Suppose $\| \x - \c \|_{\infty} > \| \x - \c_{0} \|_{\infty}$. Notice that $\| \x - \c \|_{\infty} \ge \zeta$ such that $\| \x - \c \|_{\infty}^{\alpha} \ge \zeta^{\alpha}$. Meanwhile, $2^{\beta} \tilde{M} \| \x - \c \|_{\infty}^\beta \le \frac{M}{2} \| \x - \c \|_{\infty}^\alpha$ for $\x \in [0,1]^d$ by \eqref{ineq:lb_Mprop}. We have
	\begin{align*}
		y (\c) - y (\x) 
		\le \frac{M}{2} \zeta^{\alpha} + 2^{\beta}\tilde{M} \| \x - \c_{0} \|_{\infty}^\beta 
        \le  \frac{M}{2}  \| \x - \c \|_{\infty}^{\alpha} + 2^{\beta} \tilde{M} \| \x - \c \|_{\infty}^\beta 
        \le  M \| \x - \c \|_{\infty}^{\alpha}.
	\end{align*}
	
	\item[(2)] Suppose $\| \x - \c \|_{\infty} \le \| \x - \c_{0} \|_{\infty}$ and $\| \x - \c \|_{\infty} > \rho \zeta^{\alpha/\beta}$. By definition, $\| \c_{0} - \c \|_{\infty} \le \rho \zeta^{\alpha/\beta}$. Then
	\begin{align*}
		\frac{\| \x - \c_{0} \|_{\infty}}{\| \x - \c \|_{\infty}}  \le \frac{\| \x - \c \|_{\infty} + \| \c_{0} - \c \|_{\infty}}{\| \x - \c \|_{\infty}} \le 1+\frac{\rho \zeta^{\alpha/\beta}}{\rho \zeta^{\alpha/\beta}} = 2,
	\end{align*}
	which together with \eqref{ineq:lb_Mprop} yields $2^{\beta}\tilde{M} \| \x - \c_{0} \|_{\infty}^\beta \le 2^{2\beta}\tilde{M} \| \x - \c \|_{\infty}^\beta \le \frac{M}{2} \| \x - \c \|_{\infty}^\alpha$.
	Meanwhile, $\| \x - \c \|_{\infty} \ge \zeta$ such that $\| \x - \c \|_{\infty}^{\alpha} \ge \zeta^{\alpha}$. Similar to the first case, we have
	\begin{align*}
		&y (\c) - y (\x) \le \frac{M}{2} \zeta^{\alpha} + 2^{\beta}\tilde{M} \| \x - \c_{0} \|_{\infty}^\beta
		\le  \frac{M}{2} \| \x - \c \|_{\infty}^{\alpha} + 2^{2\beta}\tilde{M} \| \x - \c \|_{\infty}^\beta          
		\le   M \| \x - \c \|_{\infty}^\alpha.
	\end{align*}
	
	\item[(3)] Suppose $\| \x - \c \|_{\infty} \le \| \x - \c_{0} \|_{\infty}$ and $\| \x - \c \|_{\infty} \le \rho \zeta^{\alpha/\beta}$. Notice that $\| \x - \c_{0} \|_{\infty} \le \| \x - \c \|_{\infty} + \| \c - \c_{0} \|_{\infty} \le \rho \zeta^{\alpha/\beta} + \rho \zeta^{\alpha/\beta} = 2\rho \zeta^{\alpha/\beta}$.
	Meanwhile, \ $\| \x - \c \|_{\infty} \ge \zeta$ such that $\| \x - \c \|_{\infty}^{\alpha} \ge \zeta^{\alpha}$. We have
	\begin{align*}
		y (\c) - y (\x) \le& \frac{M}{2} \zeta^{\alpha} + 2^{\beta}\tilde{M} \| \x - \c_{0} \|_{\infty}^\beta 
		\le \frac{M}{2}  \zeta^{\alpha} + 2^{\beta}\tilde{M} (2\rho)^{\beta}  \zeta^{\alpha} 
        \le M \| \x - \c \|_{\infty}^{\alpha}.
	\end{align*}
\end{enumerate}
We prove the lower bound by further considering three cases. Notice that $\frac{M}{2} -  2^{\beta}\tilde{M} \rho^{\beta} = \frac{M}{2} - 2^{\beta}\tilde{M} \frac{M}{2^{2\beta+1} \tilde{M}} \ge 0 $, which implies $y (\c) - y (\x) = \frac{M}{2} \zeta^{\alpha} -  2^{\beta}\tilde{M} \rho^{\beta} \zeta^{\alpha} + 2^{\beta}\tilde{M} \| \x - \c_{0} \|_{\infty}^\beta 
\ge 2^{\beta}\tilde{M} \| \x - \c_{0} \|_{\infty}^\beta$.
\begin{enumerate}
	\item[(1)] If $\| \x - \c \|_{\infty} \le \| \x - \c_{0} \|_{\infty}$, we have
	\begin{align*}
		y (\c) - y (\x) 
		\ge& 2^{\beta}\tilde{M} \| \x - \c_{0} \|_{\infty}^\beta \ge 2^{\beta}\tilde{M} \| \x - \c \|_{\infty}^\beta \ge \tilde{M} \| \x - \c \|_{\infty}^\beta.
	\end{align*}
	
	\item[(2)] If $\| \x - \c \|_{\infty} > \| \x - \c_{0} \|_{\infty}$ and $\| \x - \c \|_{\infty} > 2\rho \zeta^{\alpha/\beta}$, then by noticing that $\| \c - \c_{0} \|_{\infty} \le \rho \zeta^{\alpha/\beta}$, we have $\| \x - \c_{0} \|_{\infty}/\| \x - \c \|_{\infty}  \ge 1 - \| \c - \c_{0} \|_{\infty}/\| \x - \c \|_{\infty} \ge 1-\frac{1}{2} = \frac{1}{2}$
	which means
	\begin{align*}
		y (\c) - y (\x) 
		\ge& 2^{\beta}\tilde{M} \| \x - \c_{0} \|_{\infty}^\beta \ge \frac{2^{\beta}\tilde{M}}{2^\beta} \| \x - \c \|_{\infty}^\beta = \tilde{M} \| \x - \c \|_{\infty}^\beta.
	\end{align*}
	
	\item[(3)] If $\| \x - \c \|_{\infty} > \| \x - \c_{0} \|_{\infty}$ and $\| \x - \c \|_{\infty} \le 2\rho \zeta^{\alpha/\beta}$, then by noticing that $\| \x - \c_{0} \|_{\infty} > \rho \zeta^{\alpha/\beta}$ in all three case, we have $\| \x - \c_{0} \|_{\infty} / \| \x - \c \|_{\infty}  \ge \rho \zeta^{\alpha/\beta}/(2\rho \zeta^{\alpha/\beta}) = \frac{1}{2}$
	which means
	\begin{align*}
		y (\c) - y (\x) 
		\ge& 2^{\beta}\tilde{M} \| \x - \c_{0} \|_{\infty}^\beta \ge\frac{ 2^{\beta}\tilde{M} }{2^\beta} \| \x - \c \|_{\infty}^\beta = \tilde{M} \| \x - \c \|_{\infty}^\beta.
	\end{align*}
\end{enumerate}
This lemma is proved.

\subsection{Proof of Proposition \ref{prop:lb_local_large}}\label{subsec:proof_prop:lb_local_large}

For notational simplicity, the definition of notations introduced in this sub-subsection are valid within Section \ref{subsec:proof_prop:lb_local_large} only. Define constant $\nu = \big(   \rho^d \eta \sigma^2 \big)^{-1/(\alpha(D + 2))}$ where $\eta =  1/\big(2^{d-\alpha (D+2)-1 } 9M^2\big)$. The solution space $[0,1]^d$ is partitioned by a rectangular grid. Given the total budget $n$, let $\psi_1 = \big\lceil \nu n^{1/(\alpha(D + 2))} \big\rceil$ such that the $d$-dimensional region $[0,1]^d$ is partitioned into $\psi_1^d$ equal-size subregions. Let each $d$-dimensional subregion be indexed by the $d$-dimensional vector $\bkappa = (\kappa_1,\kappa_2,\dots,\kappa_d)$ with $\kappa_{l}=1,2,\dots,\psi_1$ and  $l=1,2,\dots,d$. 
Let $\c_{\bkappa} = \big( (2\kappa_1-1) \zeta_{1}, (2\kappa_2-1) \zeta_{1}, \dots, (2\kappa_d-1) \zeta_{1} \big)$ denote the center of the subregion indexed by $\bkappa$ where we let $\zeta_{1} = (2\psi_1)^{-1}$. The subregion indexed by $\bkappa$ is $\Zcal_{\bkappa} = \big\{\x \in [0,1]^d:  0 \le \| \x - \c_{\bkappa} \|_{\infty} < \zeta_{1} \big\}$.

Let $\bkappa_{0} = \big(\lceil \psi_1/2 \rceil,\lceil \psi_1/2 \rceil,\dots,\lceil \psi_1/2 \rceil\big)$. We consider the subregions in set $\mathcal{K}_{1}$ defined as
\begin{align*}
	\mathcal{K}_{1} = \big\{ \Zcal_{\bkappa}: \| \c_{\bkappa} - \c_{\bkappa_{0}} \|_{\infty} \le \rho \zeta_{1}^{\alpha / \beta} - \zeta_1 \big\}.
\end{align*}
Since $\| \c_{\bkappa_1} - \c_{\bkappa_2} \|_{\infty} \ge 2\zeta_{1}$ for $\bkappa_1 \ne \bkappa_2$, the number of subregions in $\mathcal{K}_{1}$ is greater than $\big(\big\lfloor \rho \zeta_{1}^{\alpha / \beta-1} -1 \big\rfloor \big)^d$. 

For each subregion $\Zcal_{\bkappa}$ in $\Kcal_{1}$, let $y_{\bkappa}(\x) = y_{\Gcal} (\x | \c_{\bkappa},\c_{\bkappa_{0}},\zeta_{1})$ where $y_{\Gcal} (\x | \c,\c_0,\zeta)$ is defined in \eqref{eq:glob_lb_fund}. Define the problem $\CSO_{\bkappa}$ for subregion $\Zcal_{\bkappa}$: the simulation observation at $\x$ is $Y(\x) = y_{\bkappa}(\x) + \varepsilon(\x)$
where $\varepsilon(\x)$ has the distribution $N(0,\sigma^2)$; we want to maximize $y_{\bkappa}(\x)$ by taking simulation observations on $y_{\bkappa}(\x)$ with budget $n$. Let $\Gcal_{\alpha,\beta,1} = \big\{ \CSO_{\bkappa}: \  \Zcal_{\bkappa} \in \Kcal_{1} \big\}$ denote the set of CSO problems $\CSO_{\bkappa}$ whose corresponding subregion is in $\Kcal_{1}$.
The set $\Gcal_{\alpha,\beta,1}$ contains more than $\big(\big\lfloor \rho \zeta_{1}^{\alpha / \beta-1} - 1 \big\rfloor \big)^d$ CSO problems.  

When $n$ large enough such that $ n \ge 4^{\alpha(D + 2)} \rho^d \eta \sigma^2$, we have 
\begin{align}\label{ineq:glob_lb_phi1_lb1}
	\psi_1 \ge  \nu n^{\frac{1}{\alpha(D + 2)}} = \big(   \rho^d \eta \sigma^2 \big)^{-\frac{1}{\alpha(D + 2)}} n^{\frac{1}{\alpha(D + 2)}} \ge 4.
\end{align} 
Since 
\begin{equation}\label{ineq:psi1bd}
\begin{aligned}
&\big(2\lceil \psi_1/2 \rceil-1\big) \zeta_{1} = \big(2\lceil \psi_1/2 \rceil-1\big) (2\psi_1)^{-1} \le ( \psi_1 + 1) (2\psi_1)^{-1} \le \frac{5}{8}, \\
&\big(2\lceil \psi_1/2 \rceil-1\big) \zeta_{1} \ge ( \psi_1 - 1) (2\psi_1)^{-1} \ge \frac{3}{8},
\end{aligned}
\end{equation}
each element of vector $\c_{\bkappa_{0}}$ is smaller than $\frac{5}{8}$ and larger than $\frac{3}{8}$. Meanwhile, when $n$ is large enough such that $ n \ge (\rho 2^{2-\alpha / \beta})^{\beta(D+2)} \rho^d \eta \sigma^2$, we have
\begin{align*}
\rho \zeta_{1}^{\frac{\alpha}{\beta}} \le \rho \big(2 \nu n^{\frac{1}{\alpha(D + 2)}}\big)^{-\frac{\alpha}{\beta}} = \rho 2^{-\frac{\alpha}{\beta}} \big(   \rho^d \eta \sigma^2 \big)^{\frac{1}{\beta(D+2)}} n^{-\frac{1}{\beta(D+2)}} \le \frac{1}{4}. 
\end{align*}
Let $b_{\Gcal,1} = \max\big\{ 4^{\alpha(D + 2)} , (\rho 2^{2-\alpha / \beta})^{\beta(D+2)} \big\} \rho^d \eta$. When $n \ge b_{\Gcal,1} \sigma^2$, we have
$\big\{ \x: \| \x - \c_{\bkappa_{0}} \|_{\infty} \le \rho \zeta_{1}^{\alpha / \beta-1} \big\} \subset [0,1]^d$. By Lemma \ref{lem:glob_fun1_prop}, the objective function of each CSO problem in $\Gcal_{\alpha,\beta,1}$ is well-defined and satisfies Assumption \ref{ass:smo} with $\alpha \le \beta$.

For notation simplicity, let $y_{\mathbf{0}}(\x) = y_{\Bcal}(\x|\c_{\bkappa_0},\zeta_1)$ where $y_{\Bcal}(\x|\c_{\bkappa_0},\zeta)$ is defined in \eqref{eq:glob_bench}. Define a benchmark problem $\CSO_{\mathbf{0}}$: the simulation observation at $\x$ is $Y(\x) = y_{\mathbf{0}}(\x_t) + \varepsilon(\x)$
where $\varepsilon(\x)$ has the distribution $N(0,\sigma^2)$; we want to maximize $y_{\mathbf{0}}(\x)$ by taking simulation observations on $y_{\mathbf{0}}(\x)$ with budget $n$.

Let $ \Eb_{\mathbf{0}} \left[ \sum_{t=1}^n \1 [ \x_t \in \Zcal_{\bkappa} ] \right]$ and $\Pb_{\mathbf{0}} \left[ \hat{\x}^*_n \in \Zcal_{\bkappa} \right]$ denote the expected number of simulated solutions that fall in $\Zcal_{\bkappa}$ and the probability that the estimated optimal point $\hat{\x}^*_n$ is in $\Zcal_{\bkappa}$ when the CSO problem is $\CSO_{\mathbf{0}}$.
Let $\bkappa_1$ and $\bkappa_2$ denote the cells having first and second smallest value of $\Eb_{\mathbf{0}} \left[ \sum_{t=1}^n \1 [ \x_t \in \Zcal_{\bkappa} ] \right]$, $ \bkappa \in \Kcal_1$. That is, 
\begin{align*}
	\bkappa_1 \in \mathop{\arg\min}_{\bkappa \in \Kcal_1} \Eb_{\mathbf{0}} \left[ \sum\nolimits_{t=1}^n \1 [ \x_t \in \Zcal_{\bkappa} ] \right], \ \bkappa_2 \in \mathop{\arg\min}_{\bkappa \in \Kcal_1\setminus\{\bkappa_1\}} \Eb_{\mathbf{0}} \left[ \sum\nolimits_{t=1}^n \1 [ \x_t \in \Zcal_{\bkappa} ] \right].
\end{align*} 
Then there must have a $\bar{\bkappa} \in \{ \bkappa_1, \bkappa_2 \}$  satisfying $\Pb_{\mathbf{0}}  (\hat{\x}^*_n \in \Zcal_{\bar{\bkappa}})  \le \frac{1}{2}$ and $\Eb_{\mathbf{0}} \big[ \sum_{t=1}^n \1 [ \x_t \in \Zcal_{\bar{\bkappa}} ] \big] \le 2n / | \Kcal_{1} | $. Notice that $\rho = \big(M/(2^{2\beta+1} \tilde{M})\big)^{1/\beta} \ge 4$ by Equation \eqref{ineq:lb_M}. Meanwhile, by Equation \eqref{ineq:glob_lb_phi1_lb1}, $\psi_1 \ge 4$. Then $\frac{\rho}{2} (2 \psi_1)^{1-\alpha / \beta} \ge 2$. Since $| \Kcal_{1} | \ge \big(\big\lfloor \rho \zeta_{1}^{\alpha / \beta-1} -1\big\rfloor \big)^d$, we have 
\begin{align}\label{ineq:Kcal_1_lb}
	| \Kcal_{1} | \ge \big(\big\lfloor \rho \zeta_{1}^{\alpha / \beta-1}-1 \big\rfloor \big)^d &\ge \big( \rho (2 \psi_1)^{1-\alpha / \beta} - 2 \big)^d  \ge \Big(\frac{ \rho}{2} (2 \psi_1)^{1 - \alpha / \beta} \Big)^d  =  \Big(\frac{\rho}{2}\Big)^d \zeta_{1}^{-\alpha D}.
\end{align}

By Lemma \ref{lem:lb},
we have that
\begin{align}\label{ineq:lb_glob_ex_measure}
	\Pb_{\bar{\bkappa}} (\hat{\x}^*_n \in \Zcal_{\bar{\bkappa}}) \le \sqrt{\Eb_{\mathbf{0}} [ L_{\mathbf{0}, \bar{\bkappa}, n} ]/2} + \Pb_{\mathbf{0}} (\hat{\x}^*_n \in \Zcal_{\bar{\bkappa}}).
\end{align}
Notice that
\begin{align*}
	\Eb_{\mathbf{0}} [ L_{\mathbf{0}, \bar{\bkappa}, n} ] = \Eb_{\mathbf{0}} \left[ \sum\nolimits_{t=1}^n \frac{ (y_{\mathbf{0}}(\x_t) - y_{\bar{\bkappa}}(\x_t))^2 }{2 \sigma^2 }  \right] 
	\le& \Eb_{\mathbf{0}} \left[ \sum\nolimits_{t=1}^n \1 [ \x_t \in \Zcal_{\bar{\bkappa}} ] \right] \max_{\x \in \Zcal_{\bar{\bkappa}}} \frac{ (  y_{\mathbf{0}}(\x) - y_{\bar{\bkappa}}(\x) )^2}{2\sigma^2} \\
	\le& \frac{2n}{ | \Kcal_{1} | } \frac{M^2}{8\sigma^2} \zeta^{2\alpha}_{1} \\
	\le& \frac{2^dn M^2}{4\rho^d\sigma^2 }  \zeta_{1}^{\alpha (D+2)}
\end{align*}
where the first inequality holds because $y_{\mathbf{0}}(\x) \ne y_{\bar{\bkappa}}(\x)$ only when $\x \in \Zcal_{\bar{\bkappa}}$, the second inequality holds because $\Eb_{\mathbf{0}} \left[ \sum_{t=1}^n \1 [ \x_t \in \Zcal_{\bar{\bkappa}} ] \right] \le 2n / | \Kcal_{1} | $, and the last inequality holds by \eqref{ineq:Kcal_1_lb}. Since $\zeta_{1} \le \big(2\nu n^{1/(\alpha(D+2))}\big)^{-1}$ and $\nu = (   \rho^d \eta \sigma^2 )^{-1/(\alpha(D+2))}$ where $\eta =  1/\big(2^{d-\alpha (D+2)-1 } 9M^2\big)$, we have
\begin{align*}
	\Eb_{\mathbf{0}} [ L_{\mathbf{0}, \bar{\bkappa}, n} ] \le& \frac{2^d nM^2  }{4 \rho^d\sigma^2} \left(2\nu n^{\frac{1}{\alpha(D+2)}}\right)^{-\alpha(D+2)}  
	= \frac{2^{d-\alpha (D+2) }  M^2  }{ 4 \rho^d\sigma^2} \nu^{-\alpha(D+2)}  
	= \frac{2^{d-\alpha (D+2) }  M^2  }{ 4 } \eta = \frac{1}{18}.
\end{align*}
Plugging the above inequality into \eqref{ineq:lb_glob_ex_measure}, we obtain
\begin{align*}
	\Pb_{\bar{\bkappa}} (\hat{\x}^*_n \in \Zcal_{\bar{\bkappa}}) \le& \sqrt{\Eb_{\mathbf{0}} [ L_{\mathbf{0}, \bar{\bkappa}, n} ]/2} + \Pb_{\mathbf{0}} (\hat{\x}^*_n \in \Zcal_{\bar{\bkappa}}) \le \frac{1}{6} + \frac{1}{2} \le \frac{2}{3}.
\end{align*}
Meanwhile, since $n \ge b_{\Gcal,1} \sigma^2 \ge \rho^d \eta \sigma^2 = \nu^{-\alpha(D + 2)}$, which leads to 
\begin{align*}
	\zeta_{1} = \big(2 \big\lceil \nu n^{\frac{1}{\alpha(D + 2)}} \big\rceil\big)^{-1} \ge \left(2  \left(\nu n^{\frac{1}{\alpha(D + 2)}}+1\right) \right)^{-1} \ge \left(4 \nu n^{\frac{1}{\alpha(D + 2)}} \right)^{-1}.
\end{align*}
If $\hat{\x}^*_n \notin \Zcal_{\bar{\bkappa}}$, we have $y_{\bar{\bkappa}}(\hat{\x}_n^*) \le 0$ such that
\begin{align*}
	\max_{\x \in \Xcal} y_{\bar{\bkappa}}(\x) - y_{\bar{\bkappa}}(\hat{\x}_n^*) \ge \frac{M}{2} \zeta^{\alpha}_{1}  \ge \frac{M}{2} \left(4 \nu n^{\frac{1}{\alpha(D + 2)}} \right)^{-\alpha} = \frac{ M \nu^{-\alpha}}{4^{\alpha} 2} n^{-\frac{1}{D + 2}}
\end{align*}
Thus,
\begin{align*}
	\Eb_{\bar{\bkappa}} \left[ \max\nolimits_{\x} y_{\bar{\bkappa}}(\x) - y_{\bar{\bkappa}}(\hat{\x}^*_n) \right]  
	& \ge \Pb_{\bar{\bkappa}} (\hat{\x}^*_n \notin \Zcal_{\bar{\bkappa}}) \Eb_{\bar{\bkappa}} \big[ \max\nolimits_{\x} y_{\bar{\bkappa}}(\x) - y_{\bar{\bkappa}}(\hat{\x}^*_n) \big| \hat{\x}^*_n \notin \Zcal_{\bar{\bkappa}}\big]  \\
	&\ge \left(1 - \frac{2}{3}\right) \frac{ M \nu^{-\alpha}}{4^{\alpha} 2} n^{-\frac{ 1 }{D + 2}} = \frac{  M (   \rho^d \eta  )^{\frac{1}{D + 2}} }{ 4^{\alpha} 6 }  \sigma^{\frac{2}{D + 2}}   n^{-\frac{1 }{D + 2}}.
\end{align*}
Notice that $\rho^d \ge 4^d$ by \eqref{ineq:lb_M}. Then
\begin{align}\label{ineq:glob_bg2_lb}
(   \rho^d \eta  )^{\frac{1}{D + 2}} \ge \big(4^d (2^{d-\alpha (D+2)-1 } 9M^2)^{-1}\big)^{\frac{1}{D + 2}} = 2^{\frac{d+1}{D + 2}+\alpha} ( 9M^2)^{-\frac{1}{D + 2}}.
\end{align}
Proposition \ref{prop:lb_local_large} holds when $b_{\Gcal,2} = 2^{\frac{d+1}{D + 2}+\alpha} 9^{-\frac{1}{D + 2}} M^{\frac{D}{D + 2}} / (4^\alpha 6) $.

\section{Proofs in Section \ref{sec:var_ind_uneq}}

\subsection{Proof of Proposition \ref{prop:lb_local_large2}}\label{subsec:lb_glob_fun2}
To show Proposition \ref{prop:lb_local_large2}, we introduce another set of CSO problems (denoted by $\Gcal_{\alpha,\beta,2}$) whose parameter values of the objective function are independent of variance. The CSO problem in $\Gcal_{\alpha,\beta,2}$ is almost the same as those in $\Gcal_{\alpha,\beta,1}$ except that the values of some parameters such as $\nu,\psi,\zeta$ are changed. Thus, we omit some explanations in this part. Notice that to distinguish the CSO problems between this part and Section \ref{subsec:proof_prop:lb_local_large}, some notations of this part are subscript by 2. 

Define constant $\nu = (   \rho \eta )^{-d/(\alpha D)}$ where $\eta = \frac{1}{4}$. Given the total budget $n$, let $\psi_2 = \big\lceil \nu n^{1/(\alpha D)} \big\rceil$ such that the $d$-dimensional region $[0,1]^d$ is partitioned into $\psi_2^d$ equal-size subregions. Let each $d$-dimensional subregion be indexed by vector $\bkappa = (\kappa_1,\kappa_2,\dots,\kappa_d)$ with $\kappa_{l}=1,2,\dots,\psi_2$ and  $l=1,2,\dots,d$. Let $\zeta_{2} = (2\psi_2)^{-1}$ and $\c_{\bkappa} = \big( (2\kappa_1-1) \zeta_{2}, \dots, (2\kappa_d-1) \zeta_{2} \big)$. The subregion indexed by $\bkappa$ is $\Zcal_{\bkappa} = \big\{\x \in [0,1]^d:  0 \le \| \x - \c_{\bkappa} \|_{\infty} < \zeta_{2} \big\}$.

Let $\bkappa_{0} = \big(\lceil \psi_2/2 \rceil,\lceil \psi_2/2 \rceil,\dots,\lceil \psi_2/2 \rceil\big)$.  Consider subregions in set $\mathcal{K}_{2}$ defined as
\begin{align*}
	\mathcal{K}_{2} = \big\{ \Zcal_{\bkappa}: \| \c_{\bkappa} - \c_{\bkappa_0} \|_{\infty} \le \rho \zeta_{2}^{\alpha / \beta - 1} - \zeta_{2} \big\}.
\end{align*}
For $\bkappa \in \Kcal_{2}$, let $y_{\bkappa}(\x) = y_{\Gcal} (\x | \c_{\bkappa},\c_{\bkappa_{0}},\zeta_{2})$ where $y_{\Gcal} (\x|\c,\c_0,\zeta)$ is defined in \eqref{eq:glob_lb_fund}. 
Define the CSO problem $\CSO_{\bkappa}$ for subregion $\Zcal_{\bkappa}$: the simulation observation at $\x$ is $Y(\x) = y_{\bkappa}(\x) + \varepsilon(\x)$
where $\varepsilon(\x)$ has the distribution $N(0,\sigma^2)$; we want to maximize $ y_{\bkappa}(\x)$ by taking simulation observations on $y_{\bkappa}(\x)$ with budget $n$. Let $\Gcal_{\alpha,\beta,2} = \big\{ \CSO_{\bkappa}: \  \bkappa \in \Kcal_{2} \big\}$ denote the set of CSO problems $\CSO_{\bkappa}$ indexed by $\bkappa \in \Kcal_{2}$.

When $n$ large enough such that $ n \ge 4^{\alpha D} (\rho \eta)^d$, we have 
\begin{align*}
	\psi_2 \ge  \nu n^{\frac{1}{\alpha D}} = (   \rho \eta )^{-\frac{d}{ \alpha D}} n^{\frac{1}{\alpha D}} \ge 4.
\end{align*} 
Analogous to \eqref{ineq:psi1bd}, we can show 
each element of vector $\c_{\bkappa_{0}}$ is smaller than $5/8$ and larger than $3/8$. Meanwhile, when $n$ large enough such that $ n \ge \big(\rho 2^{2-\alpha / \beta}\big)^{\beta D} (\rho \eta)^d$, we have
\begin{align*}
	\rho \zeta_{2}^{\frac{\alpha}{\beta}} \le \rho \big(2 \nu n^{\frac{1}{\alpha D}}\big)^{-\frac{\alpha}{\beta}} = \rho 2^{-\frac{\alpha}{\beta}} (   \rho \eta )^{\frac{d}{\beta D}} n^{-\frac{1}{\beta D}} \le \rho 2^{-\frac{\alpha}{\beta}} (   \rho \eta )^{\frac{d}{\beta D}} \big(\rho 2^{2-\alpha / \beta}\big)^{-1} (\rho \eta)^{-\frac{d}{\beta D}}  = \frac{1}{4}.
\end{align*}
Let $b_{\Gcal,3} = \max\big\{4^{\alpha D}, \big(\rho 2^{2-\alpha / \beta}\big)^{\beta D}\big\} (\rho \eta)^d$ such that $ n \ge b_{\Gcal,3}$. 
Then $\big\{ \x: \| \x - \c_{\bkappa_{0}} \|_{\infty} \le \rho \zeta_2^{\alpha / \beta} \big\} \subset [0,1]^d$ and by Lemma \ref{lem:glob_fun1_prop}, the objective function of each CSO problem in $\Gcal_{\alpha,\beta,2}$ is well-defined and satisfies Assumption \ref{ass:smo} with $\alpha < \beta$.

Let $ \Pb_{\bkappa} \big( \bigcup_{t=1}^n \{ \x_t \in \Zcal_{\bkappa} \} \big)$ and $ \Pb_{\bkappa} \left( \hat{\x}^*_n \in \Zcal_{\bkappa} \right)$ denote the probabilities when the algorithm is used to solve CSO problem $\CSO_{\bkappa}$. Since $\hat{\x}^*_n \in \{\x_1,\dots,\x_n\}$ by Assumption \ref{ass:adaptive_algo}, 
\begin{align}\label{ineq:lb_local2_prob2}
	\Pb_{\bkappa} \left( \hat{\x}^*_n \in \Zcal_{\bkappa} \right) \le \Pb_{\bkappa} \Big( \bigcup\nolimits_{t=1}^n \{ \x_t \in \Zcal_{\bkappa} \} \Big).
\end{align}

For notation simplicity, let $y_{\mathbf{0}}(\x) = y_{\Bcal}(\x|\c_{\bkappa_0},\zeta_2)$ where $y_{\Bcal}(\x|\c_{\bkappa_0},\zeta)$ is defined in \eqref{eq:glob_bench}.  
Define a benchmark problem $\CSO_{\mathbf{0}}$: the simulation observation at $\x$ is
$
Y(\x) = y_{\mathbf{0}}(\x) + \varepsilon(\x)
$
where $\varepsilon(\x)$ has the distribution $N(0,\sigma^2)$; we want to maximize $ y_{\mathbf{0}}(\x)$ by taking simulation observations on $y_{\mathbf{0}}(\x)$ with budget $n$.

Since $\zeta_{2}^{\alpha / \beta - 1} = \zeta_{2}^{- \alpha D / d} \ge \big(\nu n^{1/(\alpha D)} \big)^{\alpha D / d} \ge (\rho \eta)^{-1} n^{1/d} = 4 \rho^{-1} n^{1/d}$, we have
\begin{align*}
	| \Kcal_{2} | &\ge \big(\big\lfloor \rho \zeta_{2}^{\alpha / \beta - 1} - 1 \big\rfloor \big)^d \ge \big( \rho \zeta_{2}^{\alpha / \beta - 1} - 2 \big)^d \ge (  4   n^{1/d} - 2 )^d .
\end{align*}
Since $4   n^{1/d} \ge 4$, we have  $4   n^{1/d} - 2 \ge 2 n^{1/d}$ such that $| \Kcal_{2} |  \ge 2^d n$. By Lemma \ref{lem:lb3}, there must exist a subregion $\Zcal_{\tilde{\bkappa}} \in \Kcal_{2}$ such that $\Pb_{\mathbf{0}} \big( \tau(\Zcal_{\tilde{\bkappa}}) \le n \big) \le \frac{1}{2}$. Notice that $y_{\mathbf{0}}(\x) \ne y_{\tilde{\bkappa}} (\x )$ only when $\x \in \Zcal_{\tilde{\bkappa}}$. Combining this property with Lemma \ref{lem:lb2}, we have
\begin{align*}
	&\Pb_{\tilde{\bkappa}} \left( \bigcup\nolimits_{t=1}^n \{ \x_t \in \Zcal_{\tilde{\bkappa}} \} \right) = \Pb_{\tilde{\bkappa}} \big( \tau(\Zcal_{\tilde{\bkappa}}) \le n \big) 
	=  \Pb_{\mathbf{0}} \big( \tau(\Zcal_{\tilde{\bkappa}}) \le n \big) \le \frac{1}{2}.
\end{align*}
By \eqref{ineq:lb_local2_prob2}, we have $\Pb_{\tilde{\bkappa}} \left( \hat{\x}^*_n \in \Zcal_{\tilde{\bkappa}} \right) \le 1/2$.

Since $n \ge (   \rho \eta )^{d}$ such that $\nu n^{1/(\alpha D)} = (   \rho \eta )^{-d/(\alpha D)} n^{1/(\alpha D)} \ge 1$, which leads to $\psi_2 = \big\lceil \nu n^{1/(\alpha D)} \big\rceil \le 2 \nu n^{1/(\alpha D)} $ and thus $ (2\psi_2)^{\alpha} \le 4^{\alpha} \nu^{\alpha} n^{1/D}$. If $\hat{\x}^*_n \notin \Zcal_{\tilde{\bkappa}}$ happens, we have 
\begin{align*}
	&\max_{\x \in \Xcal} y_{\tilde{\bkappa}} (\x ) - y_{\tilde{\bkappa}}(\hat{\x}^*_n)  
	\ge y_{\tilde{\bkappa}} (\c_{\tilde{\bkappa}} ) = \frac{M}{2} \zeta_{2}^{\alpha} = \frac{M}{2(2\psi_2)^{\alpha}} \ge \frac{M}{ 2 \cdot 4^{\alpha} \nu^{\alpha}  } n^{-\frac{1}{D}}.
\end{align*}
Thus,
\begin{align*}
	\Eb_{\tilde{\bkappa}} \big[ \max\nolimits_{\x} y_{\tilde{\bkappa}}(\x) - y_{\tilde{\bkappa}}(\hat{\x}^*_n) \big] 
	\ge& \Pb_{\tilde{\bkappa}} (\hat{\x}^*_n \notin \Zcal_{\tilde{\bkappa}}) \Eb_{\tilde{\bkappa}} \big[ \max\nolimits_{\x} y_{\tilde{\bkappa}}(\x) - y_{\tilde{\bkappa}}(\hat{\x}^*_n) \big| \hat{\x}^*_n \notin \Zcal_{\tilde{\bkappa}}\big]  \\
    \ge&   \frac{M}{ 4^{\alpha+1} \nu^{\alpha}  } n^{-\frac{1}{D}}.
\end{align*}

\section{Proofs in Section \ref{sec:var_ind_eq}}

\subsection{Details of Space Partition}\label{subsec:detail_partition}
First, we provide a detailed introduction to the space partition. Let $\Delta=\lfloor 3^d/2 \rfloor$ and $\bar{a} = \lceil 4(n+1)/\Delta \rceil+2$ denote the highest level of partition. Let $n_{a} = 5^{a}$ and $\gamma_{a} = (2n_{a})^{-1} = 5^{-a}/2$, $a=0,1,2,\dots,\bar{a}$. At each level $a$, the continuous region $[0,1]^d$ is divided into $n_{a}^d$ equal-size subregions. Let each subregion at level $a$ be indexed by the vector $\bkappa_a = (\kappa_{1,a},\dots,\kappa_{d,a})$.
Let $\c_{\bkappa_a,a} = \big( (2\kappa_{1,a}-1) \gamma_{a}, \dots, (2\kappa_{d,a}-1) \gamma_{a} \big)$. The subregion $\Zcal_{\bkappa_a,a}$ indexed by $\bkappa_a$ at level $a$ is $\big\{\x \in [0,1]^d:  0 \le \| \x - \c_{\bkappa_a,a} \|_{\infty} < \gamma_{a} \big\}$.   

The subregions across different levels have a nested structure, which we have introduced in the main paper.
Below, we provide a mathematical definition of the set $\Ucal_{a}$ of subregions at each level $a$ for which functions will be defined. At level $a=0$, there is only one subregion $\Zcal_{\mathbf{1},0} = (0,1)^d$. Let $\Ucal_{0} = \{ \Zcal_{\mathbf{1},0} \}$ denote the set of the subregion at level $a=0$. At level $a=1,2,\dots,\bar{a}$, define 
\begin{align}\label{eq:lb_scv_udef}
	\Ucal_a = \big\{ \Zcal_{\bkappa_{a},a}: \Zcal_{\Pcal(\bkappa_a),a-1}  \in \Ucal_{a-1} \text{ and } \big\| \c_{\bkappa_a,a} - \c_{\Pcal(\bkappa_a),a-1} \big\|_{\infty} \le 2\gamma_{a} \big\}.
\end{align}
If subregion $\Zcal_{\bkappa_{a},a}$ is in $\Ucal_a$, then by definition, its parent subregion must be in $\Ucal_{a-1}$ and the distance between the centers of $\Zcal_{\bkappa_{a},a}$ and its parent subregion $\Zcal_{\Pcal(\bkappa_a),a-1}$ should be less than $2\gamma_a$. 

For each subregion in $\Ucal_{a}$, $a=0,1,2,\dots,\bar{a}$, we have defined the individual function in Equations \eqref{eq:scv_fun20}-\eqref{eq:scv_fun3}. Based on the individual functions, we have defined the objective functions $y_{\Ccal}(\x|\bkappa_{a},a)$ in \eqref{eq:lb_scv_benchfun} and \eqref{eq:lb_scv_fun3}. Define the CSO problem $\CSO_{\bkappa_a,a}$ for each subregion $\Zcal_{\bkappa_a,a} \in \Ucal_{a}$, $a=0,1,2,\dots,\bar{a}$: the simulation observation at $\x$ is $Y(\x) = y_{\Ccal}(\x|\bkappa_{a},a) + \varepsilon(\x)$
where $\varepsilon(\x)$ has the distribution $N(0,\sigma^2)$; we want to maximize $y_{\Ccal}(\x|\bkappa_{a},a)$ by taking simulation observations on $y_{\Ccal}(\x|\bkappa_{a},a)$ with budget $n$. 

Let $\Ccal_{2} = \big\{ \CSO_{\bkappa_{\bar{a}},\bar{a}}: \  \Zcal_{\bkappa_{\bar{a}},\bar{a}} \in \Ucal_{\bar{a}} \big\}$ denote the set of CSO problems whose corresponding subregions are in $\Ucal_{\bar{a}}$. In the following, we provide a proof for Lemma \ref{lem:scv_fun1_prop} to show that the objective function of each CSO in $\Ccal_{2}$ satisfies Assumption \ref{ass:smo} with $\alpha = \beta$. Then we provide a proof for Proposition \ref{prop:lb_scv2} to establish the variance-independent lower bound.

\subsection{Proof of Lemma \ref{lem:scv_fun1_prop}}\label{subsec:lem:scv_fun1_prop}

By Equation \eqref{ineq:lb_M3}, we have
\begin{align}\label{ineq:lb_M33}
	M/\hat{M} \ge (1-5^{-\alpha})^{-1} (5^{\alpha}-3^{\alpha}+1)  2^{\alpha}.
\end{align}
Consider any subregion $\Zcal_{\bkappa_{\bar{a},\bar{a}}}$ in $\Ucal_{\bar{a}}$. Let $y_{\CP,a} (\x|\bkappa_{\bar{a}},\bar{a}) = y_{\mathcal{D}} (\x|\bkappa_{\bar{a}},\bar{a}) + \sum_{\ell=1}^{a} y_{\CS} \big(\x|\Pcal^{\ell}(\bkappa_{\bar{a}}),\bar{a}-\ell\big)$, $a = 1,2,\dots,\bar{a}$.  In this proof, we first consider function $y_{\mathcal{D}} (\x | \bkappa_{\bar{a}},\bar{a})$ for the subregion $\Zcal_{\bkappa_{\bar{a}},\bar{a}}$ at the highest level and then analyze functions $y_{\CP,a} (\x|\bkappa_{\bar{a}},\bar{a})$, $a = 1,2,\dots,\bar{a}$, by induction. Since $ y_{\Ccal}(\x|\bkappa_{\bar{a}},\bar{a}) = y_{\CP,\bar{a}} (\x|\bkappa_{\bar{a}},\bar{a})$, the analysis of $y_{\CP,\bar{a}} (\x|\bkappa_{\bar{a}},\bar{a})$ applies to $y_{\Ccal}(\x|\bkappa_{\bar{a}},\bar{a})$ immediately.

Notice that $M \ge \hat{M}$ by Equation \eqref{ineq:lb_M33} and $\hat{M} / \tilde{M} \ge 40^{\alpha} / 5^{\alpha} = 8^{\alpha} \ge 1$.
It is straightforward to show for $\x$ in $\Zcal_{\bkappa_{\bar{a}},\bar{a}} = \big\{\x: \| \x - \c_{\bkappa_{\bar{a}},\bar{a}} \|_{\infty}   \le \gamma_{\bar{a}}\big\}$ that
\begin{align}\label{ineq:yDcal_bd}
	\tilde{M} \| \x - \c_{\bkappa_{\bar{a}},\bar{a}} \|_{\infty}^\alpha \le \big| y_{\mathcal{D}} (\c_{\bkappa_{\bar{a}},\bar{a}}|\bkappa_{\bar{a}},\bar{a}) - y_{\mathcal{D}} (\x|\bkappa_{\bar{a}},\bar{a}) \big| = \hat{M} \| \x - \c_{\bkappa_{\bar{a}},\bar{a}} \|_{\infty}^\alpha \le M \| \x - \c_{\bkappa_{\bar{a}},\bar{a}} \|_{\infty}^\alpha.
\end{align}
Consider function $y_{\CP,1} (\x | \bkappa_{\bar{a}},\bar{a}) = y_{\mathcal{D}} (\x | \bkappa_{\bar{a}},\bar{a}) + y_{\CS} \big(\x | \Pcal(\bkappa_{\bar{a}}),\bar{a}-1\big)$.  Since $\Zcal_{\bkappa_{\bar{a},\bar{a}}} \in \Ucal_{\bar{a}}$, we have 
\begin{align}\label{eq:scv_lb_cdist}
	\big\| \c_{\bkappa_{\bar{a}},\bar{a}} - \c_{\Pcal(\bkappa_{\bar{a}}),\bar{a}-1} \big\|_{\infty} \le 2 \gamma_{\bar{a}}.
\end{align}
We show by discussing two cases that if $\big\| \x - \c_{\Pcal(\bkappa_{\bar{a}}),\bar{a}-1} \big\|_{\infty}   \le 5\gamma_{\bar{a}} = \gamma_{\bar{a}-1}$, then
\begin{align}\label{ineq:scv_lb_bd}
	\tilde{M} \| \x - \c_{\bkappa_{\bar{a}},\bar{a}} \|_{\infty}^\alpha \le \big| y_{\CP,1} (\c_{\bkappa_{\bar{a}},\bar{a}} | \bkappa_{\bar{a}},\bar{a}) - y_{\CP,1} (\x|\bkappa_{\bar{a}},\bar{a}) \big| \le M \| \x - \c_{\bkappa_{\bar{a}},\bar{a}} \|_{\infty}^\alpha.
\end{align}

\textbf{Case 1.} Suppose $\| \x - \c_{\bkappa_{\bar{a}},\bar{a}} \|_{\infty}   \le \gamma_{\bar{a}}$, which together with \eqref{eq:scv_lb_cdist} implies $ \| \x - \c_{\Pcal(\bkappa_{\bar{a}}),\bar{a}-1} \|_{\infty}   \le 3\gamma_{\bar{a}}$. Then $y_{\CS} \big(\x | \Pcal(\bkappa_{\bar{a}}),\bar{a}-1\big) = \hat{M} (5^{\alpha}-3^{\alpha}) \gamma_{\bar{a}}^{\alpha}$ by  \eqref{eq:scv_fun2}, which yields 
\begin{align*}
	y_{\CP,1} (\c_{\bkappa_{\bar{a}},\bar{a}}|\bkappa_{\bar{a}},\bar{a}) &= y_{\mathcal{D}} (\c_{\bkappa_{\bar{a}},\bar{a}}|\bkappa_{\bar{a}},\bar{a}) +   \hat{M} (5^{\alpha}-3^{\alpha}) \gamma_{\bar{a}}^{\alpha}, \\
	y_{\CP,1} (\x|\bkappa_{\bar{a}},\bar{a}) &= y_{\mathcal{D}} (\x|\bkappa_{\bar{a}},\bar{a}) + \hat{M} (5^{\alpha}-3^{\alpha}) \gamma_{\bar{a}}^{\alpha}.
\end{align*}
Then $\big| y_{\CP,1} (\c_{\bkappa_{\bar{a}},\bar{a}}|\bkappa_{\bar{a}},\bar{a}) - y_{\CP,1} (\x|\bkappa_{\bar{a}},\bar{a}) \big| = \big| y_{\mathcal{D}} (\c_{\bkappa_{\bar{a}},\bar{a}}|\bkappa_{\bar{a}},\bar{a}) - y_{\mathcal{D}} (\x|\bkappa_{\bar{a}},\bar{a}) \big|$ such that \eqref{ineq:scv_lb_bd} is satisfied by \eqref{ineq:yDcal_bd}.

\textbf{Case 2.} Suppose $\| \x - \c_{\bkappa_{\bar{a}},\bar{a}} \|_{\infty}   > \gamma_{\bar{a}}$ and $\big\| \x - \c_{\Pcal(\bkappa_{\bar{a}}),\bar{a}-1} \big\|_{\infty}   \le 5 \gamma_{\bar{a}}$, which together with \eqref{eq:scv_lb_cdist} implies $ \| \x - \c_{\bkappa_{\bar{a}},\bar{a}} \|_{\infty}   \le 7\gamma_{\bar{a}} $. For $\x$ in this case, we have $y_{\mathcal{D}} (\x|\bkappa_{\bar{a}},\bar{a})=0$ such that $y_{\CP,1} (\x|\bkappa_{\bar{a}},\bar{a}) = y_{\CS} \big(\x|\Pcal(\bkappa_{\bar{a}}),\bar{a}-1\big)$, which yields
\begin{align*}
	&y_{\CP,1} (\c_{\bkappa_{\bar{a}},\bar{a}}| \bkappa_{\bar{a}},\bar{a}) - y_{\CP,1} (\x|\bkappa_{\bar{a}},\bar{a}) \\
	=& \hat{M} \gamma_{\bar{a}}^{\alpha} + \hat{M} (5^{\alpha}-3^{\alpha}) \gamma_{\bar{a}}^{\alpha} - y_{\CS} \big(\x|\Pcal(\bkappa_{\bar{a}}),\bar{a}-1\big)  \\
	=& \left\{ \begin{array}{ll}
		\hat{M} \gamma_{\bar{a}}^{\alpha} , &  \big\| \x - \c_{\Pcal(\bkappa_{\bar{a}}),\bar{a}-1} \big\|_{\infty}   \le 3\gamma_{\bar{a}} \\
		\hat{M} \gamma_{\bar{a}}^{\alpha} - \hat{M} 3^{\alpha} \gamma_{\bar{a}}^{\alpha} + \hat{M} \big\| \x - \c_{\Pcal(\bkappa_{\bar{a}}),\bar{a}-1} \big\|_{\infty}^\alpha,     &  3\gamma_{\bar{a}} < \big\| \x - \c_{\Pcal(\bkappa_{\bar{a}}),\bar{a}-1} \big\|_{\infty}   \le 5\gamma_{\bar{a}}. 
	\end{array} \right. 
\end{align*}
Then $\hat{M} \gamma_{\bar{a}}^{\alpha} \le | y_{\CP,1} (\c_{\bkappa_{\bar{a}},\bar{a}}|\bkappa_{\bar{a}},\bar{a}) - y_{\CP,1} (\x|\bkappa_{\bar{a}},\bar{a}) | \le \hat{M} (5^{\alpha} - 3^{\alpha} + 1) \gamma_{\bar{a}}^{\alpha}$.
Notice that $\hat{M} = \frac{40^{\alpha}}{5^{\alpha}-3^{\alpha}} \tilde{M}$ such that $\hat{M} \ge 40^{\alpha} \tilde{M} /5^{\alpha}  \ge 7^{\alpha} \tilde{M}$. Moreover, $ \| \x - \c_{\bkappa_{\bar{a}},\bar{a}} \|_{\infty}   \le 7\gamma_{\bar{a}} $. We have
\begin{align*}
	\hat{M} \gamma_{\bar{a}}^{\alpha} \ge 7^{\alpha} \tilde{M}  \gamma_{\bar{a}}^{\alpha} \ge \tilde{M} \| \x - \c_{\bkappa_{\bar{a}},\bar{a}} \|_{\infty}^{\alpha},
\end{align*}
which means the lower bound in \eqref{ineq:scv_lb_bd} holds. Meanwhile, since $\| \x - \c_{\bkappa_{\bar{a}},\bar{a}} \|_{\infty}   > \gamma_{\bar{a}}$ such that
\begin{align*}
	\hat{M} (5^{\alpha} - 3^{\alpha} + 1) \gamma_{\bar{a}}^{\alpha} \le M \gamma_{\bar{a}}^{\alpha} \le M \| \x - \c_{\bkappa_{\bar{a}},\bar{a}} \|_{\infty}^{\alpha},
\end{align*}
where the first inequality holds by \eqref{ineq:lb_M33}, the upper bound in \eqref{ineq:scv_lb_bd} holds.

In general, to analyze $y_{\CP,a} (\x|\bkappa_{\bar{a}},\bar{a}) = y_{\mathcal{D}} (\x|\bkappa_{\bar{a}},\bar{a}) + \sum_{\ell=1}^{a} y_{\CS} \big(\x|\Pcal^{\ell}(\bkappa_{\bar{a}}),\bar{a}-\ell\big)$, $a = 2,3,\dots,\bar{a}$, assume
\begin{align}\label{ineq:scv_lb_bdg}
	\tilde{M} \| \x - \c_{\bkappa_{\bar{a}},\bar{a}} \|_{\infty}^\alpha \le \big| y_{\CP,a-1} (\c_{\bkappa_{\bar{a}},\bar{a}} | \bkappa_{\bar{a}},\bar{a}) - y_{\CP,a-1} (\x | \bkappa_{\bar{a}},\bar{a}) \big| \le M \| \x - \c_{\bkappa_{\bar{a}},\bar{a}} \|_{\infty}^\alpha
\end{align}
for $\big\| \x - \c_{\Pcal^{a-1}(\bkappa_{\bar{a}}),\bar{a}-a+1} \big\|_{\infty}    \le \gamma_{\bar{a}-a+1}$.
For any $\ell = 2,3,\dots,a$, since $\Zcal_{\Pcal^{\ell}(\bkappa_{\bar{a}}),\bar{a}-\ell} \in \Ucal_{\bar{a}-\ell}$, we have 
\begin{align}\label{eq:scv_lb_cdist3}
	\big\| \c_{\Pcal^{\ell-1}(\bkappa_{\bar{a}}),\bar{a}-\ell+1} - \c_{\Pcal^{\ell}(\bkappa_{\bar{a}}),\bar{a}-\ell} \big\|_{\infty} \le 2 \gamma_{\bar{a}-\ell+1}.
\end{align}
We show by discussing two cases that if $\big\| \x - \c_{\Pcal^{a}(\bkappa_{\bar{a}}),\bar{a}-a} \big\|_{\infty}    \le 5\gamma_{\bar{a}-a+1}$, then
\begin{align}\label{ineq:scv_lb_bd3}
	\tilde{M} \| \x - \c_{\bkappa_{\bar{a}},\bar{a}} \|_{\infty}^\alpha \le \big| y_{\CP,a} (\c_{\bkappa_{\bar{a}},\bar{a}}|\bkappa_{\bar{a}},\bar{a}) - y_{\CP,a} (\x|\bkappa_{\bar{a}},\bar{a}) \big| \le M \| \x - \c_{\bkappa_{\bar{a}},\bar{a}} \|_{\infty}^\alpha.
\end{align}
\textbf{Case 1.} Suppose $\big\| \x - \c_{\Pcal^{a-1}(\bkappa_{\bar{a}}),\bar{a}-a+1} \big\|_{\infty}   \le \gamma_{\bar{a}-a+1}$, which together with \eqref{eq:scv_lb_cdist3} implies $ \big\| \x - \c_{\Pcal^{a}(\bkappa_{\bar{a}}),\bar{a}-a} \big\|_{\infty}   \le 3\gamma_{\bar{a}-a+1}$. Thus $y_{\CS} \big(\x|\Pcal^a(\bkappa_{\bar{a}}),\bar{a}-a\big) = \hat{M} (5^{\alpha}-3^{\alpha}) \gamma_{\bar{a}-a+1}^{\alpha}$ by \eqref{eq:scv_fun2} such that 
\begin{align*}
	y_{\CP,a} (\c_{\bkappa_{\bar{a}},\bar{a}}|\bkappa_{\bar{a}},\bar{a}) &= y_{\CP,a-1} (\c_{\bkappa_{\bar{a}},\bar{a}} | \bkappa_{\bar{a}},\bar{a}) + \hat{M} (5^{\alpha}-3^{\alpha}) \gamma_{\bar{a}-a+1}^{\alpha}, \\
	y_{\CP,a} (\x|\bkappa_{\bar{a}},\bar{a}) &= y_{\CP,a-1} (\x|\bkappa_{\bar{a}},\bar{a}) + \hat{M} (5^{\alpha}-3^{\alpha}) \gamma_{\bar{a}-a+1}^{\alpha}.
\end{align*}
Then \eqref{ineq:scv_lb_bd3} holds by \eqref{ineq:scv_lb_bdg}.

\textbf{Case 2.} Suppose $\big\| \x - \c_{\Pcal^{a-1}(\bkappa_{\bar{a}}),\bar{a}-a+1} \big\|_{\infty}   > \gamma_{\bar{a}-a+1}$ and $\big\| \x - \c_{\Pcal^{a}(\bkappa_{\bar{a}}),\bar{a}-a} \big\|_{\infty}    \le 5\gamma_{\bar{a}-a+1}$. Combining the inequality $\big\| \x - \c_{\Pcal^{a-1}(\bkappa_{\bar{a}}),\bar{a}-a+1} \big\|_{\infty}   > \gamma_{\bar{a}-a+1}$ with  \eqref{eq:scv_lb_cdist3}, we have
\begin{align*}
	\| \x - \c_{\bkappa_{\bar{a}},\bar{a}} \|_{\infty}   
	\ge& \big\| \x - \c_{\Pcal^{a-1}(\bkappa_{\bar{a}}),\bar{a}-a+1} \big\|_{\infty} - \sum_{\ell=2}^{a-1} \big\| \c_{\Pcal^{\ell-1}(\bkappa_{\bar{a}}),\bar{a}-\ell+1} - \c_{\Pcal^{\ell}(\bkappa_{\bar{a}}),\bar{a}-\ell} \big\|_{\infty} - \big\|  \c_{\bkappa_{\bar{a}},\bar{a}} - \c_{\Pcal(\bkappa_{\bar{a}}),\bar{a}-1} \big\|_{\infty} \\
	\ge&  \gamma_{\bar{a}-a+1} - \sum_{\ell=2}^{a-1} 2 \gamma_{\bar{a}-\ell+1} - 2 \gamma_{\bar{a}} \\
    =& \gamma_{\bar{a}-a+1} - \sum_{\ell=2}^{a-1} \frac{2}{5^{a-\ell}} \gamma_{\bar{a}-a+1} - \frac{2}{5^{a-1}} \gamma_{\bar{a}-a+1} \\
	\ge& \gamma_{\bar{a}-a+1} - \sum_{\ell'=1}^{\infty} \frac{2}{5^{\ell'}} \gamma_{\bar{a}-a+1}  \\ 
	=& \frac{1}{2} \gamma_{\bar{a}-a+1}
\end{align*}
Combining the inequality of $\big\| \x - \c_{\Pcal^{a}(\bkappa_{\bar{a}}),\bar{a}-a} \big\|_{\infty}    \le 5\gamma_{\bar{a}-a+1}$ with \eqref{eq:scv_lb_cdist3}, we have
\begin{align*}
	\| \x - \c_{\bkappa_{\bar{a}},\bar{a}} \|_{\infty} 
	\le& \big\| \x - \c_{\Pcal^{a}(\bkappa_{\bar{a}}),\bar{a}-a} \big\|_{\infty} + \sum_{\ell=2}^{a} \big\| \c_{\Pcal^{\ell-1}(\bkappa_{\bar{a}}),\bar{a}-\ell+1} - \c_{\Pcal^{\ell}(\bkappa_{\bar{a}}),\bar{a}-\ell} \big\|_{\infty} + \big\|   \c_{\Pcal(\bkappa_{\bar{a}}),\bar{a}-1} - \c_{\bkappa_{\bar{a}},\bar{a}} \big\|_{\infty}  \\
	\le& 5\gamma_{\bar{a}-a+1} + \sum_{\ell=2}^{a}\frac{2}{5^{a-\ell}}\gamma_{\bar{a}-a+1} + \frac{2}{5^{a-1}}\gamma_{\bar{a}-a+1} \\ 
    \le& 7\gamma_{\bar{a}-a+1} + \sum_{\ell'=1}^{\infty} \frac{2}{5^{\ell'}} \gamma_{\bar{a}-a+1} \\
	\le& 8\gamma_{\bar{a}-a+1}.
\end{align*}
For $\x$ in this case, we have $y_{\CP,a-1}  (\x|\bkappa_{\bar{a}},\bar{a}) = 0$ such that $y_{\CP,a} (\x|\bkappa_{\bar{a}},\bar{a}) = y_{\CS} \big(\x|\Pcal^a(\bkappa_{\bar{a}}),\bar{a}-a\big)$ and thus
\begin{align*}
	& y_{\CP,a} (\c_{\bkappa_{\bar{a}},\bar{a}}|\bkappa_{\bar{a}},\bar{a}) - y_{\CP,a} (\x|\bkappa_{\bar{a}},\bar{a})   \\
	=& y_{\CP,a-1} (\c_{\bkappa_{\bar{a}},\bar{a}}|\bkappa_{\bar{a}},\bar{a}) + \hat{M} (5^{\alpha}-3^{\alpha}) \gamma_{\bar{a}-a+1}^{\alpha} - y_{\CS} \big(\x|\Pcal^a(\bkappa_{\bar{a}}),\bar{a}-a\big)  \\
	=& \left\{ \begin{array}{ll}
		y_{\CP,a-1} (\c_{\bkappa_{\bar{a}},\bar{a}}|\bkappa_{\bar{a}},\bar{a}) , &  \big\| \x - \c_{\Pcal^a(\bkappa_{\bar{a}}),\bar{a}-a} \big\|_{\infty}   \le 3\gamma_{\bar{a}-a+1} \\
		y_{\CP,a-1} (\c_{\bkappa_{\bar{a}},\bar{a}}|\bkappa_{\bar{a}},\bar{a}) - \hat{M} 3^{\alpha} \gamma_{\bar{a}-a+1}^{\alpha} + \hat{M} \big\| \x - \c_{\Pcal^a(\bkappa_{\bar{a}}),\bar{a}-a} \big\|_{\infty}^\alpha,     &  3\gamma_{\bar{a}-a+1}  < \big\| \x - \c_{\Pcal^a(\bkappa_{\bar{a}}),\bar{a}-a} \big\|_{\infty}  \le 5\gamma_{\bar{a}-a+1} 
	\end{array} \right. 
\end{align*}
Note that $y_{\CP,a-1} (\c_{\bkappa_{\bar{a}},\bar{a}} | \bkappa_{\bar{a}},\bar{a}) = \hat{M} \gamma_{\bar{a}}^{\alpha} +   \sum_{\ell=1}^{a-1} \hat{M} (5^{\alpha}-3^{\alpha}) \gamma_{\bar{a}-\ell+1}^{\alpha}$. Then
\begin{align*}
	\hat{M} \gamma_{\bar{a}}^{\alpha} +   \sum_{\ell=1}^{a-1} \hat{M} (5^{\alpha}-3^{\alpha}) \gamma_{\bar{a}-\ell+1}^{\alpha} \le y_{\CP,a} (\c_{\bkappa_{\bar{a}},\bar{a}}|\bkappa_{\bar{a}},\bar{a}) - y_{\CP,a} (\x|\bkappa_{\bar{a}},\bar{a}) \le \hat{M} \gamma_{\bar{a}}^{\alpha} +   \sum_{\ell=1}^{a} \hat{M} (5^{\alpha}-3^{\alpha}) \gamma_{\bar{a}-\ell+1}^{\alpha}.
\end{align*}
Notice that $\hat{M} = \frac{40^{\alpha}}{5^{\alpha}-3^{\alpha}} \tilde{M}$ such that $8^{\alpha} \tilde{M} = \frac{5^{\alpha}-3^{\alpha}}{5^{\alpha}}\hat{M}$. Since
\begin{align*}
	\hat{M} \gamma_{\bar{a}}^{\alpha} +   \sum_{\ell=1}^{a-1} \hat{M} (5^{\alpha}-3^{\alpha}) \gamma_{\bar{a}-\ell+1}^{\alpha} \ge \hat{M} (5^{\alpha}-3^{\alpha}) \gamma_{\bar{a}-a+2}^{\alpha} \ge \tilde{M} 8^{\alpha} \gamma_{\bar{a}-a+1}^{\alpha} \ge \tilde{M} \| \x - \c_{\bkappa_{\bar{a}},\bar{a}} \|_{\infty}^{\alpha},
\end{align*}
the lower bound of \eqref{ineq:scv_lb_bd3} holds. Meanwhile, since
\begin{align*}
	\hat{M} \gamma_{\bar{a}}^{\alpha} +   \sum_{\ell=1}^{a} \hat{M} (5^{\alpha}-3^{\alpha}) \gamma_{\bar{a}-\ell+1}^{\alpha} \le \frac{\hat{M} (5^{\alpha}-3^{\alpha}+1)}{ 1-5^{-\alpha} } \gamma_{\bar{a}-a+1}^{\alpha}\le M \left(\frac{1}{2} \gamma_{\bar{a}-a+1}\right)^{\alpha} \le M \| \x - \c_{\bkappa_{\bar{a}},\bar{a}} \|_{\infty}^{\alpha},
\end{align*}
where the first inequality holds because
\begin{align*}
	\hat{M} \gamma_{\bar{a}}^{\alpha} + \sum_{\ell=1}^{a} \hat{M} (5^{\alpha}-3^{\alpha}) \gamma_{\bar{a}-\ell+1}^{\alpha} \le& \sum_{\ell=1}^{a} \hat{M} (5^{\alpha}-3^{\alpha}+1) \gamma_{\bar{a}-\ell+1}^{\alpha} 
	\le \frac{\hat{M} (5^{\alpha}-3^{\alpha}+1)}{ 1-5^{-\alpha} } \gamma_{\bar{a}-a+1}^{\alpha} . 
\end{align*}
and the second inequality holds by \eqref{ineq:lb_M33}.
Then the upper bound of \eqref{ineq:scv_lb_bd3} holds.

Notice that $y_{\CP,\bar{a}} (\x|\bkappa_{\bar{a}},\bar{a}) = y_{\Ccal} (\x|\bkappa_{\bar{a}},\bar{a})$ by definition. This lemma is proved.

\subsection{Proof of Proposition \ref{prop:lb_scv2}}\label{subsec:lb_scv_fun2}

We need the following corollary adapted from Lemma \ref{lem:lb3} during the proof of Proposition \ref{prop:lb_scv2}.

\begin{corollary}\label{coro:lem_lb3_cond}
	Let $\tilde{\Zcal}$ denote a subregion of $\Xcal$. Consider $3^d$ mutually disjoint subregions $\Zcal_{\kappa}$, $\kappa = 1,\dots,3^d$, in subregion $\tilde{\Zcal}$. When an algorithm satisfying Assumption \ref{ass:adaptive_algo} is applied to any given CSO problem and $ \Pb\big( \tau(\tilde{\Zcal}\big) \le n - \Delta ) > 0$ where $\Delta=\lfloor 3^d/2 \rfloor$, there must exist a subregion $\Zcal_{\kappa_1}$ such that 
	\begin{align}
		&\Pb\Big( \tau(\Zcal_{\kappa_1}) \ge \tau(\tilde{\Zcal}) + \Delta \Big| \tau(\tilde{\Zcal}) \le n - \Delta\Big) \ge \frac{1}{2}, \label{ineq:lb_3_a1} 
	\end{align}
	On the other hand, if $\Pb\big( \tau(\tilde{\Zcal}) \le n - \Delta \big) < 1$, then there must exist a subregion $\Zcal_{\kappa_2}$ such that 
	\begin{align}
		&\Pb\Big( \tau(\Zcal_{\kappa_2}) > n \Big| \tau(\tilde{\Zcal}) > n - \Delta \Big) \ge \frac{1}{2}. \label{ineq:lb_3_a2}
	\end{align}
\end{corollary}

Corollary \ref{coro:lem_lb3_cond} is similar to Lemma \ref{lem:lb3}, except that Corollary \ref{coro:lem_lb3_cond} is based on the conditional probability and discusses the smaller subregions in a larger subregion $\tilde{\Zcal}$. The proof of Corollary \ref{coro:lem_lb3_cond} will be provided in Section \ref{sec:pf_coro_cond}.

We have defined a problem $\CSO_{\bkappa_{a},a}$ for each $\bkappa_{a} \in \Ucal_{a}$, $a=0,1,\dots,\bar{a}_{2}$, whose objective function is
$ y_{\Ccal}(\x|\bkappa_{a},a)$. Let $\Pb_{\bkappa_a,a}$ and $\Eb_{\bkappa_a,a}$ denote the probability and expectation when the CSO problem is $\CSO_{\bkappa_a,a}$. Let $\Hcal_{\mathbf{1},0} = \Ucal_1$ denote the set of subregions in $\Ucal_{1}$ that are children of $\Zcal_{\mathbf{1},0}$. For the subregion $\Zcal_{\bkappa_a,a} \in \Ucal_{a}$ at level $a=1,\dots,\bar{a}-1$, let
\begin{align*}
	\Hcal_{\bkappa_a,a} = \big\{ \Zcal_{\bkappa_{a+1},a+1} \in \Ucal_{a+1}: \Pcal(\bkappa_{a+1}) = \bkappa_a \big\} = \big\{ \Zcal_{\bkappa_{a+1},a+1} \in \Ucal_{a+1}: \| \c_{\bkappa_{a+1},a+1} - \c_{\bkappa_a,a} \|_{\infty} \le 2 \gamma_{a+1}  \big\},
\end{align*}
which contains all subregions in $\Ucal_{a+1}$ that are children of $\Zcal_{\bkappa_a,a}$.  Obviously, $\Zcal_{\bkappa_{a+1},a+1} \subseteq \Zcal_{\bkappa_{a},a}$ if $\Zcal_{\bkappa_{a+1},a+1} \in \Hcal_{\bkappa_a,a}$. By definition, $|\Hcal_{\bkappa_a,a}| = 3^d \ge 2 \Delta$.

Suppose $n \ge b_{\Ccal,3}$ where we let $b_{\Ccal,3} = \Delta+1$. By Lemma \ref{lem:lb3}, there must exist $\Zcal_{\bar{\bkappa}_1,1} \in \Hcal_{\mathbf{1},0}$ such that $\Pb_{\mathbf{1},0} \left( \Ecal_1 \right) \ge \frac{1}{2}$ where we define $\Ecal_1 = \big\{ \tau(\Zcal_{\bar{\bkappa}_1,1}) > \Delta \big\}$. Then
\begin{align*}
	\Eb_{\mathbf{1},0} \left[ \tau(\Zcal_{\bar{\bkappa}_1,1}) \right] 
	\ge& \Eb_{\mathbf{1},0} \left[ \tau(\Zcal_{\bar{\bkappa}_1,1}) | \Ecal_1 \right] \Pb_{\mathbf{1},0} \left(  \Ecal_1 \right) \ge \frac{1}{2} \Delta.
\end{align*}
Since $y_{\bar{\bkappa}_1,1} (\x) \ne y_{\mathbf{1},0} (\x)$ only when $\x \in \Zcal_{\bar{\bkappa}_1,1}$, by Lemma \ref{lem:lb2}, $\Eb_{\bar{\bkappa}_1,1} \left[ \tau(\Zcal_{\bar{\bkappa}_1,1}) \right] = \Eb_{\mathbf{1},0} \left[ \tau(\Zcal_{\bar{\bkappa}_1,1}) \right] \ge \Delta/2$.

Let $\tilde{\Ecal}_1 = \big\{ \tau(\Zcal_{\bar{\bkappa}_1,1}) \le n-\Delta \big\}$ denotes the ``regular" event at level $a=1$. If $\Pb_{\bar{\bkappa}_1,1} (\tilde{\Ecal}_1) > 0$, by \eqref{ineq:lb_3_a1}, there must exist $\Zcal_{\bar{\bkappa}_2,2} \in \Hcal_{\bar{\bkappa}_1,1}$ such that  $\Pb_{\bar{\bkappa}_1,1}  \big(\Ecal_2 \big| \tilde{\Ecal}_1 \big) \ge \frac{1}{2}$ where we define $\Ecal_2 = \big\{ \tau(\Zcal_{\bar{\bkappa}_2,2}) \ge \tau(\Zcal_{\bar{\bkappa}_1,1})+\Delta \big\}$.  Since  $\tau(\Zcal_{\bar{\bkappa}_2,2}) \ge \tau(\Zcal_{\bar{\bkappa}_1,1})+\Delta$ under $\tilde{\Ecal}_1 \cap \Ecal_2$, we have
\begin{align*}
	\Eb_{\bar{\bkappa}_1,1} \left[ \tau(\Zcal_{\bar{\bkappa}_2,2}) \1 (\tilde{\Ecal}_1 \cap \Ecal_2) \right] \ge& \Eb_{\bar{\bkappa}_1,1} \left[ \big(\tau(\Zcal_{\bar{\bkappa}_1,1})+\Delta\big) \1 (\tilde{\Ecal}_1 \cap \Ecal_2) \right] \\
	=& \Eb_{\bar{\bkappa}_1,1} \left[ \tau(\Zcal_{\bar{\bkappa}_1,1}) \1 (\tilde{\Ecal}_1 \cap \Ecal_2) \right] + \Delta \Pb_{\bar{\bkappa}_1,1}  \big(\Ecal_2 \big| \tilde{\Ecal}_1\big)  \Pb_{\bar{\bkappa}_1,1}  (\tilde{\Ecal}_1)  \\
	\ge& \Eb_{\bar{\bkappa}_1,1} \left[ \tau(\Zcal_{\bar{\bkappa}_1,1}) \1 (\tilde{\Ecal}_1 \cap \Ecal_2) \right] + \frac{\Delta}{2} \Pb_{\bar{\bkappa}_1,1}  (\tilde{\Ecal}_1) .
\end{align*}
Meanwhile, $\Zcal_{\bar{\bkappa}_2,2} \subset \Zcal_{\bar{\bkappa}_1,1}$, which implies $\tau(\Zcal_{\bar{\bkappa}_2,2}) \ge \tau(\Zcal_{\bar{\bkappa}_1,1})$. Then
\begin{align*}
	\Eb_{\bar{\bkappa}_1,1} \left[ \tau(\Zcal_{\bar{\bkappa}_2,2}) \right] =& \Eb_{\bar{\bkappa}_1,1} \left[ \tau(\Zcal_{\bar{\bkappa}_2,2}) \1 (\tilde{\Ecal}_1^c \cup \Ecal_2^c) \right] +  \Eb_{\bar{\bkappa}_1,1} \left[ \tau(\Zcal_{\bar{\bkappa}_2,2}) \1 (\tilde{\Ecal}_1 \cap \Ecal_2) \right]  \\ 
	\ge& \Eb_{\bar{\bkappa}_1,1} \left[ \tau(\Zcal_{\bar{\bkappa}_1,1}) \1 (\tilde{\Ecal}_1^c \cup \Ecal_2^c) \right] + \Eb_{\bar{\bkappa}_1,1} \left[ \tau(\Zcal_{\bar{\bkappa}_1,1}) \1 (\tilde{\Ecal}_1 \cap \Ecal_2) \right] + \frac{\Delta}{2} \Pb_{\bar{\bkappa}_1,1}  (\tilde{\Ecal}_1)    \\
	\ge& \frac{\Delta}{2} \big(1+\Pb_{\bar{\bkappa}_1,1}  (\tilde{\Ecal}_1) \big).
\end{align*}
Since $y_{\bar{\bkappa}_2,2} (\x) \ne y_{\bar{\bkappa}_1,1} (\x)$ only when $\x \in \Zcal_{\bar{\bkappa}_2,2} \subset \Zcal_{\bar{\bkappa}_1,1}$, we have by Lemma \ref{lem:lb2} that
\begin{align*}
\Eb_{\bar{\bkappa}_2,2} \left[ \tau(\Zcal_{\bar{\bkappa}_2,2}) \right] = \Eb_{\bar{\bkappa}_1,1} \left[ \tau(\Zcal_{\bar{\bkappa}_2,2}) \right] \ge \frac{\Delta}{2} \big(1+\Pb_{\bar{\bkappa}_1,1}  (\tilde{\Ecal}_1) \big) = \frac{\Delta}{2} \big(1+\Pb_{\bar{\bkappa}_2,2} (\tilde{\Ecal}_1)\big) .
\end{align*}  
Moreover, $\Pb_{\bar{\bkappa}_2,2} (\tilde{\Ecal}_1) = \Pb_{\bar{\bkappa}_1,1} (\tilde{\Ecal}_1) > 0$.

On the other hand, if $\Pb_{\bar{\bkappa}_1,1} (\tilde{\Ecal}_1) = 0$, we can select any subregion $\Zcal_{\bar{\bkappa}_2,2} \in \Hcal_{\bar{\bkappa}_1,1}$ and have $\Pb_{\bar{\bkappa}_2,2} (\tilde{\Ecal}_2) \le \Pb_{\bar{\bkappa}_2,2} (\tilde{\Ecal}_1) = \Pb_{\bar{\bkappa}_1,1} (\tilde{\Ecal}_1) = 0$ where we define $\tilde{\Ecal}_2 = \big\{ \tau(\Zcal_{\bar{\bkappa}_2,2}) \le n-\Delta \big\}$.

In general, assume $\Zcal_{\bar{\bkappa}_{a'},a'}$, $a'=1,\dots,a$, exist where $a \le \bar{a}-2$ such that $\Pb_{\bar{\bkappa}_a,a} ( \tilde{\Ecal}_{a'} )>0$, $a'=1,\dots,a-1$, and
\begin{align}\label{ineq:lb3_induct}
	\Eb_{\bar{\bkappa}_a,a} \left[ \tau(\Zcal_{\bar{\bkappa}_a,a}) \right] \ge \frac{\Delta}{2} \Big(1+\sum\nolimits_{a'=1}^{a-1}\Pb_{\bar{\bkappa}_a,a} ( \tilde{\Ecal}_{a'} )\Big)
\end{align}
where we let $\tilde{\Ecal}_{a'} = \big\{ \tau(\Zcal_{\bar{\bkappa}_{a'},a'}) \le n-\Delta \big\}$ denote the regular event at level $a'=1,\dots,a$. If $\Pb_{\bar{\bkappa}_a,a}(\tilde{\Ecal}_{a}) > 0$, by \eqref{ineq:lb_3_a1}, there must exist $\Zcal_{\bar{\bkappa}_{a+1},a+1} \in \Hcal_{\bar{\bkappa}_a,a}$ such that $\Pb_{\bar{\bkappa}_a,a} \big(\Ecal_{a+1} \big| \tilde{\Ecal}_a\big)  \ge \frac{1}{2}$ where we let $\Ecal_{a+1} = \big\{ \tau(\Zcal_{\bar{\bkappa}_{a+1},a+1}) \ge \tau(\Zcal_{\bar{\bkappa}_a,a})+\Delta  \big\}$.  Since  $\tau(\Zcal_{\bar{\bkappa}_{a+1},a+1}) \ge \tau(\Zcal_{\bar{\bkappa}_a,a})+\Delta$ under $\tilde{\Ecal}_a \cap \Ecal_{a+1}$, we have
\begin{align*}
	\Eb_{\bar{\bkappa}_a,a} \left[ \tau(\Zcal_{\bar{\bkappa}_{a+1},a+1}) \1 (\tilde{\Ecal}_a \cap \Ecal_{a+1}) \right] \ge& \Eb_{\bar{\bkappa}_a,a} \left[ \big(\tau(\Zcal_{\bar{\bkappa}_a,a})+\Delta\big) \1 (\tilde{\Ecal}_a \cap \Ecal_{a+1}) \right] \\
	=& \Eb_{\bar{\bkappa}_a,a} \left[ \tau(\Zcal_{\bar{\bkappa}_a,a}) \1 (\tilde{\Ecal}_a \cap \Ecal_{a+1}) \right] + \Delta \Pb_{\bar{\bkappa}_a,a}  (\Ecal_{a+1} \big| \tilde{\Ecal}_a)  \Pb_{\bar{\bkappa}_a,a}  (\tilde{\Ecal}_a)  \\
	\ge& \Eb_{\bar{\bkappa}_a,a} \left[ \tau(\Zcal_{\bar{\bkappa}_a,a}) \1 (\tilde{\Ecal}_a \cap \Ecal_{a+1}) \right] + \frac{\Delta}{2} \Pb_{\bar{\bkappa}_a,a}  (\tilde{\Ecal}_a) 
\end{align*}
Meanwhile, $\Zcal_{\bar{\bkappa}_{a+1},a+1} \subset \Zcal_{\bar{\bkappa}_a,a}$, which implies $\tau(\Zcal_{\bar{\bkappa}_{a+1},a+1}) \ge \tau(\Zcal_{\bar{\bkappa}_a,a})$. Then
\begin{align*}
	\Eb_{\bar{\bkappa}_a,a} \left[ \tau(\Zcal_{\bar{\bkappa}_{a+1},a+1})  \right] 
	\ge& \Eb_{\bar{\bkappa}_a,a} \left[ \tau(\Zcal_{\bar{\bkappa}_{a},a}) \1 (\tilde{\Ecal}_a^c \cup \Ecal_{a+1}^c) \right] + \Eb_{\bar{\bkappa}_a,a} \left[ \tau(\Zcal_{\bar{\bkappa}_a,a}) \1 (\tilde{\Ecal}_a \cap \Ecal_{a+1}) \right] + \frac{\Delta}{2} \Pb_{\bar{\bkappa}_a,a}  (\tilde{\Ecal}_a)    \\
	=& \Eb_{\bar{\bkappa}_a,a} \left[ \tau(\Zcal_{\bar{\bkappa}_a,a})  \right] + \frac{\Delta}{2} \Pb_{\bar{\bkappa}_a,a}  (\tilde{\Ecal}_a) \\
    \ge& \frac{\Delta}{2} \Big(1+\sum\nolimits_{a'=1}^{a}\Pb_{\bar{\bkappa}_a,a}  (\tilde{\Ecal}_{a'}) \Big)
\end{align*}
where the last inequality holds by assumption. Since $y_{\bar{\bkappa}_{a+1},a+1} (\x) \ne y_{\bar{\bkappa}_{a},a} (\x)$ only when $\x \in \Zcal_{\bar{\bkappa}_{a+1},a+1}$, we have by Lemma \ref{lem:lb2} that
\begin{align*}
	\Eb_{\bar{\bkappa}_{a+1},a+1} \left[ \tau(\Zcal_{\bar{\bkappa}_{a+1},a+1}) \right] =& \Eb_{\bar{\bkappa}_a,a} \left[ \tau(\Zcal_{\bar{\bkappa}_{a+1},a+1})  \right] \\
	\ge& \frac{\Delta}{2} \Big(1+\sum\nolimits_{a'=1}^{a}\Pb_{\bar{\bkappa}_a,a} (\tilde{\Ecal}_{a'}) \Big) \\
	=& \frac{\Delta}{2} \Big(1+\sum\nolimits_{a'=1}^{a}\Pb_{\bar{\bkappa}_{a+1},a+1} ( \tilde{\Ecal}_{a'} )\Big).
\end{align*}

On the other hand, if $\Pb_{\bar{\bkappa}_a,a} (\tilde{\Ecal}_a) = 0$, we can select any subregion $\Zcal_{\bar{\bkappa}_{a+1},a+1} \in \Hcal_{\bar{\bkappa}_a,a}$ and have $\Pb_{\bar{\bkappa}_{a+1},a+1} (\tilde{\Ecal}_{a+1}) \le \Pb_{\bar{\bkappa}_{a+1},a+1} (\tilde{\Ecal}_{a}) = \Pb_{\bar{\bkappa}_a,a} (\tilde{\Ecal}_{a}) = 0$ where we let $\tilde{\Ecal}_{a+1} = \big\{ \tau(\Zcal_{\bar{\bkappa}_{a+1},a+1}) \le n-\Delta \big\}$.

Let $e = \bar{a} - 1$. We have obtained subregions $\Zcal_{\bar{\bkappa}_1,1},\Zcal_{\bar{\bkappa}_2,2},\dots,\Zcal_{\bar{\bkappa}_{e},e}$. We discuss by cases.
\begin{enumerate}
	\item If $\Pb_{\bar{\bkappa}_{e},e} ( \tilde{\Ecal}_{a} )>0$, $a=1,\dots,e-1$, \eqref{ineq:lb3_induct} leads to
	\begin{align*}
		\Eb_{\bar{\bkappa}_e,e} \left[ \tau(\Zcal_{\bar{\bkappa}_e,e}) \right] \ge \frac{\Delta}{2} \Big(1+\sum\nolimits_{a=1}^{e-1}\Pb_{\bar{\bkappa}_e,e} ( \tilde{\Ecal}_{a} )\Big).
	\end{align*}
	Since $\Zcal_{\bar{\bkappa}_{e-1},e-1} \subset \Zcal_{\bar{\bkappa}_{a},a}$ for $a \le e-1$, we have
	\begin{align}\label{ineq:lb_scv_tau_nest}
		\Pb_{\bar{\bkappa}_{e},e}  (\tilde{\Ecal}_{e-1})  = \Pb_{\bar{\bkappa}_{e},e} \big( \tau(\Zcal_{\bar{\bkappa}_{e-1},e-1}) \le n-\Delta \big) \le \Pb_{\bar{\bkappa}_{e},e} \big( \tau(\Zcal_{\bar{\bkappa}_{a},a}) \le n-\Delta \big) = \Pb_{\bar{\bkappa}_{e},e}  (\tilde{\Ecal}_{a}). 
	\end{align}
	We must have $\Pb_{\bar{\bkappa}_{e},e}  (\tilde{\Ecal}_{e-1})  \le \frac{1}{2}$ because otherwise, if $\Pb_{\bar{\bkappa}_{e},e}  (\tilde{\Ecal}_{e-1})  > \frac{1}{2}$, then $\Pb_{\bar{\bkappa}_{e},e}  (\tilde{\Ecal}_{a})  > \frac{1}{2}$ for all $a=1,\dots,e-1$ such that by \eqref{ineq:lb3_induct},
	\begin{align*}
		\Eb_{\bar{\bkappa}_{e},e} \left[ \tau(\Zcal_{\bar{\bkappa}_{e},e}) \right] \ge \frac{\Delta}{2} \Big(1+\sum\nolimits_{a=1}^{e-1}\Pb_{\bar{\bkappa}_{e},e}  (\tilde{\Ecal}_{a}) \Big) \ge \frac{e\Delta}{4} = \frac{(\lceil 4(n+1)/\Delta \rceil+1) \Delta }{4} > n+1
	\end{align*}
	contradicting the definition of $\tau(\Zcal_{\bar{\bkappa}_{e},e})$, which should be less than $n+1$. 
	
	\item If $\Pb_{\bar{\bkappa}_{e},e} ( \tilde{\Ecal}_{a} )=0$ for some $a=1,\dots,e-1$, we have $\Pb_{\bar{\bkappa}_{e},e} ( \tilde{\Ecal}_{e-1} )=0$ by \eqref{ineq:lb_scv_tau_nest}.
\end{enumerate}
No matter which case happens, $\Pb_{\bar{\bkappa}_{e},e} ( \tilde{\Ecal}_{e-1} ) \le \frac{1}{2}$, which implies $\Pb_{\bar{\bkappa}_{e},e} ( \tilde{\Ecal}_{e} ) \le \frac{1}{2}$ where $\tilde{\Ecal}_{e} = \big\{ \tau(\Zcal_{\bar{\bkappa}_{e},e}) \le n-\Delta \big\}$. Thus, $\Pb_{\bar{\bkappa}_{e},e} ( \tilde{\Ecal}_{e}^c ) \ge \frac{1}{2}$. 

By \eqref{ineq:lb_3_a2}, there must exist $\Zcal_{\bar{\bkappa}_{\bar{a}},\bar{a}} \in \Hcal_{\bar{\bkappa}_{e},e}$ such that the probability conditional on $\tilde{\Ecal}_{e}^c$ satisfies  $\Pb_{\bar{\bkappa}_{e},e} \big(\tau(\Zcal_{\bar{\bkappa}_{\bar{a}},\bar{a}}) > n \big| \tilde{\Ecal}_{e}^c\big)  \ge \frac{1}{2} $. Then
\begin{align*}
	\Pb_{\bar{\bkappa}_{e},e}\big( \tau(\Zcal_{\bar{\bkappa}_{\bar{a}},\bar{a}}) > n \big) \ge  \Pb_{\bar{\bkappa}_{e},e}\Big( \big\{ \tau(\Zcal_{\bar{\bkappa}_{\bar{a}},\bar{a}}) > n \big\} \cap \tilde{\Ecal}_{\tilde{a}}^c \Big) =  \Pb_{\bar{\bkappa}_{e},e}\Big( \tau(\Zcal_{\bar{\bkappa}_{\bar{a}},\bar{a}}) > n \big| \tilde{\Ecal}_{e}^c \Big) \Pb_{\bar{\bkappa}_{e},e}( \tilde{\Ecal}_{e}^c ) \ge \frac{1}{4}.
\end{align*}
By Lemma \ref{lem:lb2}, $\Pb_{\bar{\bkappa}_{\bar{a}},\bar{a}}\big( \tau(\Zcal_{\bar{\bkappa}_{\bar{a}},\bar{a}}) > n \big) = \Pb_{\bar{\bkappa}_{e},e}\big( \tau(\Zcal_{\bar{\bkappa}_{\bar{a}},\bar{a}}) > n \big) \ge 1/4$. 

Event $\big\{\tau(\Zcal_{\bar{\bkappa}_{\bar{a}},\bar{a}}) > n\big\}$ implies that there is no simulated solution in $\Zcal_{\bar{\bkappa}_{\bar{a}},\bar{a}}$. By Assumption \ref{ass:adaptive_algo}, we have $\hat{\x}^*_n \notin \Zcal_{\bar{\bkappa}_{\bar{a}},\bar{a}}$. Notice that $e+1 = \bar{a}$ and $\gamma_{\bar{a}} = \frac{1}{2} 5^{-\bar{a}} \ge \frac{1}{2} 5^{-4(n+1)/\Delta-3}$. If $\hat{\x}^*_n \notin \Zcal_{\bar{\bkappa}_{\bar{a}},\bar{a}}$, we have 
\begin{align*}
	\Big|\max_{\x \in \Xcal} y_{\Ccal}(\x|\bar{\bkappa}_{\bar{a}},\bar{a}) - y_{\Ccal}(\hat{\x}_n^*|\bar{\bkappa}_{\bar{a}},\bar{a}) \Big| \ge \hat{M} \gamma_{\bar{a}}^{\alpha}  \ge  \frac{\hat{M}5^{-3\alpha}}{2^{\alpha}} 5^{-4(n+1)\alpha/\Delta}  \ge \frac{\hat{M}5^{-3\alpha}}{2^{\alpha}} 5^{-24(n+1)\alpha/3^d}
\end{align*}
where the last inequality holds because $\Delta = \lfloor 3^d/2 \rfloor \ge 3^{d}/6$. Thus,
\begin{align*}
	&\Eb_{\bar{\bkappa}_{\bar{a}},\bar{a}} \left[ \max\nolimits_{\x} y_{\Ccal}(\x|\bar{\bkappa}_{\bar{a}},\bar{a}) - y_{\Ccal}(\hat{\x}^*_n|\bar{\bkappa}_{\bar{a}},\bar{a}) \right]  \\
	\ge& \Pb_{\bar{\bkappa}_{\bar{a}},\bar{a}}\big( \hat{\x}^*_n \notin \Zcal_{\bar{\bkappa}_{\bar{a}},\bar{a}} \big) \Eb_{\bar{\bkappa}_{\bar{a}},\bar{a}} \left[ \max\nolimits_{\x} y_{\Ccal}(\x|\bar{\bkappa}_{\bar{a}},\bar{a}) - y_{\Ccal}(\hat{\x}^*_n|\bar{\bkappa}_{\bar{a}},\bar{a}) | \hat{\x}^*_n \notin \Zcal_{\bar{\bkappa}_{\bar{a}},\bar{a}}\right]  \\
    \ge& \frac{1}{4} \frac{\hat{M}5^{-3\alpha}}{2^{\alpha}} 5^{-24(n+1)\alpha/3^d} \\ 
    =& \frac{\hat{M}}{2^{\alpha+2} 5^{(3+24/3^d)\alpha}} 5^{-24n\alpha/3^d}.
\end{align*}

\subsection{Proof of Corollary \ref{coro:lem_lb3_cond}} \label{sec:pf_coro_cond}

The proof for \eqref{ineq:lb_3_a1} almost follows that of Lemma \ref{lem:lb3} except that we should replace the unconditional probability of Lemma \ref{lem:lb3} with the conditional probability. Specifically, let
\begin{align*}
	\kappa_1 \in \mathop{\arg\min}_{\kappa = 1,\dots,3^d}  \Eb \Big[ \sum\nolimits_{t=\tau(\tilde{\Zcal})}^{\tau(\tilde{\Zcal})+\Delta-1} \1 [ \x_t \in \Zcal_{\kappa} ] \Big| \tau(\tilde{\Zcal}) \le n - \Delta \Big].
\end{align*}
We must have $ \Eb \left[ \sum_{t=\tau(\tilde{\Zcal})}^{\tau(\tilde{\Zcal})+\Delta-1} \1 (\x_t \in \Zcal_{\kappa_1})  \Big| \tau(\tilde{\Zcal}) \le n - \Delta \right] \le \frac{1}{2} $. Otherwise, 
\begin{align}\label{ineq:lb3_ass_con}
	\sum_{\kappa = 1}^{3^d} \Eb \left[ \sum\nolimits_{t=\tau(\tilde{\Zcal})}^{\tau(\tilde{\Zcal})+\Delta-1} \1 ( \x_t \in \Zcal_{\kappa} ) \Big| \tau(\tilde{\Zcal}) \le n - \Delta \right] >  n_1/2.
\end{align}
Meanwhile, since $\Zcal_{\kappa}$, $\kappa = 1,\dots,3^d$, are mutually disjoint, we have $\sum_{\kappa = 1}^{3^d} \1 ( \x_t \in \Zcal_{\kappa} ) \le 1$ for any $t=\tau(\tilde{\Zcal}),\dots,\tau(\tilde{\Zcal})+\Delta-1$, which leads to $\sum_{\kappa = 1}^{3^d} \sum_{t=\tau(\tilde{\Zcal})}^{\tau(\tilde{\Zcal})+\Delta-1} \1 ( \x_t \in \Zcal_{\kappa} ) \le \Delta \le 3^d/2$. This is contradictory to \eqref{ineq:lb3_ass_con}.
We can further show 
\begin{align*}
	\Pb \left( \bigcup\nolimits_{t=\tau(\tilde{\Zcal})}^{\tau(\tilde{\Zcal})+\Delta-1} \{ \x_t \in \Zcal_{\kappa_1} \} \big| \tau(\tilde{\Zcal}) \le n - \Delta \right) \le \frac{1}{2}.
\end{align*}
Thus
$\Pb \big( \tau (\Zcal_{\kappa_1}) < \tau(\tilde{\Zcal})+\Delta \big| \tau(\tilde{\Zcal}) \le n - \Delta \big)  \le \frac{1}{2}$.

Next we show \eqref{ineq:lb_3_a2}. Notice that $\sum_{\kappa = 1}^{3^d} \sum_{t=n-\Delta+1}^{n} \1 ( \x_{t} \in \Zcal_{\kappa} ) \le \Delta$, 
which implies that there must exist $\kappa_2$ with $\Eb [ \sum_{t=n-\Delta+1}^{n} $ $\1 ( \x_{t} \in \Zcal_{\kappa_2} ) \big| \tau(\tilde{\Zcal}) > n - \Delta ] \le \frac{1}{2}$ and thus 
\begin{align*}
	\Pb \left( \bigcup\nolimits_{t=n-\Delta+1}^{n} \{ \x_{t} \in \Zcal_{\kappa_2} \} \big| \tau(\tilde{\Zcal}) > n - \Delta \right) \le \frac{1}{2}.
\end{align*}
When $\{\tau(\tilde{\Zcal}) > n - \Delta\}$ happens, we have $ \bigcup_{t=\tau(\tilde{\Zcal})}^{n} \{ \x_{t} \in \Zcal_{\kappa_2} \} \subset \bigcup_{t=n-\Delta+1}^{n} \{ \x_{t} \in \Zcal_{\kappa_2} \}$ (if $\tau(\tilde{\Zcal}) = n+1$, then $\bigcup_{t=\tau(\tilde{\Zcal})}^{n} \{ \x_{t} \in \Zcal_{\kappa_2} \}= \emptyset$), which yields
\begin{align*}
	\Pb \Big( \tau(\Zcal_{\kappa_2}) \le n \Big| \tau(\tilde{\Zcal}) > n - \Delta \Big) = \Pb \left( \bigcup\nolimits_{t=\tau(\tilde{\Zcal})}^{n} \{ \x_{t} \in \Zcal_{\kappa_2} \} \Big| \tau(\tilde{\Zcal}) > n - \Delta \right) \le \frac{1}{2}.
\end{align*}







%

\end{document}